\documentclass[12pt]{article}
\usepackage[utf8]{inputenc}

\title{Comet-shaped quiver varieties, Weyl group actions, and modified Kostka polynomials}
\author{Mathieu Ballandras \\ \it{Université de Paris} \\ \it{Scuola Internazionale Superiore di Studi Avanzati} \\
\it{Instituto de Ciencias Matemáticas} \\ \small{mballandras@imj-prg.fr} }

\usepackage[english]{babel}
\usepackage[T1]{fontenc}

\usepackage[a4paper,top=3cm,bottom=2cm,left=3cm,right=3cm,marginparwidth=1.75cm]{geometry}

\usepackage{mathtools}

\usepackage[vcentermath]{youngtab}
\usepackage{young}

\usepackage[shortlabels]{enumitem}

\usepackage{amsthm}

\usepackage{relsize,exscale}

\usepackage{amsmath}
\usepackage{amssymb}
\usepackage{graphicx}
\usepackage[colorlinks=true, allcolors=blue]{hyperref}
\usepackage{cleveref}
\usepackage{pst-node}
\usepackage{tikz-cd}
\usetikzlibrary{angles, intersections, quotes}

\usepackage{upgreek}

\usepackage{graphicx}

\newcommand{\commentaires}[1]{}

\makeatletter
\newcommand*\bigcdot{\mathpalette\bigcdot@{.5}}
\newcommand*\bigcdot@[2]{\mathbin{\vcenter{\hbox{\scalebox{#2}{$\m@th#1\bullet$}}}}}
\renewenvironment{abstract}{
\null\vfil
\@beginparpenalty\@lowpenalty
\begin{center}
\bfseries \abstractname
\@endparpenalty\@M
\end{center}
}
\makeatother

\newtheorem{definition}{Definition}[section]
\newtheorem{theorem}[definition]{Theorem}
\newtheorem{corollary}[definition]{Corollary}
\newtheorem{lemma}[definition]{Lemma}
\newtheorem{proposition}[definition]{Proposition}

\newtheorem{remark}[definition]{Remark}
\newtheorem{example}[definition]{Example}
\newtheorem{notations}[definition]{Notations}
\newtheorem{conjecture}[definition]{Conjecture}
\newtheorem{assumption}[definition]{Assumption}

\DeclareMathOperator{\G}{G}
\DeclareMathOperator{\Spec}{Spec}
\DeclareMathOperator{\Hom}{Hom}
\DeclareMathOperator{\SymF}{Sym}

\DeclareMathOperator{\Id}{Id}

\DeclareMathOperator{\Ind}{Ind}

\DeclareMathOperator{\Rep}{Rep}
\DeclareMathOperator{\tr}{tr}
\DeclareMathOperator{\Aut}{Aut}
\DeclareMathOperator{\Imm}{Im}

\DeclareMathOperator{\op}{op}
\DeclareMathOperator{\ch}{ch}

\DeclareMathOperator{\GL}{GL}

\DeclareMathOperator{\reg}{reg}
\DeclareMathOperator{\Exp}{Exp}
\DeclareMathOperator{\Log}{Log}
\DeclareMathOperator{\PGL}{PGL}
\DeclareMathOperator{\gl}{\mathfrak{gl}}

\newcommand\gv{\mathfrak{g}_v}
\DeclareMathOperator{\supp}{Supp}

\DeclareMathOperator{\Mat}{Mat}

\DeclareMathOperator{\gen}{gen}

\newcommand{\map}[5]{
\begin{array}{ccccc}
    #1 & : & #2 & \to & #3  \\
     & & #4 & \mapsto & #5
\end{array}
}
\newcommand{\git}{/\!/}
\newcommand\KK{\mathbb K }

\newcommand{\C}{\mathbb{C}}
\newcommand{\qlbar}{\overline{\mathbb{Q}}_l}

\newcommand{\IH}{IH}

\newcommand{\Sym}{\mathfrak{S}}
\newcommand{\Fq}{\mathbb{F}_q}
\newcommand{\Fqbar}{\overline{\mathbb{F}}_q}
\newcommand{\QQ}{\mathbb Q }
\newcommand{\NN}{\mathbb N}

\newcommand{\SymX}{\SymF \left[X\right]}
\newcommand{\Part}{\mathcal{P}}

\newcommand{\Ht}{\tilde{H}}
\newcommand{\SymnX}{\SymF_n\left[X\right]}
\newcommand{\inv}{^{-1}}

\newcommand{\dcb}[1]{\mathcal{D}_c^b\left(#1\right)}
\newcommand{\hi}{\mathcal H ^i}
\newcommand{\uic}[1]{\underline{\mathcal{I C}}_{#1}^\bullet}
\newcommand{\ic}[1]{\mathcal{I C}_{#1}^\bullet}

\newcommand{\ZZ}{\mathbb{Z}}

\newcommand{\VLPs}{\mathcal{V}_{\mathbf{L},\mathbf{P},\sigma} }

\newcommand{\mlps}{\widetilde{\mathcal{M}}_{\boldsymbol{L},\boldsymbol{P},\boldsymbol{\sigma}}}

\newcommand{\qlps}{\widetilde{\mathcal{Q}}_{\boldsymbol{L},\boldsymbol{P},\boldsymbol{\sigma}}}
\newcommand{\qcmus}{\mathcal{Q}_{\overline{\mathcal{O}}_{\boldsymbol{\mu},\boldsymbol{\sigma}}}}
\newcommand{\omus}{\mathcal{O}_{\boldsymbol{\mu},\boldsymbol{\sigma}}}
\newcommand{\qlp}{\widetilde{\mathcal{Q}}_{\boldsymbol{L},\boldsymbol{P}} }
\newcommand{\omusbar}{\overline{\mathcal{O}}_{\boldsymbol{\mu},\boldsymbol{\sigma}}}
\newcommand{\vlp}{\widetilde{\mathcal{V}}_{\boldsymbol{L},\boldsymbol{P}} }

\newcommand{\ylpt}{\widetilde{\mathbb{Y}}_{\boldsymbol{L},\boldsymbol{P}}}
\newcommand{\ylp}{\mathbb{Y}_{\boldsymbol{L},\boldsymbol{P}}}

\newcommand{\qobar}{\mathcal{Q}_{ \boldsymbol{ \overline{ \mathcal{O} } } }}

\newcommand{\Htl}{\tilde{H}_{\lambda}[X;q,t]}

\begin{document}

\maketitle

\begin{abstract}
    We study an algebra spanned by modified Kostka polynomials. Particular structure coefficients of this algebra are interpreted as traces of some Weyl group actions on the intersection cohomology of comet-shaped quiver varieties.
\end{abstract}

\tableofcontents

\section{Introduction}
The modified Kostka polynomials $\widetilde{K}_{\lambda,\rho}(q,t)$ form a family of two-variable polynomials indexed by pairs of partitions of some integer $n$. They are a two-parameter deformation of the Kostka numbers and appear in the theory of symmetric functions. They were introduced by Garsia--Haiman \cite{garsia} in the expression of the modified Macdonald polynomials in terms of Schur functions. The fact that they are polynomials with non-negative integer coefficients is an important result. This is known as the Macdonald conjecture \cite{Macdonald1988} which is a consequence of the $n!$-conjecture of Garsia--Haiman \cite{garsia_haiman_n_fact}, proved by Haiman \cite{haiman_n_fact}. In this article we study an aglebra spanned by those polynomials. Its structure coefficients $c_{\mu,\nu}^{\lambda}(q,t)$ were introduced by Rodriguez-Villegas in unpublished notes, they are defined by
\[
\widetilde{K}_{\mu,\rho}\widetilde{K}_{\nu,\rho}=\sum_{\lambda \in \Part_n} c_{\mu,\nu}^{\lambda} \widetilde{K}_{\lambda,\rho}.
\]
with $\Part_n$ the set of partitions of the integer $n$. We focus in particular on the coefficients $c_{\mu,\nu}^{1^n}(q,t)$. They generalize the $q,t$-Catalan sequence of Garsia--Haiman \cite{garsia}. We give an interpretation of the coefficients $c_{\mu,\nu}^{1^n}(0,t)$ in terms of Weyl group action on the cohomology of a comet-shaped quiver variety. We also give a conjectural  interpretation of the coefficients $c_{\mu,\nu}^{1^n}(q,t)$ in terms of cohomology of a partial resolution of a character variety. We prove the $q=1$ specialization of the conjecture. 

The geometric framework is the following. The base field $\KK$ is either $\C$ or $\Fqbar$ an algebraic closure of a field $\Fq$ with $q$ elements. We consider cohomology with coefficients in $\kappa$ which is either $\C$ or $\qlbar$ with $l$ and $q$ coprime. Fix some integers $n>0$, $g\ge 0$ and $k>0$. Let $\boldsymbol{\mathcal{O}}=\left(\mathcal{O}_1,\dots,\mathcal{O}_k\right)$ be a $k$-tuple of adjoint orbits in $\gl_n$, the Lie algebra of $\GL_n$. Denote by $\overline{\mathcal{O}}_j$ the Zariski closure of the class $\mathcal{O}_j$. Some genericity conditions are imposed on the eigenvalues of the $k$-tuple $\boldsymbol{\mathcal{O}}$ (see Definition \ref{def_generic_adjoint}). The main object in this article is the following variety:
\begin{multline*}
\mathcal{Q}_{ \boldsymbol{ \overline{ \mathcal{O} } } }:=\left\lbrace\left(A_1,B_1,\dots,A_g,B_g,X_1,\dots,X_k\right)\in \gl_n^{2g}\times \overline{\mathcal{O}}_{1}\times\dots\times \overline{\mathcal{O}}_{k}\right| \\ \left. \sum_{i=1}^g \left[A_i,B_i\right] +\sum_{j=1}^k X_j = 0 \right\rbrace \git \GL_n,
\end{multline*}
with $[A_i,B_i]:= A_i B_i-B_i A_i$ the Lie bracket. The quotient is a Geometric Invariant Theory (GIT) quotient with respect to the overall adjoint action of $\GL_n$. Those varieties were studied by Crawley-Boevey \cite{crawley-boevey2003,CB_surj} in genus $g=0$. For any genus and semisimple adjoint orbits, they were studied by Letellier, Hausel and Rodriguez-Villegas \cite{aha}. Letellier \cite{letellier_gl} generalized to any type of adjoint orbits. Those varieties are called comet-shaped quiver varieties due to their interpretation as Nakajima's quiver varieties. They are additive analogues of the following character varieties:
\begin{multline*}
\mathcal{M}_{ \boldsymbol{ \overline{ \mathcal{C} } } }:=\left\lbrace\left(A_1,B_1,\dots,A_g,B_g,X_1,\dots,X_k\right)\in \GL_n^{2g}\times \overline{\mathcal{C}}_{1}\times\dots\times \overline{\mathcal{C}}_{k}\right| \\ \left. A_1 B_1 A_1\inv B_1\inv \dots A_g B_g A_g\inv B_g\inv X_1\dots X_k =\Id \right\rbrace \git \GL_n
\end{multline*}
with $\boldsymbol{\mathcal{C}}=\left(\mathcal{C}_1,\dots,\mathcal{C}_k\right)$ a $k$-tuple of conjugacy classes in $\GL_n$. The action of $\GL_n$ is by overall conjugation. Those varieties classify representations of the fundamental group of a genus $g$ Riemann surface with $k$ punctures and prescribed monodromies around those punctures. They were extensively studied by  Hausel, Letellier and Rodriguez-Villegas \cite{aha} for semisimple classes and by Letellier \cite{zclosure} for any Jordan type.  Similarly to the additive case there is a genericity condition imposed on the eigenvalues of $\boldsymbol{\mathcal{C}}$ in order for the quotient to be well behaved (see Definition \ref{def_generic_cc}). 

In general the varieties $\mathcal{Q}_{ \boldsymbol{ \overline{ \mathcal{O} } } }$ and $\mathcal{M}_{ \boldsymbol{ \overline{ \mathcal{C} } } }$ are singular, therefore it is interesting to study their intersection cohomology. Letellier \cite{letellier_gl,zclosure} studied intersection cohomology of those varieties by constructing resolutions of singularities $\qlps\to \mathcal{Q}_{ \boldsymbol{ \overline{ \mathcal{O} } } }$ (see Definition \ref{def_resol_char}). Those resolutions of singularities originate in Springer theory \cite{Springer1976,borho} and Lusztig's theory of parabolic induction  \cite{Lusztig1984,lusztig_character_sheaves,Lusztig1986}. Thanks to this origin, the resolutions of singularities (in both the additive and  multiplicative cases) come with a Weyl group action on their cohomology. This action is called the Springer action. In this article we focus on the additive version (the comet-shaped quiver varieties). We come back to character varieties, at the end, in order  to give a conjectural interpretation of the coefficient $c_{\mu,\nu}^{1^n}(q,t)$.

In the additive case, thanks to the quiver variety point of view,  Weyl group actions other than the Springer action can be constructed. The Weyl group actions on the cohomology of Nakajima's quiver varieties were studied extensively by Nakajima \cite{nakajima1994,Nakajima_00}, Lusztig \cite{lusztig_quiver} and Maffei \cite{Maffei}. They were used to prove Kac conjecture by Letellier, Hausel, Rodriguez-Villegas \cite{kac_conj} and to study unipotent characters of $\GL_n(\mathbb{F}_q)$ by Letellier \cite{Letellier_uni_12}.
Nakajima's construction \cite{nakajima1994} of the Weyl group action relies on the hyperkähler structure of the quiver varieties. We use this construction, together with technics from Lusztig (see Letellier \cite[Proof of proposition 5.5.3]{letellier_f}) and ideas from Mellit \cite[Section 8]{mellit_cell}. This allows to describe uniformly a monodromic action on the cohomology of quiver varieties with semisimple adjoint orbits and the Springer action on the cohomology of resolutions of quiver varieties.

Thanks to this construction we can obtain an interpretation of the coefficients $c_{\mu,\nu}^{1^n}(0,t)$ in terms of the Weyl group action on the compactly supported intersection cohomology of some comet-shaped quiver varieties.
\begin{theorem}
Consider a generic $4$-tuple of adjoint orbits of the following type:
\begin{itemize}
    \item $\mathcal{O}_1$ has one eigenvalue with Jordan type $\mu'\in \Part_n$,
    \item $\mathcal{O}_2$ has one eigenvalue with Jordan type $\nu'\in \Part_n$,
    \item $\mathcal{O}_3$ is semisimple regular, it has $n$ distinct eigenvalues,
    \item $\mathcal{O}_4$ is semisimple with one eigenvalue of multiplicity $n-1$ and the other of multiplicity $1$.
\end{itemize}
Then the Weyl group with respect to $\mathcal{O}_3$ is the symmetric group $\Sym_n$ and it acts on the compactly supported intersection cohomology of $\qobar$. Let $w$ be a $n$-cycle in this Weyl group, then
\[
c_{\mu,\nu}^{1^n}\left(0,{t}\right)={t}^{\frac{-d_{\boldsymbol{\mathcal{O}}}}{2}}\sum_r \tr\left(w, IH_c^{2 r} \left(\qobar,\kappa\right)\right) {t}^{r}.
\]
\end{theorem}

In general, in the multiplicative case, only the Springer action exists. The monodromic action was constructed by Mellit \cite[Section 8]{mellit_cell} only with respect to one puncture for a regular monodromy. The Springer action is defined by Letellier \cite{zclosure} on the cohomology of partial resolutions of  character varieties. The partial resolution relevant to describe the coefficient $c_{\mu,\nu}^{1^n}$ is the variety $\mathcal{M}_{\mu,\nu}$ introduced in \ref{subsect_conj_mixed_c}. The group $\Sym_n$ acts on the compactly supported intersection cohomology of $\mathcal{M}_{\mu,\nu}$. For $w\in \Sym_n$, the $w$-twisted mixed Hodge polynomial of $\mathcal{M}_{\mu,\nu}$ is

\[
IH_c^{w}\left( \mathcal{M}_{\mu,\nu}; u, v \right):=   \sum_{i,r} u^r v^i \tr\left(w, IH^{2 r,i}_c\left(\mathcal{M}_{\mu,\nu},\kappa \right)\right).
\]
with $IH^{2 r,i}_c\left(\mathcal{M}_{\mu,\nu},\kappa \right)$ the weight $2r$ graded part of the degree $i$ compactly supported intersection cohomology  of $\mathcal{M}_{\mu,\nu}$.

\begin{conjecture}\label{conjecture_c_intro}
Let $w$ be a $n$-cycle in $\Sym_n$. The coefficient $c_{\mu,\nu}^{1^n}(q,t)$ is related to the $w$-twisted mixed Hodge polynomial of $\mathcal{M}_{\mu,\nu}$ by
\[
c_{\mu,\nu}^{1^n}(q,t) = t^{\frac{{-\dim {\mathcal{M}_{\mu,\nu}}}}{2}} IH_c^{w}\left( \mathcal{M}_{\mu,\nu}, \frac{1}{q}, \sqrt{q t} \right).
\]
\end{conjecture}

In \cite{ballandras_ic} we compute Poincaré polynomial for compactly supported intersection cohomology of the character varieties $\mathcal{M}_{\overline{\boldsymbol{\mathcal{C}}}}$. This allows to prove the following theorem which is the Poincaré polynomial specialization of Conjecture \ref{conjecture_c_intro}.

\begin{theorem}
Let $w$ be a $n$-cycle in $\Sym_n$, the coefficient $c_{\mu,\nu}^{1^n}$ is related to the $w$-twisted Poincaré polynomial of $\mathcal{M}_{\mu,\nu}$ by
\[
c_{\mu,\nu}^{1^n}(1,t) = t^{\frac{{-\dim {\mathcal{M}_{\mu,\nu}}}}{2}} \sum_i t^{\frac{i}{2}} \tr \left( w , \IH_c^i\left(\mathcal{M}_{\mu,\nu},\kappa\right)\right).
\]
\end{theorem}

%However in \cite{ballandras_ic} in order to computer intersection cohomology of character varieties, we prove that the resolution $\mlps$ is diffeomorphic to some character varieties with semisimple monodromies. This allows to transfer the Springer action to the cohomology of such character variety. Therefore we obtain a multiplicative version of the previous theorem, giving an interpretation of the coefficient $c_{\mu,\nu}^{1^n}\left(1,{t}\right)$.

\section*{Aknowledgement}
This work is part of my PhD thesis under the supervision of Emmanuel Letellier and Fernando Rodriguez-Villegas. I am very grateful to both of them for introducing me to the interesting topic of Weyl group actions on the cohomology of character varieties and quiver varieties. I am aslo thankful to François Bergeron for interesting comments and suggestions.  

\section{Symmetric functions}

\subsection{Generalities}
\begin{notations}[Partitions]\label{notations_partitions}
A partition of an integer $n\in\NN$ is a decreasing sequence of non-negative integers 
\[
\lambda=(\lambda_1,\ \lambda_2,\dots, \lambda_{l(\lambda)}) \text{ with } |\lambda|:=\lambda_1+\lambda_2+\dots+\lambda_{l(\lambda)}=n.
\]
The length of $\lambda$ is the number $l(\lambda)$ of non-zero terms. The set of partitions of $n$ is denoted by $\Part_n$ and 
\[
\Part:=\bigcup_{n\in\NN } \Part_n.
\]
The dominance ordering on $\Part$ is defined by $\lambda \preceq \mu$ if and only if
$|\lambda|=|\mu|$ and 
\[
\sum_{i=1}^{k}\lambda_i \le \sum_{i=1}^{k}\mu_i \text{ for all } k\in \NN.
\]
For $\lambda=(\lambda_1,\dots,\lambda_l)$ a partition, we introduce the following notation
\[
\Part_{\lambda} := \Part_{\lambda_1}\times\dots\times \Part_{\lambda_l}.
\]
\end{notations}
\begin{notations}[Young diagrams]
To a partition $\lambda$, we associate the following set
\[
\left\lbrace (i,j) \left| 1\le i \le l(\lambda) \text{ and } 1 \le j \le \lambda_i \right. \right\rbrace.
\]
%A partition is often identified with its Young diagram so that $(i,j)\in \lambda$ means that $(i,j)$ belongs to the Young diagram of $\lambda$.
This set is called the Young diagram of $\lambda$, it gives a graphical way to think about partitions. The transpose of a Young diagram is obtained by permuting $i$ and $j$. The transpose $\lambda'$ of a partition $\lambda$ is the partition with Young diagram the transpose of the Young diagram of $\lambda$. The Young diagram of the partition $\lambda=(6,4,2)$ has the following graphical representation
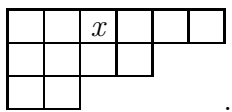
\[
\begin{Young}
      & &x &  &  & \cr
      &  &  & \cr
      & \cr
\end{Young}.
\]
The box $x$ has coordinates $(i,j)=(1,3)$. The arm length of $x$ is number of box to the right of $x$, in this case $a(x)=3$. The leg length is the number of box under $x$, we have $l(x)=1$.
\end{notations}
\begin{notations}[Symmetric functions]
Let $X=(x_1,x_2,\dots)$ be an infinite set of variable and let $\SymF[X]$ be the ring of symmetric functions in $(x_1,x_2,\dots)$. We use the usual notations from Macdonald's book \cite{mcd}. In particular, the usual basis of symmetric functions indexed by partitions: $m_{\lambda}$, $e_{\lambda}$, $h_{\lambda}$, $p_{\lambda}$ and $s_{\lambda}$ are respectively the monomial, elementary, complete, power sum and Schur symmetric functions.

The Hall pairing is denoted by $\left\langle\dots,\dots\right\rangle$ and is defined by
\begin{equation}\label{orthp}
    \left\langle p_\lambda,p_\mu \right\rangle =\delta_{\lambda,\mu}z_\lambda,
\end{equation}
the symbol $\delta_{\lambda,\mu}$ is $1$ if $\lambda$ and $\mu$ are the same partition and $0$ otherwise.  The order of the stabilizer of a permutation of cycle type $\lambda$ is denoted by $z_\lambda$, namely
\begin{equation*}
    z_\lambda=\prod_{l=1}^k i_l^{m_{l}}m_{l}! 
\end{equation*}for a partition $\lambda=(\underbrace{i_1,\dots\ i_1}_{m_{1}},\dots\ ,\underbrace{i_k,\dots\ i_k}_{m_{k}})$ .
\end{notations}

\begin{definition}[Adams operator]
For $n\in\mathbb{Z}_{>0}$, the Adams operator $p_n$ is a  ring endomorphism of $\SymF[X]$. It can be defined by its values on the generating family of power sums,
\[
p_m\left[p_n[X]\right]:=p_{mn}[X]\; \text{ for } m \in \NN_{>0} \text{ and } n \in \NN.
\]
The following notation is commonly used for Adams operators
\[
F\left[X^n\right]:=p_n\left[F[X]\right].
\]
\end{definition}
\begin{remark}
More generally, Adams operator are defined in any lambda ring. In this article, the only lambda rings appearing are rings of symmetric functions and polynomial rings such as $\KK[u]$. On such polynomial rings, the Adams operator $p_n$ is defined by $p_n[u] := u^n$. 
\end{remark}

Let $k$ be a positive integer, we consider the space of multivariate symmetric functions in $k$ infinite sets of variables over $\mathbb{Q}(q,t)$
\[
\SymF\left[X_1,\dots,X_k\right]:=\mathbb{Q}(q,t)\otimes \SymF[X_1]\otimes\dots\otimes \SymF[X_k].
\]
A series with coefficients in this ring of multivariate symmetric functions will conveniently encode cohomological information about comet shaped quiver varieties. The ring of such series is denoted by $\SymF\left[X_1,\dots,X_k\right][[s]]$, it is a lambda ring, and the Adams operators extend to ring endomorphisms of $\SymF\left[X_1,\dots,X_k\right][[s]]$ defined by
\[
p_n\left[f(q,t)F_1\left[X_1\right]\otimes\dots\otimes F_k\left[X_k\right]s^l\right]=f(q^n,t^n)F_1\left[X_1^n\right]\otimes\dots\otimes F_k\left[X_k^n\right]s^{n l}.
\]

\begin{definition}[Plethystic substitution]
    Let $F$ be a symmetric function. The ring of symmetric functions $\SymX$ is freely generated by power sums, so that $F$ can be uniquely obtained as a polynomial expression in the power sums $(p_n)_{n\in \NN}$. Interpreting $p_n$ as the Adams operator, the same polynomial expression defines an operator acting on any lambda ring $\Lambda$, this operator is denoted by $F[\dots]$. For $G\in \Lambda$, the expression $F[G]$ is called a plethystic substitution.
\end{definition}
\begin{remark}
Similarly to Adams operators, the operator $F[\dots]$ naturally extends to $\SymF\left[X_1,\dots,X_k\right][[s]]$. 
\end{remark}
\begin{definition}[Plethystic exponential and logarithm]
The plethystic exponential $\Exp : s \SymF [X_1,\dots,X_k][[s]] \to 1+\SymF [X_1,\dots,X_k][[s]]$ is defined by
\[
\Exp [G] := \exp\left( \sum_{n\ge 1} \frac{p_n[G]}{n} \right),
\]
its inverse, the plethystic logarithm 
\[
\Log : 1+ s \SymF [X_1,\dots,X_k][[s]] \to \SymF [X_1,\dots,X_k][[s]],
\]
is defined by
\[
    \Log[1+G]:=\sum_{n\ge 1}\frac{\mu(n)}{n}p_n \left[\log(1+G)\right],
\]
where $\mu$ is the usual Mobius function.
Contrarily to the ordinary ones, the plethystic exponential and logarithm are written with an uppercase character.
\end{definition}

\begin{remark}
Plethystics operations satisfy the relations
\begin{eqnarray*}
    \Exp[F+G]&=&\Exp[F]\Exp[G],\\
    \Log[(1+F)(1+G)]&=&\Log[1+F]+\Log[1+G], \\
    \Log[\Exp[G]] &=&G.
\end{eqnarray*}
\end{remark}

\subsection{Characters of the symmetric group and symmetric functions}
Let $R_n$ be the vector space spanned by characters of the symmetric group $\Sym_n$ and let $R:=\bigoplus_{n\le 0} R_n$. There is a natural ring structure on $R$ and a pairing such that $R$ is naturally isomorphic to $\SymF[X]$. Let us recall this well-known fact, see \cite{mcd} for more details.

Let $\chi$ and $\eta$ be two characters of $\Sym_n$, and let $V_{\chi}$, respectively $V_{\eta}$ be the associated representations. The pairing is defined by
\[
\left\langle \chi,\eta \right\rangle = \dim \Hom_{\Sym_n} \left(V_{\chi}, V_{\eta}\right).
\]
The spaces $R_m$ and $R_n$ are orthogonal if $m\ne n$. The product of two characters $\chi\in R_m$ and $\eta\in R_n$ is the character of the representation $\Ind_{\Sym_m\times\Sym_n}^{\Sym_{m+n}}V_\chi \otimes V_\eta$, it is denoted by $\chi.\eta$.

The irreducible characters of the symmetric group $\Sym_n$ are indexed by partitions of $n$, they are denoted by $\left(\chi_{\lambda}\right)_{\lambda \in \Part_n}$. The irreducible representation associated to the character $\chi_{\lambda}$ is denoted by $V_{\lambda}$. The indexing is the same as in Macdonald's book \cite{mcd}, so that $V_{(n)}$ is the trivial representation and $V_{(1^n)}$ the sign representation.

\begin{proposition}
    Define the characteristic map $\ch : R \to \SymF[X]$ by $\ch(\chi_{\lambda}) := s_\lambda$. It is an isomorphism between $R$ and $\SymF[X]$ compatible with the products and the pairings, $\SymF[X]$ being endowed with the Hall pairing.
\end{proposition}
\begin{proof}
See Macdonald \cite[I-7]{mcd}.
\end{proof}
\begin{remark}\label{remark_trace_mult}
The last proposition gives a representation theoretic meaning to symmetric functions. Consider $V$ a representation of $\Sym_n$ with character $\chi_V$, then
\begin{itemize}
    \item The pairing $\left\langle s_{\lambda},\ch(\chi_{V})\right\rangle$ gives the multiplicity of the irreducible representation $V_{\lambda}$ in the representation $V$.
    \item The pairing $\left\langle p_{\mu},\ch(\chi_{V})\right\rangle$ gives the trace of an element in $\Sym_n$ with cycle type $\mu$ on the representation $V$.
\end{itemize}
\end{remark}

\begin{definition}[Frobenius characteristic]\label{def_frobenius_char}
We extend the characteristic map $\ch$ to bigraded representations of $\Sym_n$ by adding variables $q$ and $t$ to keep track of the degree. To a bigraded representation of the symmetric group $V=\bigoplus_{(i,j)\in\mathbb{N}^2}V_{i,j}$ is associated a symmetric function over $\mathbb{Z}(q,t)$. This symmetric function is given by
\begin{equation}
    \ch(V)=\sum_{\lambda\in\Part_n}\sum_{(i,j)\in\mathbb{N}^2} \left\langle V_{i,j}, \chi_\lambda \right\rangle q^i t^j s_\lambda,
\end{equation}
where we have identified the representation $V_{i,j}$ with its character so that $\left\langle V_{i,j}, \chi_\lambda \right\rangle$ is the multiplicity of the irreducible representation of type $\lambda$ in $V_{i,j}$. The symmetric function $\ch(V)$ is called the $q,t$-graded Frobenius characteristic of the bigraded representation $V$.
\end{definition}

\subsection{Orthogonality and Macdonald polynomials}
In this section we recall the characterization of modified Macdonald polynomials following Mellit \cite{mellit_p,integrality}.
\subsubsection{Generalities about scalar products on \texorpdfstring{$\SymX$}{Sym[X]}}\label{subsect_pairing}
A scalar product on $\SymX$ is a non-degenerate $\QQ(q,t)$-bilinear form
\begin{equation*}
\begin{array}{ccccc}
(\dots,\dots)^S & : & \SymX\times\SymX & \to & \QQ(q,t) \\
 & & F,G & \mapsto & (F[X],G[X])^S. \\
\end{array}
\end{equation*}
It can be extended to multivariate symmetric functions by specifying the variable acted upon with a lower index
\begin{equation*}
    \begin{array}{ccccc}
(\dots,\dots)^S_X & : & \SymF[X,Y_1,\dotsi,Y_k]\times\SymF[X,Z_1,\dotsi,Z_l] & \to & \SymF[Y_1,\dotsi,Y_k,Z_1,\dotsi,Z_l], \\
\end{array}
\end{equation*}
on pure tensors it reads
\begin{equation*}
    (F[X]\otimes F'[Y_1,\dotsi,Y_k],G[X]\otimes G'[Z_1,\dotsi,Z_l])^S_X:=(F[X],G[X])^S G'[Z_1,\dotsi,Z_l]F'[Y_1,\dotsi,Y_k].
\end{equation*}

\begin{assumption}[Homogeneity]
When considering families of symmetric functions indexed by partitions such as $\left( u_\lambda\right)_{\lambda\in\Part}$, the symmetric function $u_\lambda$ is always assumed to be homogeneous of degree $|\lambda|$.
\end{assumption} 

\begin{definition}[Reproducing kernel]\label{def_kernel}
Let $(u_{\lambda})_{\lambda\in\Part}$ and $(v_{\mu})_{\mu\in\Part}$ be two basis of $\SymF[X]$ dual with respect to a scalar product $(\dots,\dots)^S$. Then the element $K_S[X,Y]\in\SymF[X,Y]$ defined by 
\begin{equation*}
    K_S[X,Y]:=\sum_{\lambda\in\Part}u_\lambda[X]v_\lambda[Y],
\end{equation*}
it is called the reproducing kernel of the scalar product $(\dots,\dots)^S$. It depends only on the scalar product but not on the choice of dual basis, it satisfies
\[
    (K_S[X,Y],F[X])^S_X=F[Y].
\]
\end{definition}

\subsubsection{Hall pairing and \texorpdfstring{$(q,t)$}{(q,t)}-deformations}

\begin{remark}\label{remark_hall_kernel}
Recall that the Hall pairing satisfies
\[
\left\langle p_{\lambda},p_{\mu}\right\rangle =  \delta_{\lambda,\mu} z_{\lambda},
\]
hence $\left( p_\lambda\right)_{\lambda\in\Part}$ and $\left( z_\mu^{-1}p_\mu \right)_{\mu\in\Part}$ are dual basis with respect to this pairing. This gives the reproducing kernel of the Hall pairing:
\begin{equation*}
    \Exp[XY]=\sum_{\lambda\in\Part}p_\lambda[X]\frac{p_\lambda[Y]}{z_\lambda}=\sum_n h_n[XY].
\end{equation*}
\end{remark}

\begin{definition}[$(q,t)$-Hall pairing]
The $(q,t)$-deformation of the Hall pairing is defined by
\begin{eqnarray*}
\left( F[X],G[X]\right)^{q,t}&:=&\left\langle F[X],G[(q-1)(1-t)X]\right\rangle.
\end{eqnarray*}
\end{definition}
\begin{remark}\label{qt_r_kernel}
The reproducing kernel of the $(q,t)$-Hall pairing is
\begin{equation*}
    \Exp\left[\frac{XY}{(q-1)(1-t)}\right].
\end{equation*}
\end{remark}
\begin{proposition}\label{prop_charact_mcd}
Let $M_{\preceq\lambda}$ be the subspace of $\SymX$ spanned by monomial symmetric functions $m_\mu[X]$ with $\mu\preceq\lambda$. The Macdonald polynomials $\left( \tilde{H}_\lambda[X;q,t]\right)_{\lambda\in\Part}$ are uniquely determined by:
\begin{itemize}
    \item Orthogonality $(\Ht_\lambda[X;q,t],\Ht_\mu[X;q,t])^{q,t}=0$ if $\lambda\ne\mu$.
    \item One of the triangularity condition $\Ht_\lambda[X(t-1)]\in M_{\preceq\lambda}$ or $\Ht_\lambda[X(q-1)]\in M_{\preceq\lambda'}$.
    \item Normalization $\Ht[1; q,t]=1$.
\end{itemize}
Moreover
\begin{equation}
   a_\lambda(q,t):= \left(\Ht_\lambda[X;q,t],\Ht_\lambda[X;q,t]\right)^{q,t}=\prod_{x\in \lambda} (q^{a(x)+1}-t^{l(x)})(q^{a(x)}-t^{l(x)+1}),
\end{equation}
where the product is over the Young diagram of $\lambda$ and $a(x)$ is the arm length and $l(x)$ the leg length (see Notations \ref{notations_partitions}).
\end{proposition}
\begin{proof}
\cite{mellit_p} corollary 2.8.
\end{proof}
The modified Macdonald polynomials $\Ht_\lambda\left[X;q,t\right]$ were first introduced by Garsia--Haiman \cite{garsia} as a deformation of other polynomials defined by Macdonald \cite{mcd}. The definition recalled here comes from \cite{mellit_p}.
\begin{definition}[Modified Kostka polynomials]\label{def_kostka_poly}
The modified Kostka polynomials $\left(\widetilde{K}_{\lambda,\rho}(q,t)\right)_{\lambda,\rho \in \Part_n}$ are defined as the coefficients of the transition matrix between the basis of Schur functions and the basis of modified Macdonald polynomials
\[
\tilde{H}_{\rho}[X;q,t] = \sum_{\lambda\in\Part_n} \widetilde{K}_{\lambda,\rho}(q,t) s_{\lambda}.
\]
\end{definition}
\begin{notations}\label{notation_qt_drop}
The variables $(q,t)$ will often be omitted and the modified Kostka polynomial denoted by $\widetilde{K}_{\lambda,\rho}$ and the modified Macdonald polynomial by $\tilde{H}_{\lambda}[X]$.
\end{notations}

\subsection{A result of Garsia--Haiman}\label{subsect_garsia_haiman}

The remaining of this section is devoted to the presentation of a result of Garsia--Haiman \cite[Theorem 3.4]{garsia}. This result is important to study the coefficients $c_{\mu,\nu}^{1^n}$ in the next sections. Even though there are no new results, some proofs are included for convenience of the reader.

\begin{proposition}\label{prop_delta1}
The operator $\Delta_1$ is defined by
\begin{equation*}
    \Delta_1 F\left[X\right]:=F[X]-F\left[X+\frac{(1-q)(1-t)}{z}\right]\Exp\left[-zX \right]\left|_{z^0}\right. ,
\end{equation*}
where $_{|z^0}$ means taking the coefficient in front of $z^0$. This operator acts on modified Macdonal polynomials by
\begin{equation*}
    \Delta_1 \Ht_\lambda \left[X;q,t\right] = (1-t)(1-q)\sum_{(i,j)\in\lambda}q^{j-1} t^{i-1}\Ht_\lambda \left[X;q,t\right].
\end{equation*}
Moreover the following relation holds
\begin{equation}\label{pleth_u}
    \Ht_\lambda \left[1-u;q,t\right]=\prod_{(i,j)\in\lambda}\left(1-uq^{j-1}t^{i-1}\right).
\end{equation}
\end{proposition}
\begin{proof}
\cite[Corollary 3.1 and Theorem 3.2]{garsia}
\end{proof}
\begin{lemma}\label{lemma_uplus}
At first order in $u$
\begin{equation}\label{uplus}
\Ht_\lambda\left[1+u;q,t\right]=1+u\sum_{(i,j)\in\lambda}q^{j-1}t^{i-1}+\mathcal{O}(u^2).
\end{equation}
\end{lemma}
\begin{proof}
One should be careful with plethystic substitutions, to compute the left hand side of \eqref{uplus} one cannot just substitute $-u$ for $u$ in \eqref{pleth_u}. Indeed $p_n[1-u]=1-u^n$ and $p_n[1+u]=1+u^n$ so that substituting $-u$ for $u$ in the latter gives back the former only when $n$ is odd. We denote by $d_{\lambda,\mu}$ the coefficient of $p_\mu$ in the expansion of $\Ht_\lambda$ in the basis of power sums $(p_{\kappa})_{\kappa\in\Part_n}$, then
\begin{eqnarray*}
    \Ht_\lambda \left[1-u;q,t\right]&=&\sum_{|\mu|=|\lambda|}d_{\lambda,\mu} \prod_i (1-u^{\mu_i}), \\
     \Ht_\lambda \left[1+u;q,t\right]&=&\sum_{|\mu|=|\lambda|}d_{\lambda,\mu} \prod_i (1+u^{\mu_i}).
\end{eqnarray*}
We conclude by comparing the coefficient in front of $u$ and using \eqref{pleth_u}.
\end{proof}
\begin{lemma}\label{h_pairing}
Let $F\in\SymnX$ be a symmetric function of degree $n\ge 2$. Then the coefficient in front of $u$ in $F[1+u]$ is given by the Hall pairing with a complete symmetric function
\begin{equation*}
    \left. F[1+u]_{}\right|_{u}=\left\langle h_{(n-1,1)}[X],F[X]\right\rangle.
\end{equation*}
\end{lemma}
\begin{proof}
The coefficient of $m_\lambda$ in the monomial expansion of $F$ is denoted by $c_{\lambda}$. The plethystic substitution $F[1+u]$ corresponds to the evaluation of the symmetric function $F$ on the set of variables $(1,u,0,\dots)$, moreover
\begin{equation*}
    F[1+u]=\sum_{|\lambda|=n}c_\lambda m_\lambda[1+u].
\end{equation*}
Therefore the only $m_\lambda$ contributing are the one with $\lambda$ of length at most two and the coefficient in front of $u$ is $c_{(n-1,1)}$. The conclusion follows as complete symmetric functions and monomial symmetric functions are dual with respect to the Hall pairing.
\end{proof}
\begin{lemma}\label{pn_scal}
Let $F\in\SymnX$ be a symmetric function of degree $n$ then
\begin{equation*}
    \left.\frac{F[1-u]}{1-u}\right|_{u=1}=\left\langle F[X],p_n[X]\right\rangle,
\end{equation*}
where $\left|_{u=1}\right.$ means setting $u=1$.
\end{lemma}
\begin{proof}
Let $d_\lambda$ be the coefficient in front of $p_\lambda$ in the power sum expansion of $F$,
\begin{eqnarray*}
F[1-u]&=&\sum_{|\lambda|=n}d_\lambda p_\lambda[1-u]\\
&=&\sum_{|\lambda|=n}d_\lambda \prod_i (1-u^{\lambda_i}).
\end{eqnarray*}
When dividing by $(1-u)$ and setting $u=1$ all terms coming from partitions of length at least two will vanish as $(1-u)^2$ divides them. Therefore we have
\begin{equation*}
    \left.\frac{F[1-u]}{1-u}\right|_{u=1}=d_{(n)} \left.\frac{1-u^n}{1-u}\right|_{u=1}=n d_{(n)}.
\end{equation*}
The size of the centralizer of a $n$-cycle in $\Sym_n$ is $z_{(n)}=n$, the conclusion follows by orthogonality of power sums \eqref{orthp}.
\end{proof}
Now we can state and recall the proof of an important theorem of Garsia--Haiman.
\begin{theorem}[Garsia--Haiman \cite{garsia} Theorem 3.4]\label{garsia_m_th}
We denote by $\prod_{(i,j)\in\lambda}'$ a product over the young diagram of a partition $\lambda$ omitting the top left corner with $(i,j)=(1,1)$. The following identity holds
\begin{equation}\label{garsia_m}
    (-1)^{n-1} s_{(1^n)}[X]=(q-1)(1-t)\sum_{|\lambda|=n}\frac{\sum_{(i,j)\in\lambda}q^{j-1}t^{i-1}\prod_{(i,j)\in\lambda}'(1-q^{j-1}t^{i-1})\Ht_\lambda[X]}{a_\lambda(q,t)}.
\end{equation}
\end{theorem}
\begin{proof}
The reproducing kernel of the $(q,t)$-Hall pairing is given in Remark \ref{qt_r_kernel}. The degree $n$ term of $\Exp[Z]$ is $h_n[Z]$.
The basis $\left(\Ht_\lambda[X]\right)_{\lambda\in\Part}$ and $\left(\frac{\Ht_\lambda[X]}{a_\lambda}\right)_{\lambda\in\Part}$ are dual with respect to this scalar product. Following Definition \ref{def_kernel} and Remarks \ref{remark_hall_kernel}, \ref{qt_r_kernel}, the degree $n$ term of the reproducing kernel of the $(q,t)$-Hall pairing is
\begin{equation*}
    h_n\left[\frac{XY}{(q-1)(1-t)}\right]=\sum_{|\lambda|=n}\frac{\Ht_\lambda[X]\Ht_\lambda[Y]}{a_\lambda}.
\end{equation*}
Now expand $h_n$ in the basis of power sums, proceed to the substitution $Y=1-u$ and apply \eqref{pleth_u}
\begin{equation*}
    \sum_{|\mu|=n}z_\mu^{-1}p_\mu\left[\frac{X(1-u)}{(q-1)(1-t)}\right]=\sum_{|\lambda|=n}\frac{\Ht_\lambda[X]\prod_{(i,j)\in\lambda}\left(1-uq^{j-1}t^{i-1}\right)}{a_\lambda}.
\end{equation*}
Divide by $(1-u)$, set $u=1$, apply Lemma \ref{pn_scal} to the left hand side and compute explicitly the right hand side:
\begin{equation*}
    \sum_{|\mu|=n}z_\mu^{-1}\left(p_\mu\left[\frac{XY}{(q-1)(1-t)}\right],p_{(n)}[Y]\right)_Y=\sum_{|\lambda|=n}\frac{\Ht_\lambda[X]\prod_{(i,j)\in\lambda}'\left(1-q^{j-1}t^{i-1}\right)}{a_\lambda}.
\end{equation*}
As Adams operator are ring morphisms, we have
\begin{equation*}
   p_\mu\left[\frac{XY}{(q-1)(1-t)}\right]=p_\mu\left[\frac{X}{(q-1)(1-t)}\right]p_\mu [Y],
\end{equation*}
then by orthogonality of power sums \eqref{orthp}
\begin{equation}\label{proof_g1}
    p_{(n)}\left[\frac{X}{(q-1)(1-t)}\right]=\sum_{|\lambda|=n}\frac{\Ht_\lambda[X]\prod_{(i,j)\in\lambda}'\left(1-q^{j-1}t^{i-1}\right)}{a_\lambda}.
\end{equation}
Apply the operator $\Delta_1$ to \eqref{proof_g1}. According to Proposition \ref{prop_delta1}, the operator $\Delta_1$ is diagonal in the basis of Macdonald polynomials and we obtain, up to a sign, the right hand side of \eqref{garsia_m}. Let us compute the left hand side
\begin{eqnarray*}
    \Delta_1 p_{(n)}\left[\tfrac{X}{(q-1)(1-t)}\right]=p_{(n)}\left[\tfrac{X}{(q-1)(1-t)}\right]-p_{(n)}\left[\tfrac{X}{(q-1)(1-t)}-\tfrac{1}{z}\right]\Exp[-zX]\left|_{z^0}\right. \\
    =p_{(n)}\left[\tfrac{X}{(q-1)(1-t)}\right]-p_{(n)}\left[\tfrac{X}{(q-1)(1-t)}\right]\Exp[-z X]\left|_{z^0}\right. +p_{(n)}\left[\tfrac{1}{z}\right]\Exp[-zX]\left|_{z^0}\right. \\
    =\tfrac{1}{z^n}\Exp[-z X]\left|_{z^0}\right. .
\end{eqnarray*}
In the second line we used that the Adams operator $p_n$ is a ring morphism and in the last line that it acts on $z$ by raising to the power $n$. Now $\Exp[-z X]$ is the inverse of $\Exp[z X]$ so that if $X$ is the infinite set of variables $(x_1, x_2, \dots)$, then
\begin{equation*}
    \Exp[-zX]=\prod_i(1-z x_i).
\end{equation*}
The coefficient in front of $z^n$ is $(-1)^n e_n[X]$ so that
\begin{equation*}
    (-1)^n e_n[X]=-(q-1)(1-t)\sum_{|\lambda|=n}\frac{\sum_{(i,j)\in\lambda}q^{j-1}t^{i-1}\prod_{(i,j)\in\lambda}'(1-q^{j-1}t^{i-1})\Ht_\lambda[X]}{a_\lambda}.
\end{equation*}
To conclude, notice that $e_n=s_{(1^n)}$.
\end{proof}

\section{Geometric background}\label{sect_perv_ic}
In this section we recall classical results about intersection cohomology. The main reference is Beilinson--Bernstein--Deligne--Gabber  \cite{bbd}.

\subsection{Notations and generalities on the bounded derived category of constructible sheaves}

The field $\KK$ is either $\C$ or an algebraic closure $\Fqbar$ of a finite field $\Fq$ with $q$ elements. Let $X$ be an algebraic variety over $\KK$ and let $l$ be a prime different from the characteristic of $\KK$. We denote by $\kappa_X$ the constant $l$-adic sheaf on $X$ with coefficients in $\qlbar$. When there are no risk of confusion we just write $\kappa$ instead of $\kappa_X$. For $\KK=\C$ we can also consider the constant sheaf with complex coefficients, in the analytic topology.
\begin{notations}
The bounded derived category of $\kappa$-constructible sheaves on $X$ is denoted by $\dcb{X}$. Its objects are represented by complexes of sheaves $K$ such that the cohomology sheaves $\mathcal H ^i K$ are  $\kappa$-constructible sheaves on $X$ and finitely many of them are non-zero. Let $Y$ be a variety over $\KK$ and let $f:X\to Y$ be a morphism, we have the usual four functors 
\begin{eqnarray*}
f^* , f^! : \dcb Y  \to \dcb X, \\
f_*, f_! : \dcb X \to \dcb Y.
\end{eqnarray*}
\end{notations}
\begin{theorem}[Base change]
Let $K\in\dcb{Y'}$ and consider a cartesian square
\begin{equation}\label{diagram_base_change}
\begin{tikzcd}
 X' \arrow[r,"g"]\arrow[d,"b"'] & \arrow[d,"a"] Y' \\
 X \arrow[r,"f"'] & Y,
\end{tikzcd}
\end{equation}
then the natural morphism $f^* a_! K \to b_! g^* K$ is an isomorphism.
\end{theorem}
\begin{remark}\label{remark_basechange_fiber}
Let $\alpha\hookrightarrow X$ be a geometric point of $X$ and let $\beta$ be its image by $f$. 
The variety $X_{\alpha}:= X'\times_{X} \alpha$ is the fiber of $b$ over $\alpha$ and the variety $Y_{\beta}:= Y'\times_Y \beta$ is the fiber of $a$ over $\beta$.
The cartesian square \eqref{diagram_base_change} induces the following cartesian square where $h$ is an isomorphism
\[
\begin{tikzcd}
 X_{\alpha} \arrow[r,"h"]\arrow[d] & \arrow[d] Y_{\beta} \\
 \alpha \arrow[r] & \beta.
\end{tikzcd}
\]
The base change isomorphism for this diagram identifies with the stalk at $\alpha$ of the base change isomorphism of Diagram \eqref{diagram_base_change},
\[
f^* a_! K_{\alpha} \to b_! g^* K_{\alpha}.
\]
This morphism is nothing but the morphism obtained by functoriality of the compactly supported cohomology
\[
H_c^{\bigcdot}(Y_{\beta},K)\xrightarrow{h^*} H_c^{\bigcdot}(X_{\alpha}, h^*K).
\]
\end{remark}
\begin{definition}\label{def_W_equiv}
Let $W$ be a finite group acting from the left on a variety $X$. For all $w\in W$ there is a morphism $w: X\to X$. An action of $W$ on an element $K\in \dcb{X}$ is the data of isomorphisms $\phi_w : w^* K\cong K$ such that for all $w,w'\in W$,
\begin{equation}\label{eq_W_equiv}
\phi_{w' w}=\phi_w w^*(\phi_{w'}) ,
\end{equation}
and such that $\phi_{1}=\Id$. We say that the complex $K$ is $W$-equivariant.
\end{definition}
\begin{remark}
When the action of $W$ on $X$ is trivial, an action of $W$ on $K\in \dcb X$ is just a group morphism from the opposite group $W^{\op}$ to the group of automorphism $\Aut(K)$.
\end{remark}
\begin{proposition}\label{prop_equivariant_push}
Let $f:X\to Y$ be a $W$-equivariant morphism between varieties with left $W$-action. Let $W$ act on $K$ by morphisms $\phi_w : w^*K\cong K$, then $W$ acts on $f_! K$.
\end{proposition}
\begin{proof}
Base change formulas allow to define the action, for all $w\in W$ they provide an isomorphism $w^*f_! K \to f_! w^* K$. Compose this isomorphism with  $f_! \phi_w$ to obtain an isomorphism $\widetilde{\phi}_w : w^* f_! K\to f_! K$. The compatibility \eqref{eq_W_equiv} follows from functoriality of base change.
\end{proof}

\subsection{Intersection cohomology}

\begin{definition}[Intersection complex]\label{def_intersection_complex}
Let $Y\hookrightarrow X$ be a closed embedding and let $j:U\hookrightarrow Y$ be an open embedding. Assume that $U$ is smooth, irreducible and that $\overline U = Y$. Let $\xi$ be a local system on $U$. The intersection complex $\underline{\mathcal{I C}}_{Y,\xi}^\bullet$ is the unique (up to isomorphism)  element $K$ in $\dcb{Y}$ characterized by
\begin{eqnarray*}
\hi K &=& 0 \text{ \ if \ } i < -\dim Y, \\
\mathcal{H}^{-\dim Y} K_{\left| U\right.}&=&\xi, \\
\dim \left(\supp \hi K \right) &<& -i \text{ \ if \ } i > -\dim Y ,\\
\dim \left(\supp \hi D_Y K \right) &<& -i \text{ \ if \ } i > -\dim Y.
\end{eqnarray*}
We also denote by $\underline{\mathcal{I C}}_{Y,\xi}^\bullet$ its extension $j_*\underline{\mathcal{I C}}_{Y,\xi}^\bullet$.  
\end{definition}
\begin{remark}[Continuation principle]\label{remark_continuation_principle}
The intersection complex of $\xi$ can also be defined as the intermediate extension $\underline{\mathcal{I C}}_{Y,\xi}^{\bullet}= j_{! *} \xi $. Moreover the functor $j_{! *}$ is fully faithful (see Kiehl-Weissauer \cite[III - Corollary 5.11]{Kiehl_Weissauer}).  
\end{remark}
\begin{remark}
The intersection complex does not depend on the choice of smooth open subset in $Y$. When the local system $\xi$ is not specified, it is chosen to be the constant sheaf $\kappa_U$ and $\uic X := \uic {X,\kappa_U}$. We denote by $\mathcal{I C}_X^{\bullet}$
\end{remark}
\begin{notations}
The shifted intersection complexes are
\[
\mathcal{I C}_{X,\xi}^{\bullet} := \underline{\mathcal{I C}}_{X,\xi}^{\bullet - \dim X} \text{ and }
\mathcal{I C}_X^{\bullet} := \underline{\mathcal{I C}}_X^{\bullet - \dim X}.
\]
\end{notations}

\begin{definition}[Intersection cohomology]
Let $p:X\to \Spec \KK$ be the structural morphism and $k$ an integer. The $k$-th intersection cohomology space of $X$ is
\[
IH^k(X,\kappa) := \mathcal{H}^{k-\dim X} p_* \uic X = \mathcal{H}^k p_* \mathcal{I C}_X^{\bullet}
\]
and the $k$-th compactly supported intersection cohomology space of $X$ is
\[
\IH^k_c (X,\kappa) := \mathcal{H}^{k-\dim X} p_! \uic X = \mathcal{H}^k p_! \mathcal{I C}_X^{\bullet}
\]

\end{definition}

For $\KK=\C$, Saito \cite{saito1986} proved that the intersection cohomology spaces carry a mixed Hodge structure. Thus there exists on $I H_c^k(X,\mathbb{Q})$ an increasing finite filtration called the weight filtration and denoted by $W^k_{\bigcdot}$  such that the complexified quotient $\mathbb{C}\otimes_{\mathbb{Q}}W^k_{r}/ W^k_{r-1}$ carries a pure Hodge structure of weight $r$. The Hodge numbers of this structure are denoted $h_c ^{i,j,k}(X)=\dim I H_c^{i,j,k}(X,\C)$ and satisfy $i+j=r$.
\begin{definition}\label{def_polynomial_coho}
The mixed Hodge structure is encoded in the mixed Hodge polynomial,
\begin{equation}
    I H_c\left(X;x,y,v\right):=\sum_{i,j,k}  h_c^{i,j,k}(X) x^i y^j v^k.
\end{equation}
This polynomial has an important specialisation, the Poincaré polynomial
\begin{equation}
    P_c(X;v):=I H_c\left(X;1,1,v\right)=\sum_k \dim I H_c^k(X,\kappa) v^k.
\end{equation}
In this article \textit{"Poincaré polynomial"} always refers to \textit{"Poincaré polynomial for compactly supported intersection cohomology"}.
\end{definition}

\section{Main objects and notations}

\subsection{Adjoint orbits in \texorpdfstring{$\gl_n$}{gl n}}\label{sect_conjugacy_classes}
\subsubsection{Notations for adjoint orbits}\label{subsect_notations_adjoint_conj}
The goal of this article is to relate some geometric objects to combinatorial data. The first step is the well-known labelling of adjoint orbits by their Jordan types, which we recall in order to fix the notations.

For an integer $r$ and for $ z \in \KK$, we denote by $J_r (z)\in\gl_r$ the Jordan block of size $r$ with eigenvalue $z$ so that $J_r(z)-z \Id_r$ is nilpotent of order $r$. Let $\mu=(\mu_1,\mu_2,\dots,\mu_s)$ be a partition of an integer $m$ and let $z \in \KK$. Let $J_{\mu} (z)$ be the matrix with eigenvalue $z$ and Jordan blocks of sizes given by $(\mu_j)_{1\le j\le s}$,
\[
J_{\mu} (z) :=      
\begin{pmatrix}
    J_{\mu_1}(z) & &  &\\
    & J_{\mu_2}(z) & & \\
    & & \ddots & \\
     & & & J_{\mu_s}(z)
  \end{pmatrix} \in \gl_{m}.
\]
Let $\nu=(\nu_1,\dots,\nu_l) \in \Part_n$ be a partition of $n$, introduce the following notation
\[
\Part_{\nu} := \Part_{\nu_1}\times\Part_{\nu_2} \times\dots\times \Part_{\nu_l}.
\]
Consider a diagonal matrix $\sigma$,
\begin{equation}\label{eq_def_sigma}
\sigma = \begin{pmatrix}
    \sigma_1 \Id_{\nu_1} & & & \\
    & \sigma_2 \Id_{\nu_2} & & \\
    &  & \ddots & \\
    &  &   &  \sigma_l \Id_{\nu_l}
\end{pmatrix},
\end{equation}
with $\sigma_i \ne \sigma_j$ for $i\ne j$, so that $\nu_i$ is the multiplicity of the eigenvalue $\sigma_i\in \KK$. 
\begin{notations}\label{notations_jmus}
Consider an element $\underline{\mu}=\left(\mu^{1},\dots,\mu^{l}\right)$ in $\Part_{\nu}$,
we denote by $\mathcal{O}_{\underline{\mu},\sigma}$ the adjoint orbit of the matrix
\[
J_{\underline{\mu},\sigma}:=\begin{pmatrix}
   J_{\mu^{1}}(\sigma_1 ) & & & \\
    & J_{\mu^{2}}(\sigma_2)  & & \\
    &  & \ddots & \\
    &  &   &  J_{\mu^{l}}(\sigma_l).
\end{pmatrix}.
\]

\end{notations}
Let us recall a well-known proposition.
\begin{proposition}
The Zariski closure of the adjoint orbit $\mathcal{O}_{\underline{\mu},\sigma}$ is
\[
\overline{\mathcal{O}}_{\underline{\mu},\sigma} = \bigsqcup_{\underline{\rho}\preceq\underline{\mu} } \mathcal{O}_{\underline{\rho},\sigma},
\]
the union is over the set of $l$-tuples $\underline{\rho}=\left(\rho^{1},\dots,\rho^{l} \right)$ with
$
\rho^{j}\preceq \mu^{j} \; \text{ for  } 1\le j \le l
$. The dominance order on partition was recalled in \ref{notations_partitions}. \end{proposition}

Over $\Fqbar$, there is a more precise description of adjoint orbits. Denote by $F$ the Frobenius endomorphism of $\gl_n(\Fqbar)$ raising the coefficients to the power $q$.

\begin{definition}[Type of an $F$-stable adjoint orbit]\label{def_type_fq_cc_ao}
Let $\mathcal{O}$ be an $F$-stable adjoint orbit in $\gl_n(\Fqbar)$, i.e., an orbit such that $F(\mathcal{O})\subset \mathcal{O}$. Then the set $\mathcal{O}^F$ of fixed points in $\mathcal{O}$ under the Frobenius is not empty. The characteristic polynomial of $\mathcal{O}$ has its coefficients in $\Fq$ so that its eigenvalues, which live in $\Fqbar$, are permuted by the Frobenius. The spectrum of $\mathcal{O}$, with multiplicities, reads 
\begin{equation*}
\left(\underbrace{\left(\gamma_1,\dots,\gamma_1^{q^{d_1-1}}\right),\dots,\left(\gamma_1,\dots,\gamma_1^{q^{d_1-1}}\right)}_{m_1},\dots,\underbrace{\left(\gamma_l,\dots,\gamma_l^{q^{d_l-1}}\right),\dots,\left(\gamma_l,\dots,\gamma_l^{q^{d_l-1}}\right)}_{m_l}\right),
\end{equation*}
where $\gamma_i\in\Fqbar$ is such that $\gamma_i^{q^{d_i-1}}\ne\gamma_i$, $\gamma_i^{q^{d_i}}=\gamma_i$ and $\gamma_i\ne\gamma_j$ for $i\ne j$. Then the orbit ${\mathcal{O}}$ determines partitions $\omega^i\in \Part_{m_i}$ giving the size of the Jordan blocks of the Frobenius orbit of eigenvalues $\left(\gamma_i,\dots,\gamma_i^{q^{d_i-1}}\right)$. Up to ordering, it defines a sequence $\omega=(d_1,\omega_1)\dots(d_l,\omega_l)$ in $\mathbb{Z}_{>0}\times \Part$ called the type of the adjoint orbit.
\end{definition}

\subsubsection{Resolutions of Zariski closures of adjoint orbits}\label{subsect_resolution_conjugacy}
In this section we recall the construction of resolutions of closures of adjoint orbits. 
The references for this construction are Kraft--Procesi \cite{Kraft_Procesi_80}, Nakajima \cite{Nakajima_1998,nakajima_01}, Crawley-Boevey \cite{Crawley-Boevey_parabolic, crawley-boevey2003} and Shmelkin \cite{shmelkin2009} (see also Letellier \cite{letellier_gl}).

Using the notations from the previous section, consider an adjoint orbit $\mathcal{O}_{\underline{\mu},\sigma}$. The matrix $\sigma\in\gl_n$ is diagonal as in \eqref{eq_def_sigma} and we denote by $M$ its stabilizer in $\GL_n$,
\[
M=\begin{pmatrix}
   \GL_{\nu_1} & 0 & \\
   0 & \GL_{\nu_2} &  \\
   \vdots & 0 & \ddots 
\end{pmatrix}.
\]
Let $\underline{\mu}=(\mu^{1},\dots,\mu^{l})\in \Part_\nu$ so that $\mu^{i}$ is a partition of the integer $\nu_i$. The transposed of the partition $\mu^i$ is denoted by ${\mu^{i}}'=\left({\mu^{i}_1}',{\mu^{i}_2}',\dots \right) $. Let $L$ be the subgroup of $\GL_n$ formed by block diagonal matrices with blocks of size ${\mu^{i}_r}'$, it is a subgroup of $M$ with the following form
\[
L=\left(
\begin{array}{ccccc}
\overbrace{
\begin{array}{ccccc}
  \GL_{{\mu^{1}_1}' }  & 0 &   \\
   0  &  \GL_{{\mu^{1}_2}' } &  \\
 \vdots & 0 & \ddots 
\end{array}
}^{\nu_1}  &     &   \\
  &  \overbrace{
\begin{array}{ccccc}
  \GL_{{\mu^{2}_1}' }  & 0 &   \\
   0  &  \GL_{{\mu^{2}_2}' } &  \\
 \vdots & 0 & \ddots 
\end{array}
}^{\nu_2}  &   \\
  &    &  \ddots 
\end{array} 
\right).
\]
\begin{notations}\label{notations_group_part}
For a  partition $\nu=(\nu_1,\dots,\nu_l)$,  define 
\[
\Sym_{\nu} := \Sym_{\nu_1}\times\dots\times \Sym_{\nu_l} \; \text{ and } \; \GL_{\nu}:=\GL_{\nu_1}\times\dots\times\GL_{\nu_l}.
\]
Now for $\underline{\rho}=(\rho^1,\dots,\rho^l)\in \Part_{\nu}$, we use the following notations,
\[
\GL_{\underline{\rho}}:=\GL_{\rho^1}\times\dots\GL_{\rho^l}=\prod_{r,s}\GL_{\rho^r_s}
\]
and
\[
\Sym_{\underline{\rho}} :=\Sym_{\rho^1}\times\dots\Sym_{\rho^l}=\prod_{r,s}\Sym_{\rho^r_s}.
\]
Then the previously introduced Levi subgroups satisfy $M\cong \GL_{\nu}$ and $L\cong \GL_{\underline{\mu}'}$.
\end{notations}
Denote by $P$ the parabolic subgroup of blocks upper triangular matrices having $L$ as a Levi factor, then $P=L U_P$ with
\[
U_P =\left(
\begin{array}{ccccc}
\overbrace{
\begin{array}{ccccc}
  \Id_{{\mu^{1}_1}' }  & * &   \\
   0  &  \Id_{{\mu^{1}_2}' } &  \\
 \vdots & 0 & \ddots 
\end{array}
}^{\nu_1}  &  \text{\Large{*}} &   \\
  &  \overbrace{
\begin{array}{ccccc}
  \Id_{{\mu^{2}_1}' }  & * &   \\
   0  &  \Id_{{\mu^{2}_2}' } &  \\
 \vdots & 0 & \ddots 
\end{array}
}^{\nu_2}  &   \\
  &    &  \ddots 
\end{array} 
\right),
\]
and the Lie algebra counterpart of this Levi decomposition is $\mathfrak{p}=\mathfrak{l}\oplus \mathfrak{u}_P$.
\begin{proposition}[Resolutions of Zariski closures of conjugacy classes]\label{prop_resolu_adjoint}
Consider
\[
\mathbb{\widetilde{Y}}_{L,P,\sigma} :=\left\lbrace (X, g P) \in \gl_n \times  \GL_n / P \left| g\inv X g \in \sigma + \mathfrak{u}_P \right.\right\rbrace .
\]
The image of the projection to the first factor $\mathbb{\widetilde{Y}}_{L,P,\sigma}\to \gl_n   $ is the Zariski closure of the adjoint orbit $\mathcal{O}_{\underline{\mu},\sigma}$. Moreover the following map is a resolution of singularities
\[
\map{p^{\sigma}}{\mathbb{\widetilde{Y}}_{L,P,\sigma}}{\overline{\mathcal{O}}_{\underline{\mu},\sigma}}{(X,g P)}{X}.
\]
\end{proposition}

\begin{remark}\label{remark_ss_resol}
If $M$ (which is defined as the stabilizer of $\sigma$ in $\GL_n$) is exactly $L$ then the adjoint orbit $\mathcal{O}_{\underline{\mu},\sigma}$ is semisimple. Semisimple orbits in $\gl_n$ are closed and smooth, for such orbits, $p^{\sigma}$ is an isomorphism.
\end{remark}

\begin{remark}\label{remark_springer_orbit}
The decomposition theorem and Springer theory (or more generally Lusztig's parabolic induction \cite{Lusztig1984,lusztig_character_sheaves,Lusztig1986}) provide more information about the previous resolution of singularities in terms of Weyl group representations. For $G$ a reductive group and for $T$ a maximal torus in $G$ the Weyl group is denoted by
\[
W_G := N_G(T)/T.
\]
The Weyl group of $M$ is  $W_M\cong \prod_i \Sym_{\nu_i}$ and let $V_{\underline{\rho}}$ be the representation $\bigotimes_i V_{\rho^i}$ of $W_M$. The Weyl group of $L$ is $W_L\cong \prod_{i,j} \Sym_{{\mu_j^i}'}$ and let $\epsilon$ be its sign representation. Lusztig's parabolic induction provides the following decomposition
\[
H_c^{\bullet + d_{\underline{\mu}}}\left( \mathbb{\widetilde{Y}}_{L,P,\sigma} ;\kappa \right)= \bigoplus_{\underline{\rho} \preceq \underline{\mu} }  \Hom_{W_M}\left(\Ind_{W_L}^{W_M} \epsilon , V_{\underline{\rho}}\right)\otimes
H_c^{\bullet+ d_{\underline{\rho}}}\left(\overline{\mathcal{O}}_{\underline{\rho},\sigma};\kappa\right),
%\uic{\overline{\mathcal{O}}_{\underline{\rho},\sigma}}
\]
with $d_{\underline{\mu}}$ the dimension of $\overline{\mathcal{O}}_{\underline{\mu},\sigma}$. Letellier \cite{letellier_gl} constructed an action of the relative Weyl group 
\[
W_M(L)=N_M(L)/L
\]
on the spaces $\Hom_{W_M}\left(\Ind_{W_L}^{W_M} \epsilon , V_{\underline{\rho}}\right)$.
\end{remark}

\subsection{The  \texorpdfstring{varieties     $\mathcal{ Q }_{\overline{\mathcal{O}}_{\boldsymbol{\mu},\boldsymbol{\sigma}}}$}{comet-shaped quiver varieties} and their resolutions}

The varieties we are interested in are additive analogues of the character varieties classifying representations of the fundamental group of a compact Riemann surface with $k$ punctures.  They were studied by Crawley-Boevey \cite{crawley-boevey2003,CB_surj} in the case $g=0$, by Hausel, Letellier and Rodriguez-Villegas \cite{aha} for semisimple adjoint orbits and by Letellier \cite{letellier_gl} in general.

For each puncture   an adjoint orbit  $\mathcal{O}_{\underline{\mu}^j,\sigma^j}\subset\gl_n$ is fixed.
A bold symbol is used to represent $k$-tuple:
\begin{equation}\label{eq_bold_kuple}
    \begin{array}{ccc}
\boldsymbol{\mu}& := & \left( \underline{\mu}^1,\dots, \underline{\mu}^k\right),  \\
\boldsymbol{\sigma} & := & \left(\sigma^1,\dots,\sigma^k\right), \\
\mathcal{O}_{\boldsymbol{\mu},\boldsymbol{\sigma}} & := & \left(\mathcal{O}_{\underline{\mu}^1,\sigma^1}, \dots, \mathcal{O}_{\underline{\mu}^k,\sigma^k} \right).
\end{array}
\end{equation}
\begin{definition}[Comet-shaped quiver varieties]
Consider the variety
\begin{multline*}
    \mathcal{V}_{\overline{\mathcal{O}}_{\boldsymbol{\mu},\boldsymbol{\sigma}} }:=\left\lbrace \left(A_1,B_1,\dots,A_g,B_g,X_1,\dots,X_k\right)\in \gl_n^{2 g}\times \overline{\mathcal{O}}_{\underline{\mu}^1,\sigma^1}\times\dots \times \overline{\mathcal{O}}_{\underline{\mu}^k, \sigma^k  } \; \right| \\ 
    \left. \sum_{i=1}^g [A_i,B_i] +\sum_{j=1}^k X_j = 0  \right\rbrace.
\end{multline*}
This is an affine variety acted upon by $\GL_n$ by coordinate-wise adjoint action. The center of $\GL_n$ acts trivialy so that the action factors through a $\PGL_n$-action. The main focus of the article is the following GIT quotient,
\begin{equation}
    \mathcal{ Q }_{\overline{\mathcal{O}}_{\boldsymbol{\mu},\boldsymbol{\sigma}}}:=\left.\left.\mathcal{V}_{\overline{\mathcal{O}}_{\boldsymbol{\mu},\boldsymbol{\sigma}}} \  \right/\right/ \PGL_n=\Spec\left(\KK\left[\mathcal{V}_{\overline{\mathcal{O}}_{\boldsymbol{\mu},\boldsymbol{\sigma}}}\right]^{\GL_n}\right).
\end{equation}
We call the varieties $\mathcal{Q }_{\overline{\mathcal{O}}_{\boldsymbol{\mu},\boldsymbol{\sigma}}}$ the comet-shaped quiver varieties because of their interpretation as Nakajima's quiver varieties that we will recall in Section \ref{sect_Nakajima_quiver_varieties}.
\end{definition}

\begin{definition}[Generic adjoint orbits]\label{def_generic_adjoint}
Denote by $\Delta(\sigma^j)$ the multiset of eigenvalues of $\sigma^j$ repeated according to multiplicities, $\sigma^j_r$ appears exactly $\nu^j_r$ times in the multiset $\Delta(\sigma^j)$. The $k$-tuple of adjoint orbits $\mathcal{O}_{\boldsymbol{\mu},\boldsymbol{\sigma}}$ is generic if and only if it satisfies the two following conditions:
\begin{enumerate}
    \item \[ \sum_{j=1}^k\sum_{\alpha \in \Delta (\sigma^j )}\alpha = 0, \]
    \item for any $r \le n-1$ and for all $(R_1,\dots, R_k)$ with $R_j\subset \Delta(\sigma^j)$ of size $r$
    \[
    \sum_{j=1}^k \sum_{\alpha\in R_j} \alpha \ne 0.
    \].
\end{enumerate}
\end{definition}
\begin{definition}[Generic conjugacy classes]\label{def_generic_cc}
If all the eigenvalues in $\boldsymbol{\sigma}$ are non-zero, then the adjoint orbits are also conjugacy classes in $\GL_n$. They are denoted by $\mathcal{C}_{\boldsymbol{\mu},\boldsymbol{\sigma}}=\left(\mathcal{C}_{\underline{\mu}^1,\sigma^1}, \dots, \mathcal{C}_{\underline{\mu}^k,\sigma^k}\right)$ instead of $\mathcal{O}_{\boldsymbol{\mu},\boldsymbol{\sigma}}$. A $k$-tuple of conjugacy classes $\mathcal{C}_{\boldsymbol{\mu},\boldsymbol{\sigma}}$ is generic if it satisfies the two following conditions:
\begin{enumerate}
    \item \[ \prod_{j=1}^k\prod_{\alpha \in \Delta (\sigma^j )}\alpha = 1, \]
    \item for any $r \le n-1$, for all $(R_1,\dots, R_k)$ with $R_j\subset \Delta(\sigma^j)$ of size $r$
    \[
    \prod_{j=1}^k \prod_{\alpha\in R_j} \alpha \ne 1.
    \].
\end{enumerate}

\end{definition}

\begin{proposition}[\cite{letellier_gl} Proposition 5.2.4, Corollary 5.2.5]
Let $\mathcal{V}_{ \mathcal{O}_{\boldsymbol{\mu},\boldsymbol{\sigma}} }:= \mathcal{V}_{\overline{\mathcal{O}}_{\boldsymbol{\mu},\boldsymbol{\sigma}}}\cap \gl_n^{2 g}\times \mathcal{O}_{\underline{\mu}^{1},\sigma^1}\times\dots\times \mathcal{O}_{\underline{\mu}^{k},\sigma^k}$, and let $\mathcal{Q}_{ \mathcal{O}_{\boldsymbol{\mu},\boldsymbol{\sigma}} }$ be the image of $\mathcal{V}_{ \mathcal{O}_{\boldsymbol{\mu},\boldsymbol{\sigma}} }$ in $\mathcal{Q}_{ \overline{\mathcal{O}}_{\boldsymbol{\mu},\boldsymbol{\sigma}} }$.
Assume that $\overline{\mathcal{O}}_{\boldsymbol{\mu},\boldsymbol{\sigma}}$ is generic, then
\[
\mathcal{Q}_{\overline{\mathcal{O}}_{\boldsymbol{\mu},\boldsymbol{\sigma}}} = \bigsqcup_{\boldsymbol{\rho}\preceq\boldsymbol{\mu}} \mathcal{Q}_{{\mathcal{O}}_{\boldsymbol{\rho},\boldsymbol{\sigma}}}
\]
is a stratification of $\mathcal{Q}_{\overline{\mathcal{O}}_{\boldsymbol{\mu},\boldsymbol{\sigma}}}$ with smooth strata. Moreover, if it is non-empty,
\[
\dim \mathcal{Q}_{\overline{\mathcal{O}}_{\boldsymbol{\mu},\boldsymbol{\sigma}}} = d_{\boldsymbol{\mu}} = n^2 (2 g -2) + 2 +\sum_{j=1}^k \dim \mathcal{O}_{\underline{\mu}^j,\sigma^j}.
\]
\end{proposition}

As before, $\sigma^j$ is a diagonal matrix with stabilizer $M^j := Z_{\GL_n}(\sigma^j)$ such that with Notations \ref{notations_group_part},
\[
M^j \cong \GL_{\nu^j}
\]
for some partition $\nu^j\in \Part_n$ given by the multiplicities of the eigenvalues of $\sigma^j$. The Jordan type of the eigenvalue $\sigma^j_i$ in the adjoint orbit $\mathcal{O}_{\underline{\mu}^j,\sigma^j}$ is $\mu^{j,i}\in \Part_{\nu^{j}_i}$. Denote by ${\mu^{j,i}}'=\left({\mu^{j,i}_1}',{\mu^{j,i}_2}',\dots \right)$ the transposed partition. Let $L^j\subset M^j$ be the subgroup of block diagonal matrices as in \ref{subsect_resolution_conjugacy},
\[
L^j \cong \underbrace{\GL_{{\mu^{j,1}_1}'}\times \GL_{{\mu^{j,1}_2}'} \times \dots }_{\subset \GL_{\nu^j_1}}\times \dots \times \underbrace{\GL_{{\mu^{j,l_j}_1}'}\times \GL_{{\mu^{j,l_j}_2}'} \times \dots}_{\subset \GL_{\nu^j_{l_j}}}.
\]
Let $\mathbb{\widetilde{Y}}_{L^j,P^j,\sigma^j}$ be a resolution of $\overline{\mathcal{O}}_{\underline{\mu}^j,\sigma^j}$ as constructed in \ref{subsect_resolution_conjugacy}.  Let $\boldsymbol{L}:=\prod_{j=1}^k L^j$ and let $\boldsymbol{P}:=\prod_{j=1}^k P^j$, then define
\[
\mathbb{\widetilde{Y}}_{\boldsymbol{L},\boldsymbol{P},\boldsymbol{\sigma}} := \prod_{1\le j\le k} \mathbb{\widetilde{Y}}_{L^j,P^j,\sigma^j}.
\]
Letellier \cite{letellier_gl} constructed resolutions of singularities of $\mathcal{Q}_{\overline{\mathcal{O}}_{\boldsymbol{\mu},\boldsymbol{\sigma}}}$
\begin{definition}[Resolutions of $\mathcal{Q}_{\overline{\mathcal{O}}_{\mu^j,\sigma^j}}$]\label{def_resol_char}
Define
\begin{multline}
   \qlps := \left\lbrace(A_i, B_i)_{1\le i \le g}, (X_j, g_j P^j)_{1\le j \le k} \in \gl_n^{2 g}\times \mathbb{\widetilde{Y}}_{\boldsymbol{L},\boldsymbol{P},\boldsymbol{\sigma}} \right. \\
   \left| \sum_{i=1}^g [A_i, B_i] + \sum_{j=1}^k X_j = 0 \right\rbrace\git\PGL_n.
\end{multline}
The action of $\PGL_n$ on $g_j P^j$ is by left multiplication. The maps $p^{\sigma^j} : \mathbb{\widetilde{Y}}_{L^j,P^j,\sigma^j} \to \overline{\mathcal{O}}_{\underline{\mu}^j,\sigma^j}$ induce a map
\[
p^{\boldsymbol{\sigma}} :    \qlps \to \mathcal{Q}_{\overline{\mathcal{O}}_{\boldsymbol{\mu},\boldsymbol{\sigma} }  },
\]
this morphism is a resolution of singularities.
\end{definition}
\begin{remark}\label{remark_qss_resol}
Similarly to Remark \ref{remark_ss_resol}, if $L^j=M^j$ for $1\le j \le k$, then the adjoint orbit $\mathcal{O}_{\boldsymbol{\mu},\boldsymbol{\sigma}}$ are semisimple and $p^{\boldsymbol{\sigma}}$ is an isomorphism.
\end{remark}

\begin{remark}\label{remark_springer_comet}
Similarly to Remark \ref{remark_springer_orbit},
\[
H_c^{\bullet +d_{\boldsymbol{\mu}}}\left(\qlps, \qlbar\right)=
\bigoplus_{\boldsymbol{\rho}\preceq \boldsymbol{\mu}}V_{\boldsymbol{\mu},\boldsymbol{\rho}} \otimes IH_c^{\bullet +d_{\boldsymbol{\rho}}}\left({\mathcal{Q}_{\overline{\mathcal{O}}_{\boldsymbol{\rho},\boldsymbol{\sigma}}}}, \qlbar \right),
\]
with $V_{\boldsymbol{\mu},\boldsymbol{\rho}}:=\bigotimes_{j=1}^k \Hom_{W_{M^j}}\left(\Ind_{W_{L^j}}^{W_{M^j}} \epsilon , V_{\underline{\rho}^j}\right) $. Therefore, as in Remark \ref{remark_springer_orbit}, Letellier constructed an action of the relative Weyl group $W_{\boldsymbol{M}}(\boldsymbol{L})=\prod_{j=1}^k W_{M^j}(L^j)$ on the cohomology of the resolution $\qlps$, we call this action the Springer action.

\end{remark}

\begin{definition}\label{def_char_var}
Similar constructions exist in the multiplicative case (see Letellier \cite{zclosure}). If all the eigenvalues in $\boldsymbol{\sigma}$ are non-zero, then the adjoint orbits are actually conjugacy classes in $\GL_n$, they are denoted with a $\mathcal{C}$ instead of an $\mathcal{O}$. For a generic $k$-tuple of conjugacy classes $\overline{ \mathcal{C}  }_{ \boldsymbol{\mu} ,\boldsymbol{\sigma}  }$, the character variety is defined by
\begin{multline*}
\mathcal{M}_{\overline{ \mathcal{C}  }_{ \boldsymbol{\mu} ,\boldsymbol{\sigma}  } }:=\left\lbrace\left(A_1,B_1,\dots,A_g,B_g,X_1,\dots,X_k\right)\in \GL_n^{2g}\times \overline{\mathcal{C}}_{\underline{\mu}^1,\sigma^1}\times\dots\times \overline{\mathcal{C}}_{\underline{\mu}^k,\sigma^k}\right| \\ \left. A_1 B_1 A_1\inv B_1\inv \dots A_g B_g A_g\inv B_g\inv X_1\dots X_k =\Id \right\rbrace\git\PGL_n.
\end{multline*}
It admits a resolution of singularities
\begin{multline*}
    \mlps:=\left\lbrace(A_i, B_i)_{1\le i \le g}, (X_j, g_j P^j)_{1\le j \le k} \in \GL_n^{2 g}\times \mathbb{\widetilde{Y}}_{\boldsymbol{L},\boldsymbol{P},\boldsymbol{\sigma}} \right. \\
   \left|  A_1 B_1 A_1\inv B_1\inv \dots A_g B_g A_g\inv B_g\inv X_1\dots X_k =\Id \right\rbrace\git\PGL_n,
\end{multline*}
and the cohomology of the resolution admits the following decomposition
\begin{equation}\label{eq_springer_char}
H_c^{\bullet +d_{\boldsymbol{\mu}}}\left(\mlps, \qlbar\right)= 
\bigoplus_{\boldsymbol{\rho}\preceq \boldsymbol{\mu}}V_{\boldsymbol{\mu},\boldsymbol{\rho}} \otimes IH_c^{\bullet +d_{\boldsymbol{\rho}}}\left({\mathcal{M}_{\overline{\mathcal{C}}_{\boldsymbol{\rho},\boldsymbol{\sigma}}}}, \qlbar \right).
\end{equation}
\end{definition}

\subsection{Intersection cohomology of the varieties \texorpdfstring{$\mathcal{Q}_{\overline{\mathcal{O}}_{\boldsymbol{\mu},\boldsymbol{\sigma} }} $ }{ Q}}

\begin{definition}[Hausel-Letellier-Villegas kernel]\label{def_omega_hlv}
Let $k\in\ZZ_{>0}$ and let $g\in\ZZ_{\ge0}$, the $k$-points, genus $g$ Cauchy function is defined by
\begin{equation}
    \Omega^g_k(z,w):=\sum_{\lambda\in\Part}\mathcal{H}_\lambda(z,w)\prod_{i=1}^k\Ht_\lambda\left[X_i,z^2,w^2\right]s^{|\lambda|},
\end{equation}
with
\begin{equation}
\mathcal{H}_\lambda(z,w):=\prod\frac{\left(z^{2a+1}-w^{2l+1}\right)^{2g}}{\left(z^{2a+2}-w^{2l}\right)\left(z^{2a}-w^{2l+2}\right)}.
\end{equation}
The degree $n$ Hausel--Letellier--Villegas kernel is defined by
\[
\mathbb{H}^{H L V}_n(z,w):=\left.(z^2-1)(1-w^2) \Log \Omega_k^g(z,w)\right|_{s^n}.
\]
\end{definition}
This kernel was introduced to describe the cohomology of character varieties for genus $g$ Riemann surfaces with $k$ punctures. It also describes the cohomology of the varieties $\mathcal{Q}_{\overline{\mathcal{O}}_{\boldsymbol{\mu},\boldsymbol{\sigma} }} $. Consider a generic $k$-tuple of adjoint orbits $\mathcal{O}_{\boldsymbol{\mu},\boldsymbol{\sigma}}$ , with $\boldsymbol{\mu}=(\underline{\mu}^1,\dots,\underline{\mu}^k)$ and $\underline{\mu}^j=\left(\mu^{j,1},\dots,\mu^{j,l_j}\right)$. The transposed of the partition $\mu^{j,i}\in \Part_{\nu^j_i}$ is denoted by ${\mu^{j,i}}'$ and we define the following symmetric function
\begin{equation}\label{eq_s_boldmu}
s_{\boldsymbol{\mu'}}:=\prod_{j=1}^k \prod_{i=1}^{l_j} s_{ {\mu^{j,i}}' }[X_j].
\end{equation}

\begin{theorem}\label{th_poincare_quiver}
Let $\mathcal{O}_{\boldsymbol{\mu},\boldsymbol{\sigma}}$ be a generic $k$-tuple of adjoint orbits. The Poincaré polynomial for compactly supported intersection cohomology  of $\qcmus$ is
\[
P_c\left(\qcmus,v\right) = v^{d_{\boldsymbol{\mu}}}\left\langle  s_{\boldsymbol{\mu'}}, \mathbb{H}_n^{H L V}(0,v)  \right\rangle.
\]
\end{theorem}
\begin{proof}
For semisimple adjoint orbits, the variety is smooth, intersection cohomology coincides with usual cohomology and the theorem is proved by Hausel, Letellier and Rodriguez-Villegas \cite{aha}. The general case is proved by Letellier \cite{letellier_gl}.
\end{proof}

Hausel--Letellier--Rodriguez-Villegas \cite{aha} proposed a conjecture for the mixed Hodge polynomial of character varieties with semisimple monodromies. It was generalized by Letellier \cite{zclosure} to monodromies with any Jordan type.
\begin{conjecture}[Letellier \cite{zclosure}, Conjecture 1.5]\label{conjecture_ic_char}
Let $\mathcal{C}_{\boldsymbol{\mu},\boldsymbol{\sigma}}$ be a generic $k$-tuple of conjugacy classes, the mixed Hodge polynomial (see Definition \ref{def_polynomial_coho}) of the character variety $\mathcal{M}_{\overline{\mathcal{C}}_{\boldsymbol{\mu},\boldsymbol{\sigma}}}$ is
\[
IH_c (  \mathcal{M}_{\overline{\mathcal{C}}_{\boldsymbol{\mu},\boldsymbol{\sigma}}}; q, v) = \left( v \sqrt{q} \right)^{d_{\boldsymbol{\mu}}} \left\langle s_{\boldsymbol{\mu'}}, \mathbb{H}_n^{H L V}\left(\frac{-1}{\sqrt{q}},v\sqrt{q}\right) \right\rangle
\]
with $q=x y$.
\end{conjecture}

\section{Construction in terms of Nakajima's quiver varieties}\label{sect_Nakajima_quiver_varieties}
In order to study monodromic Weyl group action on the cohomology of the varieties $\qcmus$ and $\qlps$, they need to be put in a family for varying eigenvalues $\boldsymbol{\sigma}$. One way to do that is to construct them as Nakajima's quiver varieties \cite{nakajima1994}. The family obtained is a fibration by the moment map. In this section we recall the construction of $\qcmus$ as a comet-shaped quiver variety. In genus zero the construction is due to Crawley-Boevey \cite{crawley-boevey2003}, for any genus it is due to Hausel, Letellier, Rodriguez-Villegas \cite{aha} and Letellier \cite{letellier_gl}.
\subsection{Generalities about Nakajima's quiver varieties}
\label{subsect_generalities}
In this section we recall the construction of Nakajima's quiver varieties \cite{nakajima1994} in order to fix the notations. Consider a quiver $\Gamma$ with a set of vertices $\Omega_0$ and a set of edges $\Omega_1$. We denote by $t(\gamma)$ the tail and by $h(\gamma)$ the head of an edge $\gamma \in \Omega_1$. A dimension vector for $\Gamma$ is an element $v\in \mathbb{Z}_{\ge 0}^{\Omega_0}$. The space of quiver representations with dimension vector $v$ is identified with a space of matrices,
\[
\Rep\left(\Gamma,v\right):=\bigoplus_{\gamma\in\Omega_1}\Mat_{\KK}(v_{h(\gamma)},v_{t(\gamma)}).
\]
Its cotangent bundle $T^*\Rep\left(\Gamma,v\right)$ can be identified with the space of representations of an extended quiver $\widetilde{\Gamma}$. This extended quiver has the same set of vertices as $\Gamma$. It is obtained by adding an inverse $\overline{\gamma}$ to each edge $\gamma \in \Omega_1$:
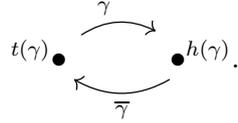
\begin{equation*}
    \begin{tikzcd}
    {}^{t(\gamma)}\bullet \arrow[r,bend left, "\gamma"] &\arrow[l,bend left, "\overline{\gamma}"] \bullet^ {h(\gamma)}.
    \end{tikzcd}
\end{equation*}
We denote by $\overline{\Omega}_1$ the set of such inverted edges, then the set of edges of $\widetilde{\Gamma}$ is $\widetilde{\Omega}:=\Omega_1\sqcup \overline{\Omega}_1$. We have $T^*\Rep\left(\Gamma,v\right)\cong \Rep\left(\widetilde{\Gamma},v\right)$.

For a dimension vector $v\in \ZZ_{\ge 0}^{\Omega_0}$ consider the reductive group $\GL_v:= \prod_{i\in\Omega_0} \GL_{v_i}$. This group acts on $\Rep\left(\widetilde{\Gamma},v\right)$ by
\[
\left(g_{v_i}\right)_{i\in \Omega_0}.(\phi_{\gamma})_{\gamma \in \widetilde{\Omega}} := \left(g_{h(\gamma)} x_{\gamma} g_{t(\gamma)}\inv \right). 
\]
The diagonal embedding of the multiplicative group $\KK^*$ in $\GL_v$ acts trivially so that the action goes down to an action of the group $\G_v := \GL_v / \KK^*$.  The Lie algebra of $\GL_v$ is
\[
\gl_v := \bigoplus_{i\in \Omega_0} \gl_{v_i},
\]
and the Lie algebra of $\G_v$ is
\[
\gv = \left\lbrace (x_j)_{j\in {\Omega_0}}\in \gl_v \left|\sum_{j\in {\Omega_0}} \tr x_j = 0\right. \right\rbrace.
\]
The center of this Lie algebra is given by
\[
 Z(\gv)=\left\lbrace (\xi_j \Id_{v_j})_{j\in {\Omega_0}} \left|(\xi_j)_{j\in\Omega_0}\in {\KK}^{\Omega_0} \text{ with } \sum_{j\in\Omega_0} v_j \xi_j =0 \right. \right\rbrace.
\]
Consider an element $\theta\in\ZZ^{\Omega_0}$ such that $\sum_{i \in \Omega_0} \theta_i v_i = 0$. Such an element is called a stability parameter, it defines a character $\chi_{\theta}$ of the group $\G_v$ by
\begin{equation} \label{chi_theta_quiver}
    \chi_\theta\left((g_j)_{j\in\Omega_0}\right)=\prod_{j\in\Omega_0}\det(g_j)^{-\theta_j}.
\end{equation}
We denote by $\Rep\left(\widetilde{\Gamma},v\right)^{\theta\text{-ss}}$ and by $\Rep\left(\widetilde{\Gamma},v\right)^{\theta\text{-s}}$, the $\theta$-semistable, respectively the $\theta$-stable locus, in the sense of GIT, for the linearization $\chi_{\theta}$. Consider the moment map 
\[
   \map{\mu}{\Rep(\widetilde{\Gamma},v)}{\gv}{\left(\phi_{\gamma}\right)_{\gamma\in\widetilde{\Omega}}}{\sum_{\gamma\in\widetilde{\Omega}}\epsilon(\gamma)\phi_{\gamma}\phi_{\overline{\gamma}}  }
\]
where $\epsilon(\gamma)=1$ for $\gamma \in \Omega_1$ and $\epsilon(\gamma)=-1$ for $\gamma\in \overline{\Omega}_1$.
\begin{definition}[Nakajima's quiver variety]\label{def_quiver_varieties}
Let $\xi$ be an element in $Z(\gv)$, the center of the Lie algebra of $\G_v$. The Nakajima's quiver variety $\mathcal{M}_v^\theta(\xi)$ is defined as the GIT quotient
\begin{equation*}
    \mathcal{M}_v^\theta(\xi):=\mu\inv(\xi)\cap\Rep\left(\widetilde{\Gamma},v\right)^{\theta\text{-ss}}\git \G_v.
\end{equation*}
\end{definition}

We also need another kind of quiver varieties, the Nakajima's framed quiver varieties. Fix a second dimension vector $w\in\NN^
{\Omega_0}$ and define
\[
\Rep\left(v,w\right):=\bigoplus_{j\in {\Omega_0}}\Mat_{\KK}(v_i,w_i),
\]
\[
\Rep\left(w,v\right):=\bigoplus_{j\in {\Omega_0}}\Mat_{\KK}(w_i,v_i).
\]
An element $g\in \GL_v$ acts on $a=(a_j)_{j\in {\Omega_0}}\in \Rep\left(v,w\right) $ by 
\[
g.a:=( a_j g_j \inv)_{j\in {\Omega_0}}
\]
and it acts on $b=(b_j)_{j\in {\Omega_0}}\in \Rep\left(v,w\right) $ by 
\[
g.b:=(g_j b_j )_{j\in {\Omega_0}}.
\]
Introduce the space of framed quiver representations
\[
\Rep\left(\widetilde{\Gamma},v,w\right):=\Rep\left(v,w\right)\oplus\Rep\left(w,v\right)\oplus \Rep\left(\widetilde{\Gamma},v\right).
\]
In this context the moment map is
\[
\begin{array}{ccccc}
     \mu'& : &\Rep(\widetilde{\Gamma},v,w) & \to & \gl_v  \\
     & & (a, b, \phi) & \mapsto & (\mu(\phi)_j-b_j a_j)_{j\in {\Omega_0}}
\end{array}.
\]
For $\theta\in\ZZ^{\Omega_0}$, consider the linearization 
\[
\begin{array}{ccccc}
     \chi_\theta& : & \GL_v & \to & \KK^*  \\
     & & (g_i)_{i\in \Omega_0} &\mapsto & \prod_{i\in \Omega_0} \det (g_i) ^{-\theta_i}
\end{array}.
\]
\begin{definition}[Nakajima's framed quiver varieties]\label{def_framed_quiver}
For $\xi$ in the center of $\gl_v$ and for $\theta\in\ZZ^{\Omega_0}$, the Nakajima's framed quiver variety $\mathcal{M}_{v,w}^\theta(\xi)$ is defined as a GIT quotient with respect to the linearization $\chi_\theta$,
\[
\mathcal{M}_{v,w}^\theta(\xi):={\mu'}\inv(\xi)\cap\Rep\left(\widetilde{\Gamma},v,w\right)^{\theta\text{-ss}}\git\GL_v.
\]
\end{definition}

\subsection{Resolutions of  Zariski closures of adjoint orbits as Nakajima's framed quiver varieties}

In this section we recall the construction of resolutions of closures of adjoint orbits as Nakajima's framed quiver varieties. Those results come from Kraft-Procesi \cite{Kraft_Procesi_80}, Nakajima \cite{Nakajima_1998,nakajima_01}, Crawley-Boevey \cite{Crawley-Boevey_parabolic, crawley-boevey2003}, Shmelkin \cite{shmelkin2009} and Letellier \cite{letellier_gl}. In this subsection and in \ref{subsect_comet_quiver}, we fix the base field $\KK=\C$.

Let $\mathcal{O}_{\underline{\mu},\sigma}$ be an adjoint orbit with semisimple part $\sigma$ and Jordan type $\underline{\mu}\in \Part_{\nu}$ as in \ref{subsect_notations_adjoint_conj}. Consider the resolution $\mathbb{\widetilde{Y}}_{L,P,\sigma}\to \overline{\mathcal{O}}_{\underline{\mu},\sigma}$ as in \ref{prop_resolu_adjoint}. There is a Nakajima's framed quiver variety realizing this resolution. Let $d:=\sum_{i=1}^l \mu^i_1$ and recall that
\[
L\cong \prod_{i=1}^l \prod_{r=1}^{\mu^i_1} \GL_{{\mu^i_r}'}.
\]
The indices $\left({\mu^i_r}'\right)_{ \substack{1\le i\le l \\ 1 \le r \le \mu^i_1  }}$ are relabelled $\left( c_s\right)_{1\le s \le d}$ so that
\[
L\cong \prod_{s=1}^d \GL_{c_s}.
\]
Introduce the parameter $\zeta=\left(\zeta_s\right)_{1\le s\le d}$ such that $\zeta_s=\sigma_i$ if $c_s$ corresponds to ${\mu^i_r}'$ for some $r$.
Consider the quiver $\Gamma_{\mathcal{O}_{\underline{\mu},\sigma}}$ of type $A_{d-1}$ with summit indexed by integers between $1$ and $d-1$ and arrows going in the decreasing direction. Introduce the dimension vector $v_{\mathcal{O}_{\underline{\mu},\sigma}}:=(v_1,...,v_{d-1})$ with 
\[
v_1:=n-c_1, \; \; v_i:=v_{i-1}-c_i \; \text{ for } i>1.
\]
and  $w:=(n,0,\dots,0)$. Define the parameter $\xi_{\mathcal{O}_{\underline{\mu},\sigma}}=(\xi_1,...,\xi_{d-1})$ by $\xi_i = \zeta_i - \zeta_{i+1}$ so that
\begin{equation*}
\xi_i:=\left\lbrace
\begin{array}{ccc}
    \sigma_{k}-\sigma_{k+1} & \text{ if }& i= \mu^{1}_1+\cdots+\mu^{k}_1 \\
     0 & \text{ otherwise } &
\end{array} \right. .
\end{equation*}
The parameter $\xi_{\mathcal{O}_{\underline{\mu},\sigma}}$ will be identified with the element $(\xi_j\Id_{v_j})_{1\le j\le d-1}$.

 We summarize everything in the following diagram showing the quiver, the dimension vector, the parameter $\zeta$  and the parameter $\xi$.
\[
\begin{tikzcd}
 \bullet^1  &\arrow[l] \bullet^2  & \arrow[l] \cdots& \arrow[l] \bullet^{\mu^{1}_1+\cdots+\mu^{k}_1} & \arrow[l] \cdots & \arrow[l] \bullet^{d-1} \\
n-c_1 & n-c_1-c_2 &\cdots & n-\nu_1-\cdots-\nu_k &\cdots & c_r
 \\
 \sigma_1 & \sigma_1 & \cdots & \sigma_k & \cdots & \sigma_r \\
 0 & 0 &\cdots & \sigma_{k}-\sigma_{k+1} &\cdots & 0
%  0 & 0 &\cdots & \theta_{\mu^{(1)}_1+\cdots+\mu^{(k)}_1}>0 &\cdots & 0
\end{tikzcd}
\]
\begin{remark}
When writing the dimension vector under the quiver, we used the fact that $|\mu^i|=\nu_i$.
% sum of the number of box in the young diagram of $\mu^{(i)}$ is equal to the multiplicity $m_i$ of the eigenvalue $\alpha_i$.
\end{remark}

Consider a second dimension vector $w=(n,0,\dots,0)$ and an extended representation $(a,b,\phi)\in\Rep\left(\widetilde{\Gamma}_{\mathcal{O}_{\underline{\mu},\sigma}} ,v_{\mathcal{O}_{\underline{\mu},\sigma}},w\right)$. As $w_i=0$ unless $i=1$, the component $a$ is just a linear map $a : V_1\to W_1$ and $b : W_1 \to V_1$ with $W_1=\C^n$. For $1\le i \le d-2$, denote by $\phi_{i+1,i}$ the linear map associated to the edge from $i+1$ to $i$ and by $\phi_{i,i+1}$ the map associated to the reverse edge from $i$ to $i+1$. Such a representation belongs to ${\mu'}\inv(\xi_{\mathcal{O}_{\underline{\mu},\sigma}})$ if and only if
\begin{equation}\label{prepro}
\left\lbrace
\begin{array}{cccc}
\phi_{2,1}\phi_{1,2}-b a&=&\left(\zeta_1-\zeta_2\right)\Id_{v_1}& \\
\phi_{i+1,i}\phi_{i,i+1}-\phi_{i-1,i}\phi_{i,i-1}&=&\left(\zeta_{i}-\zeta_{i+1}\right)\Id_{v_i} & \text{ for }2\le i \le d-2\\
-\phi_{d-1,d-2}\phi_{d-1,d-2}&=&\left(\zeta_{d-1}-\zeta_{d}\right)\Id_{v_{d-1}} &
\end{array}\right. .
\end{equation}
Those equations are called the preprojective relations.
\begin{example}
The adjoint orbit of the matrix
\[
\begin{pmatrix} 
 \sigma_1 & 1 & 0 & 0 & 0 & 0 \\
 0 & \sigma_1 & 1 & 0 & 0 & 0 \\
 0 & 0 & \sigma_1 & 0 & 0 & 0 \\
 0 & 0 & 0 & \sigma_1 & 0 & 0 \\
 0 & 0 & 0 & 0 & \sigma_2 & 0 \\
 0 & 0 & 0 & 0 & 0 & \sigma_2
\end{pmatrix}
\]
has Jordan type $\underline{\mu}=\left((3,1),(1,1)\right)\in \Part_4 \times \Part_2 $ and we obtain
\[
\begin{tikzcd}
& W_1\arrow[d,bend right, swap, "b"]
\\
 & V_1 \arrow[u,bend right, swap, "a"] \arrow[r,bend right, swap, "\phi_{1,2}"] &\arrow[l,bend right, swap, "\phi_{2,1}"] V_2\arrow[r,bend right, swap, "\phi_{2,3}"]   & \arrow[l,bend right, swap, "\phi_{3,2}"] V_3\\
v_{{\mathcal{O}_{\underline{\mu},\sigma}}} : & 4 & 3 & 2 \\
\zeta : & \sigma_1 & \sigma_1 & \sigma_1 \\
 \xi_{{\mathcal{O}_{\underline{\mu},\sigma}}} : & 0 & 0 & \sigma_1-\sigma_2 .
\end{tikzcd}
\]
\end{example}
\begin{theorem}\label{theorem_resolution_adjoint}
First consider the Nakajima's framed quiver variety $\mathcal{M}^0_{v_{\mathcal O},w}(\xi_{{\mathcal{O}_{\underline{\mu},\sigma}}})$ obtained from the previous data and stability parameter $\theta=0$. The following map is well-defined and is an isomorphism
%The following map is well-defined and is a bijection (it is an isomorphism when $\KK=\C$)
\begin{equation*}
    \map{\Psi_0}{\mathcal{M}^0_{v_{\mathcal{O}_{\underline{\mu},\sigma}},w}\left(\xi_{{\mathcal{O}_{\underline{\mu},\sigma}}}\right)}{\overline{\mathcal{O}}_{\underline{\mu},\sigma}}{(a,b,\phi)}{a b-\sigma_1 \Id_n}.
\end{equation*}
Now take a stability parameter $\theta\in\ZZ_{>0}^{d-1}$, the following map is an isomorphism
%the following map is a bijection (an isomorphism when $\KK=\C$).
\[
\map{\Psi_\theta }{\mathcal{M}^{\theta}_{v_{\mathcal{O}_{\underline{\mu},\sigma}},w}\left(\xi_{{\mathcal{O}_{\underline{\mu},\sigma}}}\right)}{\widetilde{\mathbb{Y}}_{L,P,\sigma}}{(a,b,\phi)}{\left(ab+\sigma_1 \Id_{n}, f_{a,b,\phi} \right)},
\]
where $f_{a,b,\phi}$ is the flag $0\subset E_{d-1}\subset \dots \subset E_1\subset \C^n$ defined by
\[
\begin{array}{cccc}
E_1&:=& \Imm(a), & \\
E_i&:=& \Imm(a\circ \phi_{2,1}\circ \phi_{3,2}\circ  \cdots \circ \phi_{i,i-1}) & \text{ for } 2\le i \le d-1.
\end{array}
\]
Moreover, the following diagram commutes
\begin{equation*}
\begin{tikzcd}
\mathcal{M}^{\theta}_{v_{\mathcal{O}_{\underline{\mu},\sigma}},w}\left(\xi_{{\mathcal{O}_{\underline{\mu},\sigma}}}\right)\arrow{r}{\Psi_\theta} \arrow[swap]{d}{\pi} & \widetilde{\mathbb{Y}}_{L,P,\sigma}\arrow{d}{p^{\sigma}} \\
\mathcal{M}^0_{v_{\mathcal{O}_{\underline{\mu},\sigma}},w }\left(\xi_{{\mathcal{O}_{\underline{\mu},\sigma}}}\right)\arrow[swap]{r}{\Psi_0} & \overline{\mathcal{O}}_{\underline{\mu},\sigma},
\end{tikzcd}
\end{equation*}
where $p^{\sigma}$ is  the resolution of $\overline{\mathcal{O}}_{\underline{\mu},\sigma} $  from Proposition \ref{prop_resolu_adjoint} and $\pi$ is the natural map from GIT theory.

\end{theorem}
\subsection{Comet-shaped quiver varieties}\label{subsect_comet_quiver}
As in the previous subsection, we fix the base field $\KK=\C$. Let $\mathcal{O}_{\boldsymbol{\mu},\boldsymbol{\sigma}}=\left(\mathcal{O}_{\underline{\mu}^1,\sigma^1},\dots,\mathcal{O}_{\underline{\mu}^k,\sigma^k}\right)$ be a generic $k$-tuple of adjoint orbits in $\gl_n$. We recall Crawley-Boevey's result relating the variety $\mathcal{ Q }_{\overline{\mathcal{O}}_{\boldsymbol{\mu},\boldsymbol{\sigma}}}$ to a quiver variety. The idea is to glue together $k$ quivers of type $A$ corresponding to each adjoint orbit $\mathcal{O}_{\underline{\mu}^j,\sigma^j}$ to a central vertex $0$ and to add $g$ loops to this central vertex. We obtain the following comet-shaped quiver $\Gamma_{\mathcal{O}_{\boldsymbol{\mu},\boldsymbol{\sigma}}}$,
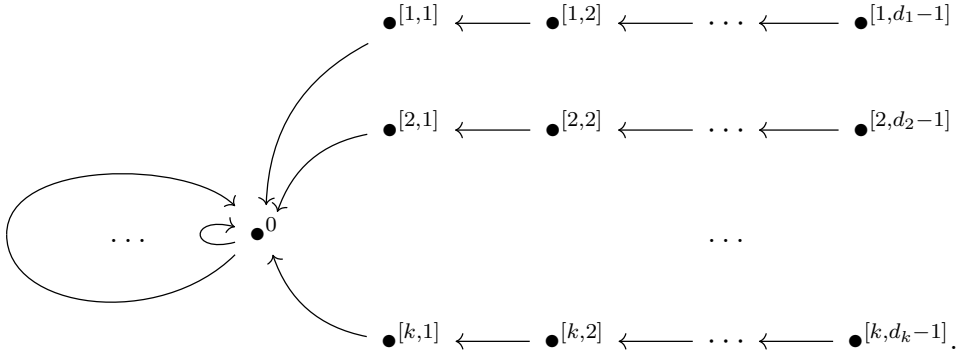
\[
\begin{tikzcd}
 & {} & {} &\arrow[bend right]{ldd} \bullet^{\left[1,1\right]} &\arrow{l} \bullet^{\left[1,2\right]}  & \arrow[l] \cdots& \arrow[l] \bullet^{\left[1,d_1-1\right]} \\
& {} & {} &\arrow[bend right,end anchor={[yshift=-0.5ex]}]{ld}\bullet^{\left[2,1\right]}  &\arrow[l] \bullet^{\left[2,2\right]}  & \arrow[l] \cdots& \arrow[l] \bullet^{\left[2,d_2-1\right]} \\
{} \arrow[rr, start anchor = center, controls={+(0,1) and +(-1,1) }] \arrow[rr, start anchor = center, end anchor = south west, no head, controls={+(0,-1) and +(-1,-1) }] &\dots& \bullet^0 \arrow[loop left] & {} & {} & \dots & {} \\
 & {} & {} &\arrow[bend left]{lu} \bullet^{\left[k,1\right]} &\arrow{l} \bullet^{\left[k,2\right]}  & \arrow[l] \cdots& \arrow[l] \bullet^{\left[k,d_k-1\right]} .
\end{tikzcd}
\]
The $j$-th leg is a quiver of type $A$ with vertices labelled from $[j,1]$ to $[j,d_j-1]$. The dimension vector $v_{\omus}$ is defined such that its coordinate at the central vertex is $n$ and its coordinates on the $j$-th leg coincide with the dimension vector $v_{\mathcal{O}_{\underline{\mu}^j,\sigma^j}}$ described in the previous section.
% \begin{equation*}
%     v_{\omus ,[j,i]}:=v_{\mathcal{O}_{{\underline{\mu}^j,\sigma^j}},i}
% \end{equation*}
% For $1\le i \le k$ and $1\le j \le s_i$. 
Similarly, the parameter $\xi_{\omus}$ is defined such that its coordinates on the $j$-th leg coincide with the parameter $\xi_{ \mathcal{O}_{\underline{\mu}^j,\sigma^j} }$.
%\begin{equation*}
 %   \xi_{\omus,[i,j]}:=\xi_{\mathcal{O}_i,j}
%\end{equation*}
The component at the central vertex $\xi_{\omus,0}$ is defined such that $ v_{\omus}.\xi_{\omus}=0$ hence
\begin{equation*}
   n \xi_{\omus,0}=-\sum_{j=1}^k \sum_{i=1}^{d_j- 1}  v_{\omus,[j,i]}  \xi_{\omus,[j,i]}.
\end{equation*}
Consider a representation of the extended quiver $\phi\in\Rep\left(\widetilde{\Gamma}_{\omus},v_{\omus}\right)$.
\begin{itemize}
    \item Denote by $\phi_{[j,i]}$ the linear map associated to the arrow with tail $[j,i]$ and $\phi_{\overline{[j,i]}}$ the linear map associated to the reversed arrow with head $[j,i]$.
    \item For $1\le i\le g$ the map associated to the $i$-th loop is denoted $\phi_i$ and the one associated to the reverse loop is denoted $\phi_{\overline{i}}$.
\end{itemize} 
As usual $\mu$ is the moment map and $\xi_{\omus}$ is identified with an element in the center of the Lie algebra $\mathfrak{g}_{v_{\omus}}$ and we let
\begin{equation*}
    X_j:=\phi_{[j,1]}\phi_{\overline{[j,1]}}-\zeta_{[j,1]}.
\end{equation*}
If $\phi$ belongs to $\mu\inv(\xi_{\omus})$ then $X_j\in\overline{\mathcal{O}}_{\underline{\mu}^j,\sigma^j}$. Indeed it follows from the previous description of closures of adjoint orbits as framed quiver varieties and the identification, for each leg, of the vector space at the central vertex with the framing vector space $W_1$ from the previous section.

Now if $A_i$ is the linear map associated to the $i$-th loop of the quiver and $B_i$ the map associated to the reversed loop, the preprojective relation at the central vertex is exactly the equation defining $\mathcal{V}_{\omusbar}$. Then the following map is well-defined
\begin{equation*}
    \map{\Psi^{\omus}}{\mu\inv(\xi_{\omus}) }{\mathcal{V}_{\omusbar}}{\phi}{\left(A_1,B_1,\dots,A_g,B_g,X_1,\dots,X_k\right)}.
\end{equation*}
\commentaires{\gl_n^{2 g}\times \overline{\mathcal{O}}_1\times\dots\times\overline{\mathcal{O}}_k}
\begin{theorem}
In the following diagram where the vertical arrows are  quotient maps, the application $\Psi^{\omus}$ goes down to the quotient to an isomorphism $\Phi^{\omus}$,
%a bijective morphism $\Phi^{\omus}$ (when $\KK=\C$ it is an isomorphism).
\begin{equation*}
 \begin{tikzcd}
 \mu\inv(\xi_{\omus}) \arrow{r}{\Psi^{\omus}}\arrow{d}{}& \mathcal{V}_{\omusbar}\arrow{d}{}\\
  \mathcal{M}^0_{v_{\omus}}(\xi_{\omus})\arrow{r}[swap]{\Phi^{\omus}}&\mathcal{ Q }_{\omusbar}.
 \end{tikzcd}
\end{equation*}
\end{theorem}
\begin{proof}
It is proved by Crawley-Boevey \cite{crawley-boevey_2001, crawley-boevey2003}, see also Letellier \cite[Proposition 5.2.2]{letellier_gl}  for any genus.
\end{proof}

The resolution $\qlps$ of $\qcmus$, as introduced in \ref{def_resol_char}, is also interpreted as a Nakajima's quiver variety for the quiver $\Gamma_{\omus}$.
\begin{theorem}\label{th_quiver_resol}
Consider a stability parameter $\theta$ associated to the quiver $\mathcal{ Q }_{\omusbar}$ such that $\theta_{[j,i]}>0$ for each vertex $[j,i]$. There is an isomorphism $\Phi^{\omus,\theta}:\mathcal{M}^{\theta}_{v_{\omus}}(\xi_{\omus})\to \qlps $ and the following diagram commutes,

\begin{equation*}
     \begin{tikzcd}
\mathcal{M}^{\theta}_{v_{\omus}}(\xi_{\omus}) \arrow{r}{\Phi^{\omus,\theta}}\arrow[swap]{d}{\pi}& \qlps \arrow{d}{p^{\sigma}}\\
  \mathcal{M}^0_{v_{\omus}}(\xi_{\omus})\arrow{r}[swap]{\Phi^{\omus}}&\mathcal{ Q }_{\omusbar},
 \end{tikzcd}
\end{equation*}
where $\pi$ is the natural projection from GIT theory.
\end{theorem}
\begin{proof}
The map $\Phi^{\omus,\theta}$ is constructed by Letellier \cite[Section 5.3]{letellier_gl}.
%It follows from Letellier's article \cite{letellier_gl}, where the construction of the map $\Phi^{\omus,\theta}$ is given in Section 5.3.
This map is induced by the map $\Psi_\theta$ of Theorem \ref{theorem_resolution_adjoint}. Contrarily to Letellier's article,  we do not consider partial resolutions so that our parameter $\theta$ only has non-zero components. Therefore the dimension vector for the quiver variety $\mathcal{M}^{\theta}_{v_{\omus}}(\xi_{\omus})$ describing the resolution $\qlps $ is the same as the dimension vector of the quiver variety describing $\mathcal{ Q }_{\omusbar}$. 
\end{proof}

The quiver variety point of view gives a criteria for non-emptiness. The question of emptiness of $\qcmus$ (and its analogous character variety) is known as the Deligne-Simpson problem. See Kostov \cite{Kostov_ds} for a survey about this problem. Crawley-Boevey gave a solution to the problem in the generic case in terms of roots of quivers \cite{crawley-boevey2003}, see also Letellier \cite[Section 5.2]{letellier_gl}. Those results are summarized in the following theorem.
\begin{theorem}
Let $\mathcal{O}_{\boldsymbol{\mu},\boldsymbol{\sigma}}$ be a generic $k$-tuple of adjoint orbits. The variety $\qcmus$ is non-empty if and only  $\mathcal{Q}_{\mathcal{O}_{\boldsymbol{\mu},\boldsymbol{\sigma}}}$ is not empty. This happens if and only if the dimension vector $v_{\mathcal{O}_{\boldsymbol{\mu},\boldsymbol{\sigma}}}$ is a root of the quiver $\Gamma_{\mathcal{O}_{\boldsymbol{\mu},\boldsymbol{\sigma}} }$. This is always the case for $g>0$.
\end{theorem}

\subsection{Family of comet-shaped quiver varieties}\label{subsect_family_qlpt}
Now the field $\KK$ is again either $\C$ or $\Fqbar$. When the eigenvalues $\boldsymbol{\sigma}$ are varying, the varieties $\qlps$ fit in a family. First we describe explicitly this family, then we give an interpretation in terms of Nakajima's quiver varieties and moment map.
% is described explicitly, before being interpreted in terms of Nakajima's quiver varieties and moment map.
\begin{notations}
From now on the pair $\mathbf{L},\mathbf{P}$ is fixed. For short, let
\[
Z(\boldsymbol{l}):=Z(\mathfrak{l}_1)\times\dots\times Z(\mathfrak{l}_k).
\]
Denote by $\mathcal{B}$ the subset of elements $\boldsymbol{\sigma}\in Z(\boldsymbol{l})$ such that the $k$-tuple of adjoint orbits $\omus$ is generic. Note that the genericity condition depends only the semisimple part $\boldsymbol{\sigma}$ and not on the type $\boldsymbol{\mu}$. 
%We assume that $\boldsymbol{L}$ is chosen so that $\mathcal{B}$ is not empty. 
Then $\mathcal{B}$ is either empty or Zariski open in the hyperplane of $Z(\boldsymbol{l})$ defined by the vanishing of the sum of the traces. %Identifying  $Z(\boldsymbol{l})$ with an affine space, $\mathcal{B}$ is either empty or the complementary of a finite union of hyperplanes in the codimension one subspace.
\end{notations}
\begin{definition}[Family of varieties $\qlps$]
Define
\begin{multline*}
    \vlp:=\left\lbrace \left(\boldsymbol{\sigma}, \; (A_i,B_i)_{1\le i \le g}, \; (X_j,g_j P^j)_{1\le j \le k} \right) \; \right| \\ \left. \boldsymbol{\sigma}\in \mathcal{B}, \text{ and } (A_i,B_i)_{1\le i \le g}, \; (X_j,g_j P^j)_{1\le j \le k}\in \VLPs    \right\rbrace,
\end{multline*}
\begin{equation*}
    \qlp:= \vlp \git \GL_n.
\end{equation*}
Denote by $\eta$ the natural map $\eta : \qlp \to \mathcal{B}$. The varieties $\qlps = \eta\inv (\boldsymbol{\sigma}) $ fit in a family $\qlp$ over $\mathcal{B}$.
\end{definition}
The choice of $\boldsymbol{L}$ determines a unique quiver $\Gamma_{\omus}$ and a unique dimension vector $v_{\omus}$ independent of the choice of $\boldsymbol{\sigma}$. Assume that the dimension vector is indivisible so that $\mathcal{B}$ is not empty. Then we can make the following assumption:
\begin{assumption}[Genericity of the stability parameter $\theta$]
The stability parameter $\theta$ is generic, \textit{i.e.}, it is a stability parameter for the quiver $\Gamma_{\omus}$ with dimension vector $v_{\omus}$ such that $\theta.v_{\omus}=0$ and $\theta.v \ne 0$ for a smaller dimension vector $v$.
%$(\theta,0,0)\in H_{v_{\omus}}^{\reg}$ with notations from \ref{def_hk_reg_locus}.
\end{assumption} 
The construction of  Theorem \ref{th_quiver_resol} extends to this family. It provides the following commutative diagram (the left vertical arrows is induced by the moment map $\mu$)
\begin{equation}\label{diagram_family_varieties}
     \begin{tikzcd}
\mu\inv(\mathfrak{z}_{v_{\omus}}^{\gen})^{\theta\text{-ss}}\git G_{v_{\omus}} \arrow{r}{\Phi}\arrow[swap]{d}{}& \qlp\arrow{d}{\eta}\\
\mathfrak{z}_{v_{\omus}}^{\gen}\arrow{r}[swap]{}& \mathcal{B},
\end{tikzcd}
\end{equation}
where $\theta$ is a fixed generic stability parameter and $\mathfrak{z}_{v_{\omus}}^{\gen}$ is the subset of the center of the Lie algebra $\mathfrak{g}_{v_{\omus}}$ corresponding to the subset $\mathcal{B}$ under the correspondence between the parameters $\xi_{\omus}$ and the eigenvalues $\boldsymbol{\sigma}$. Note that the correspondence between parameters of the quiver variety $\xi_{\omus}\in Z(\mathfrak{g}_{v_{\omus}})$ and $Z(\boldsymbol{l})$ is not bijective, only difference of successive eigenvalues appear in the construction of the quiver variety. Thus the previous diagram relies on a choice of $k-1$ eigenvalues.
To $\boldsymbol{\sigma}\in Z(\boldsymbol{l})$ associate the element $(\xi_{\omus},\sigma^1_1,\dots,\sigma^{k-1}_1)$ in $Z(\mathfrak{g}_{v_{\omus}})\times \KK^{k-1}$ this defines a bijective map
\begin{equation}
   h : Z(\boldsymbol{l}) \xrightarrow[]{\sim} \mathfrak{z}_{v_{\omus}}\times \KK^{k-1}.
\end{equation}
Note that for a given parameter $\xi_{\omus}$ the genericity conditions is independant of the choice of the $k-1$  eigenvalues, namely $h\inv(\xi_{\omus},\sigma^1_1,\dots,\sigma^{k-1}_1)$ is generic if and only if $h\inv(\xi_{\omus},0,\dots,0)$ is generic.
Therefore Diagram
\eqref{diagram_family_varieties} can be modified in order to account for various choices of eigenvalues, then the horizontal arrows are isomorphisms and %isomorphisms when $\KK=\C$.

\begin{equation}\label{diagram_family_varieties_good}
     \begin{tikzcd}
\KK^{k-1}\times\mu\inv(\mathfrak{z}_{v_{\omus}}^{\gen})^{\theta\text{-ss}}\git G_{v_{\omus}} \arrow{r}{\Phi}\arrow[swap]{d}{\Id\times \mu}& \qlp\arrow{d}{\eta}\\
\KK^{k-1}\times\mathfrak{z}_{v_{\omus}}^{\gen}\arrow{r}[swap]{}& \mathcal{B}.
 \end{tikzcd}
\end{equation}
\begin{theorem}\label{th_qlp_locally_constant}
If $\KK=\C$, or if $\KK=\Fqbar$ and the characteristic is large enough, the cohomology sheaves $\mathcal{H}^i \eta_! \kappa$ are constant sheaves.
\end{theorem}
\begin{proof}
When $\KK = \C$ this follows from the quiver variety point of view of Diagram \eqref{diagram_family_varieties_good}, it is a well-known fact used by Nakajima \cite{nakajima1994} to construct a Weyl group action on the cohomology of quiver varieties (see also \cite{ballandras20}).
%When $\KK=\C$, this is a consequence of Chapter \ref{chapter_trivialization} Corollary \ref{coro_locally_constant_hk} and diagram \eqref{diagram_family_varieties_good}. 
To prove the result for $\KK=\Fqbar$ we can change characteristic as in \cite[proof of Theorem 2.3]{kac_conj}.  This implies the result in large enough characteristic.
\end{proof}
\section{Monodromic Weyl group action}\label{sect_weyl_quiver}
The goal of this section is to construct and study the monodromic Weyl group action on the cohomology of the quiver varieties $\qlps$. We use technics from Nakajima \cite{nakajima1994}, Lusztig (see Letellier \cite[Proof of proposition 5.5.3]{letellier_f}), Mellit \cite[Section 8]{mellit_cell} and Hausel--Letellier--Rodriguez-Villegas \cite{kac_conj}. The construction relies essentially on Theorem \ref{th_qlp_locally_constant}. 

\subsection{Family of resolutions of closures of adjoint orbits}\label{subsect_YLP}
%The resolutions of adjoint orbits previously introduced fit in a family. This will be useful to construct relative Weyl group actions when considering family of comet-shaped quiver varieties.

In this section we study a family formed by the varieties $\widetilde{\mathbb{Y}}_{L,P,\sigma}$ when $\sigma$ is varying. It will be usefull in the next section to study a family of comet-shaped quiver varieties and to obtain some dimension estimates to prove Lemma \ref{lemma_c_ic}.

Let $P$ be a parabolic subgroup of $\GL_n$ and let $L$ be a Levi factor of $P$, then $L$ is isomorphic to a group of block diagonal matrices $\GL_{c_1}\times\dots\times\GL_{c_r}$. The Lie algebra of $L$ and of $U_P$ are denoted by $\mathfrak{l}$ respectivly by $\mathfrak{u}_P$. At the level of the Lie algebras the Levi decomposition reads $\mathfrak{p}=\mathfrak{l}\oplus \mathfrak{u}_P$. The center of this Lie algebra $\mathfrak{l}$ is denoted by $Z(\mathfrak{l})$ and its regular locus is
\begin{equation*}
    Z(\mathfrak{l})^{\reg} =\left\lbrace x\in Z(\mathfrak{l}) \left| Z_G(x)= L\right. \right\rbrace.
\end{equation*}
Define
\begin{equation*}
    \widetilde{\mathbb{Y}}^{\reg}_{L,P}=\left\lbrace (x,g L)\in \gl_n\times \GL_n / L\left|g\inv x g \in Z(\mathfrak{l})^{\reg}\right. \right\rbrace.
\end{equation*}
Consider the projection on the first factor $p^{\reg} : \widetilde{\mathbb{Y}}^{\reg}_{L,P} \to \gl_n$, denote by $\mathbb{Y}^{\reg}_{L,P}$ its image. This image consists of semisimple elements with $r$ distinct eigenvalues with multiplicities $c_1,\dots,c_r$. Consider the relative Weyl group $W_{\GL_n}(L)=N_{\GL_n}(L)/L$, and for each $w\in W_{\GL_n}(L)$ chose a representative $ \dot{w}\in N_{\GL_n}(L)$. This relative Weyl group acts on $Z(\mathfrak{l})$ by
\[
w.\sigma := \dot{w} \sigma\dot{w}\inv.
\]
Consider the fiber product
\[
\begin{tikzcd}
 Z(\mathfrak{l})^{\reg} \arrow[d] & \arrow[l] \mathbb{Y}^{\reg}_{L,P} \times_{ Z(\mathfrak{l})^{\reg} / W_{\GL_n}(L)}  Z(\mathfrak{l})^{\reg}   \arrow[d]\\
 Z(\mathfrak{l})^{\reg} / W_{\GL_n}(L) & \arrow[l,"\chi"] \mathbb{Y}^{\reg}_{L,P},
\end{tikzcd}
\]
with $\chi$ the characteristic polynomial.
Note that the following map is an isomorphism
\begin{equation}\label{eq_regular_fiber_prod}
\begin{array}{ccc}
    \widetilde{\mathbb{Y}}^{\reg}_{L,P} & \to &  \mathbb{Y}^{\reg}_{L,P} \times_{ Z(\mathfrak{l})^{\reg} / W_{\GL_n}(L)}  Z(\mathfrak{l})^{\reg}  \\
  (x, g L)   & \mapsto & (x ,g\inv x g) .
\end{array}
\end{equation}
Therefore the $W_{\GL_n}(L)$-action on $Z(\mathfrak{l})^{\reg}$ induces an action on $\widetilde{\mathbb{Y}}^{\reg}_{L,P}$.  It is given explicitly by
\begin{equation*}
    w.(x,g L) = (x,g \dot{w}\inv L).
\end{equation*}
Then the morphism
\begin{equation*}
    \widetilde{\mathbb{Y}}^{\reg}_{L,P} \xrightarrow{p^{\reg}} \mathbb{Y}^{\reg}_{L,P}
\end{equation*}
is a Galois cover with group $W_{\GL_n}(L)$. This relative Weyl group acts on the push forward of the constant sheaf $p^{\reg}_{*}\kappa$. Define
\begin{equation*}
    \widetilde{\mathbb{Y}}_{L,P}=\left\lbrace (x,g P)\in \gl_n\times \GL_n/P \left|g\inv x g \in Z(\mathfrak{l})\oplus \mathfrak{u}_P\right.\right\rbrace.
\end{equation*}

\begin{remark}\label{remark_partial_flag}
An element $g P\in \GL_n /P$ is identified with a partial flag 
\[
0=E_r\subset E_{r-1}\subset \dots \subset E_1 \subset \KK^n
\]
such that $\dim E_{i-1}/E_i = c_i$ for all $1\le i \le r$. Indeed $\GL_n$ acts transitively on such flags and the stabilizer of any of them is isomorphic to $P$. Then a point $(x,g P)$ in $\widetilde{\mathbb{Y}}_{L,P}$ consists of an endomorphism $x\in\gl_n$ and a partial flag $g P$ preserved by $x$ such that $x$ acts as a scalar on $E_{i-1}/E_i$ for all $1\le i \le r$.
\end{remark}

% \begin{remark}
% An element $g P\in \GL_n /P$ identifies with a partial flag $0=E_l\subset E_{l-1}\subset \dots \subset E_1 \subset \KK^n$ such that $\dim E_{i-1}/E_i = c_i$ for all $1\le i \le l$. Then a point $(x,g P)$ in $\widetilde{\mathbb{Y}}_{L,P}$ consists of an endomorphism $x\in\gl_n$ and a partial flag $g P$ preserved by $x$ such that $x$ acts as a scalar on $E_{i-1}/E_i$ for all $1\le i \le l$.
% \end{remark}
Denote by $\mathbb{Y}_{L,P}$ the image of the projection to the first factor $p:\widetilde{\mathbb{Y}}_{L,P}\to\gl_n$. Note that the map $p$ is proper. The following theorem is a particular case of \cite[Lemma 4.3 and Proposition 4.5]{Lusztig1984}. It can be seen as a generalization of Borho--MacPherson \cite{borho} approach to Springer theory.

\begin{theorem}
The subvariety $\mathbb{Y}^{\reg}_{L,P}$ is open, smooth and dense in  $\mathbb{Y}_{L,P}$. The following square is cartesian
\begin{equation}\label{square_g_springer}
    \begin{tikzcd}
    \widetilde{\mathbb{Y}}^{\reg}_{L,P} \arrow[r,"i"]\arrow[d,swap,"p^{\reg}"] & \widetilde{\mathbb{Y}}_{L,P} \arrow[d,"p"] \\
    \mathbb{Y}^{\reg}_{L,P} \arrow[hookrightarrow]{r} & \mathbb{Y}_{L,P},
    \end{tikzcd}
\end{equation}
with $i$ the map $(x,g L)\to (x,g P)$. Moreover $p_! \kappa = \ic{\mathbb{Y}_{L,P}, \; p^{\reg}_! \kappa }$ so that $W_{\GL_n} (L)$ acts on $p_! \kappa$.
\end{theorem}

\begin{remark}\label{remark_dim_YLP}
The morphism $p^{\reg}$ is a Galois cover and $i$ is an open embedding so that the dimension can be easily computed
\begin{equation}
    \dim \mathbb{Y}_{L,P} = \dim \widetilde{\mathbb{Y}}_{L,P} = \dim \widetilde{\mathbb{Y}}^{\reg}_{L,P} = \dim \GL_n - \dim L + \dim Z(L) .
\end{equation}
\end{remark}
Let us describe the relation with the resolutions of closures of adjoint orbits introduced in \ref{prop_resolu_adjoint}. Let $\sigma \in Z(\mathfrak{l})$ and let $M:=Z_{\GL_n}(\sigma)$ be the stabilizer of $\sigma$ in $\GL_n$. The notations from \ref{subsect_resolution_conjugacy} are used so that $M\cong  \GL_{\nu}$ for $\nu$ a partition of $n$. Moreover $L\subset M$ and the integers $(c_1,c_2,\dots,c_r)$ are relabelled $({\mu^1_1}',{\mu^1_2}',\dots)$ so that ${\mu^i}'$ is a partition of $\nu_i$. The inclusion $L\subset M$ comes from the inclusions 
\[
\GL_{{\mu^i_1}'}\times \dots \times \GL_{{\mu^i_{l_i}}'} \subset \GL_{\nu_i}.
\]
The resolution of the closure of $\overline{\mathcal{O}}_{\underline{\mu},\sigma}$ fits in the following diagram
\begin{equation}\label{diag_springer}
    \begin{tikzcd}
 \widetilde{\mathbb{Y}}_{L,P} \arrow[d,"p"] \arrow[r,hookleftarrow]& \widetilde{\mathbb{Y}}_{L,P,\sigma} \arrow[d,"p^{\sigma}"]\\
 \mathbb{Y}_{L,P} \arrow[r,hookleftarrow] & \overline{\mathcal{O}}_{\underline{\mu},\sigma}=\bigsqcup_{\underline{\rho}\preceq\underline{\mu}} \mathcal{O}_{\underline{\rho},\sigma}
    \end{tikzcd}
\end{equation}
The decomposition $\overline{\mathcal{O}}_{\underline{\mu},\sigma}=\bigsqcup_{\underline{\rho}\preceq\underline{\mu}} \mathcal{O}_{\underline{\rho},\sigma}$ actualy comes from a decomposition of $\mathbb{Y}_{L,P}$. Define
\[
\mathbb{Y}_{L,P}^{M,\underline{\rho}}:=\bigsqcup_{\sigma'\in Z(\mathfrak{m})^{\reg}} \mathcal{O}_{\underline{\rho},\sigma'}.
\]
This decomposition is similar to the one introduced by Shoji \cite{Shoji}.
\begin{proposition}\label{prop_decomposition_YLP}
The variety $\mathbb{Y}_{L,P}^{M,\underline{\rho}}$ is smooth of dimension
\[
\dim \mathbb{Y}_{L,P}^{M,\underline{\rho}} = \dim \mathcal{O}_{\underline{\rho},\sigma} + \dim Z(\mathfrak{m}).
\]
The variety $\mathbb{Y}_{L,P}$ admits the following decomposition
\[
\mathbb{Y}_{L,P} = \bigsqcup_{M} \bigsqcup_{\underline{\rho}\preceq \underline{\mu}} \mathbb{Y}_{L,P}^{M,\underline{\rho}}.
\]
The first union is over the set of stabilizers of elements $\sigma \in Z(\mathfrak{l})$. In the second union, $\underline{\mu}$ depends on $M$ as previously described.
The unique part indexed by $M=L$ is $\mathbb{Y}^{\reg}_{L,P}$.

\end{proposition}
\begin{proof}
Denote by $Z_{\underline{\rho}}$ the stabilizer in $\GL_n$ of the element $J_{\underline{\rho},\sigma}$ in $\mathcal{O}_{\underline{\rho},\sigma}$ (see Notations \ref{notations_jmus}). There is a finite cover
\[
\begin{array}{ccc}
    Z(\mathfrak{m})^{\reg}\times \GL_n / Z_{\underline{\rho}} & \to & \mathbb{Y}_{L,P}^{M,\underline{\rho}}  \\
 \left( \sigma', g Z_{\underline{\rho}}\right)  & \mapsto & g J_{\underline{\rho},\sigma'} g\inv , 
\end{array}
\]
hence $\mathbb{Y}_{L,P}^{M,\underline{\rho}} $ is smooth and
\[
\dim \mathbb{Y}_{L,P}^{M,\underline{\rho}} = \dim \mathcal{O}_{\underline{\rho},\sigma} + \dim Z(\mathfrak{m}).
\]
\end{proof}

\subsection{Decomposition of the family \texorpdfstring{$\mathcal{Q}_{\boldsymbol{L},\boldsymbol{P}}$}{ }}
In this section we study a family related to the family $\widetilde{\mathcal{Q}}_{\boldsymbol{L},\boldsymbol{P}}$ introduced briefly in \ref{subsect_family_qlpt}. We compute some dimensions which will be useful to prove Lemma \ref{lemma_c_ic}. 
\begin{notations}
First we recall the notations from \ref{subsect_YLP} in this context. For $1\le j \le k$,
\[
\widetilde{\mathbb{Y}}_{L^j,P^j}:=\left\lbrace (X,g_j P^j) \in \gl_n\times \GL_n /P^j \left| g_j\inv X g_j \in Z(\mathfrak{l}^j)\oplus \mathfrak{u}_{P^j} \right. \right\rbrace 
\]
and define
\[
\ylpt := \widetilde{\mathbb{Y}}_{L^1,P^1}\times\dots \times \widetilde{\mathbb{Y}}_{L^k,P^k}.
\]
Then $\ylp$ is the image in $\gl_n^k$ of the map $p$ forgetting the partial flags $g_j P^j$,
\[
\begin{array}{ccccc}
    p & : & \ylpt & \to & \gl_n^k  \\
   & & (X_j, g_j P^j)_{1\le j\le k}  &\mapsto & (X_j)_{1\le j\le k}. 
\end{array}
\]
Similarly $\mathcal{V}_{\boldsymbol{L},\boldsymbol{P}}$, respectively $\mathcal{Q}_{\boldsymbol{L},\boldsymbol{P}}$, is obtained from $\vlp$, respectively $\qlp$, by forgetting the partial flags.
\end{notations}

In this section a decomposition of the family $\mathcal{Q}_{\boldsymbol{L},\boldsymbol{P}}$ is deduced from the  decomposition $\overline{\mathcal{O}}_{\underline{\mu},\sigma}= \bigsqcup_{\rho \preceq \underline{\mu} } \mathcal{O}_{\underline{\rho},\sigma}$ and from the decomposition introduced in Proposition \ref{prop_decomposition_YLP},
\[
\mathbb{Y}_{L,P} = \bigsqcup_{M} \bigsqcup_{\underline{\rho} \preceq \underline{\mu}} \mathbb{Y}^{M,\underline{\rho}}_{L,P}.
\] 
The decomposition is used in the next section (Lemma \ref{lemma_c_ic}) in order to define a Weyl group action.
%to extend the monodromic action from $\mathcal{B}^{\reg}$ to the whole $\mathcal{B}$.

Let $\mathbb{Y}^{\mathcal{B}}_{\boldsymbol{L},\boldsymbol{P}}$ be the subset of elements in $\mathbb{Y}_{\boldsymbol{L},\boldsymbol{P}}$ with generic semisimple parts, \textit{i.e.}, with a $k$-tuple of semisimple parts belonging to $\mathcal{B}$. The set $\mathbb{Y}^{\mathcal{B}}_{\boldsymbol{L},\boldsymbol{P}}$ is assumed to be non-empty. The dimension of $\mathbb{Y}^{\mathcal{B}}_{\boldsymbol{L},\boldsymbol{P}}$ is computed similarly to  $\dim \mathbb{Y}_{L,P}$ in Remark \ref{remark_dim_YLP},
\[
\dim \mathbb{Y}^{\mathcal{B}}_{\boldsymbol{L},\boldsymbol{P}} = k n^2+ \dim{\mathcal{B}} -\sum_{j=0}^k \dim L^j.
\]
The decomposition $\mathbb{Y}_{L,P}=\bigsqcup_{M}\bigsqcup_{\underline{\rho}\preceq\underline{\mu}} \mathbb{Y}_{L,P}^{M,\underline{\rho}} $ induces a similar decomposition for $\mathbb{Y}^{\mathcal{B}}_{\boldsymbol{L},\boldsymbol{P}}$,
\[
\mathbb{Y}^{\mathcal{B}}_{\boldsymbol{L},\boldsymbol{P}} = \bigsqcup_{\boldsymbol{M}}\bigsqcup_{\boldsymbol{\rho}\preceq\boldsymbol{\mu}} \mathbb{Y}^{\mathcal{B},\boldsymbol{M},\boldsymbol{\rho}}_{\boldsymbol{L},\boldsymbol{P}},
\]
where $\boldsymbol{M}=(M^1,\dots, M^l)$ and $\mathbb{Y}^{\mathcal{B},\boldsymbol{M},\boldsymbol{\rho}}_{\boldsymbol{L},\boldsymbol{P}}$ is the subset of elements in
\[
\mathbb{Y}_{L^1,P^1}^{M^1,\underline{\rho}^1}\times \dots \times \mathbb{Y}_{L^k,P^k}^{M^k,\underline{\rho}^k}
\]
with a generic $k$-tuple of semisimple parts. From the computation of the dimension of $\mathbb{Y}_{L,P}^{M,\underline{\rho}}$ in Proposition \ref{prop_decomposition_YLP}, we deduce that when $Z(\boldsymbol{m})\cap \mathcal{B}$ is not empty
\begin{equation}\label{eq_dim_YLPBMr}
\dim \mathbb{Y}^{\mathcal{B},\boldsymbol{M},\boldsymbol{\rho}}_{\boldsymbol{L},\boldsymbol{P}} = \sum_{j=1}^n \dim \mathcal{O}_{\underline{\rho}^j,\sigma^j} + \dim Z(\boldsymbol{m})\cap \mathcal{B}.
\end{equation}

Now the decomposition of $\mathbb{Y}^{\mathcal{B}}_{\boldsymbol{L},\boldsymbol{P}}$ induces a decomposition of the family of quiver varieties $\mathcal{Q}_{\boldsymbol{L},\boldsymbol{P}}$ and we define
\[
\mathcal{Q}^{\boldsymbol{M},\boldsymbol{\rho}}_{\boldsymbol{L},\boldsymbol{P}}:= \left(\mathcal{V}_{\boldsymbol{L},\boldsymbol{P}}\times_{\mathbb{Y}^{\mathcal{B}}_{\boldsymbol{L},\boldsymbol{P}}} \mathbb{Y}^{\mathcal{B},\boldsymbol{M},\boldsymbol{\rho}}_{\boldsymbol{L},\boldsymbol{P}}\right) \git \PGL_n.
\]
\begin{proposition}\label{prop_decomposition_qlp}
The variety $\mathcal{Q}_{\boldsymbol{L},\boldsymbol{P}}$ admits the following decomposition
\[
\mathcal{Q}_{\boldsymbol{L},\boldsymbol{P}} = \bigsqcup_{\boldsymbol{M}}\bigsqcup_{\boldsymbol{\rho}\preceq\boldsymbol{\mu}} \mathcal{Q}^{\boldsymbol{M},\boldsymbol{\rho}}_{\boldsymbol{L},\boldsymbol{P}}.
\]
When non-empty, the dimension of a part is
\begin{equation}\label{eq_dim_qlpmr}
    \dim \mathcal{Q}^{\boldsymbol{M},\boldsymbol{\rho}}_{\boldsymbol{L},\boldsymbol{P}} = n^2(2g-2)+2 +\dim Z(\boldsymbol{m})\cap\mathcal{B} + \sum_{j=1}^k \dim \mathcal{O}_{\underline{\rho}^j,\sigma^j}.
\end{equation}
\end{proposition}
\begin{proof}
% The decomposition is defined by
% \[
% \mathcal{Q}^{\boldsymbol{M},\boldsymbol{\rho}}_{\boldsymbol{L},\boldsymbol{P}}:= \left(\mathcal{V}_{\boldsymbol{L},\boldsymbol{P}}\times_{\mathbb{Y}_{\boldsymbol{L},\boldsymbol{P}}} \mathbb{Y}^{\mathcal{B},\boldsymbol{M},\boldsymbol{\rho}}_{\boldsymbol{L},\boldsymbol{P}}\right) \git \PGL_n.
% \]
The dimension of $\mathcal{Q}^{\boldsymbol{M},\boldsymbol{\rho}}_{\boldsymbol{L},\boldsymbol{P}}$ can be computed just like the dimension of $\mathcal{Q}_{\mathcal{O}_{\boldsymbol{\mu},\boldsymbol{\sigma}}}$ (see Hausel, Letellier, Rodriguez-Villegas \cite[Theorem 2.2.4]{aha} and Letellier \cite[Corollary 5.2.3]{letellier_gl}). The computation relies on the smoothness of $\mathbb{Y}^{\mathcal{B},\boldsymbol{M},\boldsymbol{\rho}}_{\boldsymbol{L},\boldsymbol{P}}$ which follows from the smoothness of $\mathbb{Y}^{\mathcal{B},M^j,\underline{\rho}^j}_{L^j,P^j}$ and on the expression \eqref{eq_dim_YLPBMr} for the dimension of $\mathbb{Y}^{\mathcal{B},\boldsymbol{M},\boldsymbol{\rho}}_{\boldsymbol{L},\boldsymbol{P}}$.
\end{proof}

\subsection{\texorpdfstring{$W$-equivariant}{W-equivariant} structure on the cohomology of the fibers of the family \texorpdfstring{$\qlp$}{ }}
%\subsection{Construction of a Weyl group action on the cohomology of the quiver varieties in the family \texorpdfstring{$\qlp$}{}}

In this section we use the family $\qlp \to \mathcal{B}$ in order to construct a Weyl group action on the cohomology of the varieties $\qlps$ for $\boldsymbol{\sigma}\in \mathcal{B}$.
%The family $\qlp \to \mathcal{B}$ is used to construct a Weyl group action on the cohomology of the varieties $\qlps$ for $\boldsymbol{\sigma}\in \mathcal{B}$. 
The Weyl group studied in this section is
\[
W:=W_{\GL_n}(L^1)\times\dots\times W_{\GL_n}(L^k).
\]
Each $W_{\GL_n}(L^j)$ is isomorphic to a symmetric group and acts on $Z(\mathfrak{l}^j)$ by permuting the distinct eigenvalues with the same multiplicities. Then the Weyl group $W$ acts on $\mathcal{B}$, for $w=(w_1,\dots,w_k)\in W$ and $\boldsymbol{\sigma}=\left(\sigma^{1},\dots,\sigma^{k}\right) \in \mathcal{B}$,
\[
w.\boldsymbol{\sigma}:=\left( \Dot{w}_1\sigma^{1}\Dot{w}_1\inv,\dots,\Dot{w}_k \sigma^{k}\Dot{w}_{k}\inv\right),
\]
where $\Dot{w}_j$ is a representative in $\GL_n$ of $w_j\in W_{\GL_n}(L^j)$. Consider the diagram
\begin{equation}
\begin{tikzcd} \label{diagram_family_Mellit_qr}
  \mathcal{B}  \arrow[d,"\pi^0"'] & \arrow[l," \eta ",swap] \widetilde{\mathcal{Q}}_{\boldsymbol{L},\boldsymbol{P}} \arrow[d,"p"]   \\
  \mathcal{B} /W & \arrow[l,"\chi"]  \mathcal{Q}_{\boldsymbol{L},\boldsymbol{P}}.
\end{tikzcd}
\end{equation}
Thanks to the quiver variety point of view, the cohomology sheaves $\mathcal{H}^i \eta_! \kappa$ are constant (Theorem \ref{th_qlp_locally_constant}). In this section a $W$-equivariant structure on those cohomology sheaves is constructed. The Weyl group actions on the cohomology of quiver varieties with such constant sheaves were introduced by Nakajima \cite{nakajima1994}. Here we also use a method from Lusztig (see \cite[Proof of Proposition 5.5.3]{letellier_f}), this method is also applied by Laumon-Letellier \cite[Section 5.2]{laumon_letellier}. This approach allows to extend the equivariant structure away from a regular locus. Mellit obtained a similar result with a different construction for character varieties \cite[Section 8]{mellit_cell}.

Before constructing the equivariant structure, let us define the regular locus. Denote by ${\mathcal{B}}^{\reg}$ the subset of regular elements, i.e. elements $\left(\sigma^{1},\dots,\sigma^{k}\right)\in \mathcal{B}$ such that $Z_{\GL_n}(\sigma^{j})= L^j$. It is the locus of $\mathcal{B}$ where the $W$-action is free.
Diagram \eqref{diagram_family_Mellit_qr} is pulled back to the regular locus
\begin{equation}
\begin{tikzcd} \label{diagram_family_reg_Mellit_qr}
  {\mathcal{B}}^{\reg}  \arrow[d,"\pi^{\reg}"'] & \arrow[l," \eta^{\reg} ",swap] \widetilde{\mathcal{Q}}^{\reg}_{\boldsymbol{L},\boldsymbol{P}} \arrow[d,"p^{\reg}"]    \\
  {\mathcal{B}}^{\reg} /W & \arrow[l,"\chi^{\reg}"]  \mathcal{Q}^{\reg}_{\boldsymbol{L},\boldsymbol{P}}.
\end{tikzcd}
\end{equation}
Similarly to \eqref{eq_regular_fiber_prod}, notice that 
\begin{equation}\label{eq_isom_regular_prod}
    \mathcal{Q}^{\reg}_{\boldsymbol{L},\boldsymbol{P}}\times_{\mathcal{B}^{\reg}/W}  \mathcal{B}^{\reg}\cong \widetilde{\mathcal{Q}}^{\reg}_{\boldsymbol{L},\boldsymbol{P}}.
\end{equation}

\begin{theorem}\label{th_local_equiv}
The cohomology sheaves $\mathcal{H}^i\eta_! \kappa$ admit a $W$-equivariant structure over $\mathcal{B}$.
\end{theorem}
\begin{proof}
Consider the diagram
\begin{equation}
    \begin{tikzcd}[column sep = tiny, row sep = tiny]
    &  & \arrow[lld,bend right,"\eta"'] \arrow[dddl,bend left,"p"] \arrow[dl,"c"'] \widetilde{\mathcal{Q}}_{\boldsymbol{L},\boldsymbol{P}} \\
      \mathcal{B}  \arrow[dd,"\pi^0"'] & \mathcal{Q}_{\boldsymbol{L},\boldsymbol{P}}\times_{\mathcal{B}/W} \mathcal{B} \arrow[l,"a",swap]  \arrow[dd,"b"] &  \\
      &  &  \\
  \mathcal{B} /W & \arrow[l,"\chi"]  \mathcal{Q}_{\boldsymbol{L},\boldsymbol{P}} &
    \end{tikzcd},
\end{equation}
the group $W$ acts on $\mathcal{Q}_{\boldsymbol{L},\boldsymbol{P}}\times_{\mathcal{B}/W} \mathcal{B}$ and the morphism $a$ is $W$-equivariant.
The variety $\mathcal{Q}^{\reg}_{\boldsymbol{L},\boldsymbol{P}}\times_{\mathcal{B}^{\reg}/W}  \mathcal{B}^{\reg}$ is smooth, dense and open in $\mathcal{Q}_{\boldsymbol{L},\boldsymbol{P}}\times_{\mathcal{B}/W} \mathcal{B}$. The constant sheaf $\kappa$ over $\mathcal{Q}^{\reg}_{\boldsymbol{L},\boldsymbol{P}}\times_{\mathcal{B}^{\reg}/W}  \mathcal{B}^{\reg}$ is $W$-equivariant. Indeed for $w\in W$ we can define a morphism
\[
\phi_w : w^* \kappa \to \kappa
\]
which is the identity on the stalks. It satisfies the conditions of Definition \ref{def_W_equiv}. Applying the continuation principle from Remark \ref{remark_continuation_principle}, this $W$-equivariant structure extends to a $W$-equivariant structure on $\ic{\mathcal{Q}_{\boldsymbol{L},\boldsymbol{P}}\times_{\mathcal{B}/W} \mathcal{B}}$. %corollary 5.11 continuation principle Weil Conjectures, Perverse Sheaves and l’adic Fourier Transform
Notice that $\eta_! \kappa \cong a_! c_! \kappa$. We shall see in Lemma \ref{lemma_c_ic} that 
\[
c_!\kappa \cong \ic{\mathcal{Q}_{\boldsymbol{L},\boldsymbol{P}}\times_{\mathcal{B}/W} \mathcal{B}}.
\]
Then the $W$-equivariant structure on $c_! \kappa$ induces a $W$-equivariant structure on $\eta_! \kappa$.
% Restricted to the regular locus this $W$-equivariant structure coincide with the one from previous section as $c$ induces an isomorphism
% \begin{equation}\label{eq_isom_regular_prod}
% \mathcal{Q}^{\reg}_{\boldsymbol{L},\boldsymbol{P}}\times_{\mathcal{B}^{\reg}/W}  \mathcal{B}^{\reg} \cong \widetilde{\mathcal{Q}}^{\reg}_{\boldsymbol{L},\boldsymbol{P}}.
% \end{equation}
Up to the isomorphism $c_! \kappa \cong \ic{\mathcal{Q}_{\boldsymbol{L},\boldsymbol{P}}\times_{\mathcal{B}/W} \mathcal{B}}$, the theorem is proved.
\end{proof}

It remains to prove the following lemma.
\begin{lemma}\label{lemma_c_ic}
There is an isomorphism $c_! \kappa \cong \ic{\mathcal{Q}_{\boldsymbol{L},\boldsymbol{P}}\times_{\mathcal{B}/W} \mathcal{B}}$. 
%The map $c$ is small, so that $c_! \kappa\cong \ic{\mathcal{Q}_{\boldsymbol{L},\boldsymbol{P}}\times_{\mathcal{B}/W} \mathcal{B}}$.
\end{lemma}
\begin{proof}

Because of the isomorphism \eqref{eq_isom_regular_prod}, the restriction of $c_! \kappa$ to the smooth locus $\mathcal{Q}^{\reg}_{\boldsymbol{L},\boldsymbol{P}}\times_{\mathcal{B}^{\reg}/W} \mathcal{B}^{\reg}$ is the constant sheaf $\kappa$. In order to verify the hypothesis of Definition \ref{def_intersection_complex} it remains to prove that the map $c$ is small, i.e. that it satisfies the following inequality 
\[
\dim \left\lbrace x \in \mathcal{Q}_{\boldsymbol{L},\boldsymbol{P}}\times_{\mathcal{B}/W} \mathcal{B} \left| \dim c\inv (x) \ge d    \right.  \right\rbrace \le \dim \mathcal{Q}_{\boldsymbol{L},\boldsymbol{P}}\times_{\mathcal{B}/W} \mathcal{B} - 2 d \; \text{ for all }d>0.
\]
We use dimension estimates from Lusztig \cite[1.2]{Lusztig1984}, see also \cite[Theorem 1.4]{Shoji}. In the Lie algebra $\gl_n$ the estimate becomes, for $X$ in an adjoint orbit $\mathcal{O}$, 
\begin{equation}\label{eq_dimension_estimate}
\dim \left\lbrace g P \in \GL_n / P \left| g\inv X g \in \sigma + \mathfrak{u}_P \right. \right\rbrace \le \frac{1}{2} \left(n^2 - \dim L -\dim \mathcal{O}\right).
\end{equation}
The proof is then standard in Springer theory. Let $d>0$ and let $x$ be such that 
\[
\dim c\inv (x) \ge d,
\]
the element $x$ belongs to $\mathcal{Q}_{\mathcal{O}_{\boldsymbol{\rho},\boldsymbol{\sigma}}}$ for an element $\boldsymbol{\sigma}\in \mathcal{B}$ and for some adjoint orbits $\mathcal{O}_{\underline{\rho}^1,\sigma^1},\dots ,\mathcal{O}_{\underline{\rho}^k,\sigma^k}$. The dimension estimate \eqref{eq_dimension_estimate} implies
\[
d \le \frac{1}{2} \left( k n^2 -\sum_{j=1}^k \dim L^j-\dim \mathcal{O}_{\underline{\rho}^j,\sigma^j} \right),
\]
so that
\begin{equation*}\label{eq_estimake_kuple_adjoint}
\sum_{j=1}^k \dim \mathcal{O}_{\underline{\rho}^j,\sigma^j} \le  k n^2 - \sum_{j=1}^k \dim L^j -2 d.
\end{equation*}
Using the decomposition from Proposition \ref{prop_decomposition_qlp}, we have that $x\in \mathcal{Q}^{\mathcal{B},\boldsymbol{M},\boldsymbol{\rho}}_{\boldsymbol{L},\boldsymbol{P}}$. The previous inequality together with the expression \eqref{eq_dim_qlpmr} for the dimension of $\mathcal{Q}^{\mathcal{B},\boldsymbol{M},\boldsymbol{\rho}}_{\boldsymbol{L},\boldsymbol{P}}$ give
\begin{eqnarray}\label{eq_estimate_dim_qlpbmr}
\dim \mathcal{Q}^{\mathcal{B},\boldsymbol{M},\boldsymbol{\rho}}_{\boldsymbol{L},\boldsymbol{P}} \le n^2(2 g-2) + 2 + \dim Z(\boldsymbol{m})\cap \mathcal{B} + k n^2 - \sum_{j=1}^k \dim L^j-2 d.
\end{eqnarray}
Moreover
\begin{equation}\label{eq_prod_qlpbmr}
    \dim \mathcal{Q}^{\mathcal{B},\boldsymbol{M},\boldsymbol{\rho}}_{\boldsymbol{L},\boldsymbol{P}}\times_{\mathcal{B}/W}\mathcal{B} = \dim \mathcal{Q}^{\mathcal{B},\boldsymbol{M},\boldsymbol{\rho}}_{\boldsymbol{L},\boldsymbol{P}}
\end{equation}
and
\begin{equation}\label{eq_dim_prod_qlp}
\dim \mathcal{Q}_{\boldsymbol{L},\boldsymbol{P}}\times_{\mathcal{B}/W}\mathcal{B} = \dim \mathcal{Q}_{\boldsymbol{L},\boldsymbol{P}} =  n^2 (2 g -2) + 2 + \dim\mathcal{B} + k n^2 - \sum_{j=1}^k \dim L^j.
\end{equation}
Combining \eqref{eq_estimate_dim_qlpbmr},\eqref{eq_prod_qlpbmr}  and \eqref{eq_dim_prod_qlp},
\begin{equation}\label{eq_estimate_zb}
    \dim \mathcal{Q}^{\mathcal{B},\boldsymbol{M},\boldsymbol{\rho}}_{\boldsymbol{L},\boldsymbol{P}}\times_{\mathcal{B}/W}\mathcal{B}  \le \dim \mathcal{Q}_{\boldsymbol{L},\boldsymbol{P}}\times_{\mathcal{B}/W}\mathcal{B} + 2 d + \dim Z(\boldsymbol{m})\cap \mathcal{B}-\dim\mathcal{B}.
\end{equation}
As $d$ is assumed to be strictly positive, necessarily the inclusion $\boldsymbol{L}\subsetneq \boldsymbol{M}$ is strict, so that
\begin{equation}\label{eq_dim_zb}
\dim Z(\boldsymbol{m})\cap \mathcal{B}<\dim\mathcal{B}.
\end{equation}
Now \eqref{eq_estimate_zb} and \eqref{eq_dim_zb} provide the estimate
\begin{equation}\label{eq_final_estimate}
    \dim \mathcal{Q}^{\mathcal{B},\boldsymbol{M},\boldsymbol{\rho}}_{\boldsymbol{L},\boldsymbol{P}}\times_{\mathcal{B}/W}\mathcal{B} < \dim \mathcal{Q}_{\boldsymbol{L},\boldsymbol{P}}\times_{\mathcal{B}/W}\mathcal{B} - 2 d.
\end{equation}
% To conclude  notice that 
% \[
% \dim \mathcal{Q}^{\mathcal{B},\boldsymbol{M},\boldsymbol{\rho}}_{\boldsymbol{L},\boldsymbol{P}}\times_{\mathcal{B}/W}\mathcal{B} = \dim \mathcal{Q}^{\mathcal{B},\boldsymbol{M},\boldsymbol{\rho}}_{\boldsymbol{L},\boldsymbol{P}}.
% \]
To conclude, the set $\left\lbrace x\in \mathcal{Q}_{\boldsymbol{L},\boldsymbol{P}}\times_{\mathcal{B}/W}\mathcal{B} \left| \dim c\inv(x) \ge d  \right. \right\rbrace$ is a finite union of varieties $\mathcal{Q}^{\mathcal{B},\boldsymbol{M},\boldsymbol{\rho}}_{\boldsymbol{L},\boldsymbol{P}}\times_{\mathcal{B}/W}\mathcal{B}$ with dimensions satisfying the previous estimate \eqref{eq_final_estimate}.
\end{proof}

\begin{remark}
Let us study the restriction of the $W$-equivariant sheaves $\mathcal{H}^i\eta_! \kappa$ to the regular locus. Recall that 
$\mathcal{Q}^{\reg}_{\boldsymbol{L},\boldsymbol{P}}\times_{\mathcal{B}^{\reg}/W}  \mathcal{B}^{\reg}\cong \widetilde{\mathcal{Q}}^{\reg}_{\boldsymbol{L},\boldsymbol{P}}$, then for $\boldsymbol{\sigma}\in\mathcal{B}^{\reg}$
\[
\mathcal{H}^i_{\boldsymbol{\sigma}}\eta_!\kappa \cong H_c^i(\qlps,\kappa).
\]
For $w\in W$, the $W$-equivariant structure is given by the functoriality of the compactly supported cohomology (see Proposition \ref{prop_equivariant_push} and Remark \ref{remark_basechange_fiber})
\[
w^* : H_c^i\left(\widetilde{\mathcal{Q}}_{\boldsymbol{L},\boldsymbol{P},w.\boldsymbol{\sigma}},\kappa \right) \to   H_c^i\left(\qlps,\kappa\right).
\]
Therefore the construction of Theorem \ref{th_local_equiv} gives a canonical extension over $\mathcal{B}$ of this natural Weyl group action over $\mathcal{B}^{\reg}$. 
\end{remark}

\subsection{Monodromic Weyl group action on the cohomology \texorpdfstring{of $\qlps$}{}}

We saw in \ref{th_qlp_locally_constant} that the cohomology sheaves $\mathcal{H}^i \eta_! \kappa$ are constant sheaves over $\mathcal{B}$. Together with the $W$-equivariant structure, this allows to construct a Weyl group action on the cohomology of the varieties $\qlps$ for any $\boldsymbol{\sigma}\in\mathcal{B}$, this is called the monodromic Weyl group action. Note that the fiber over $\boldsymbol{\sigma}$ of this constant sheaf is $H_c^i( \qlps ; \kappa)$. Thus for any $\boldsymbol{\sigma},\boldsymbol{\boldsymbol{\tau}}\in \mathcal{B}$, there is an isomorphism
\[
f_{\boldsymbol{\sigma},\boldsymbol{\boldsymbol{\tau}}} :  H_c^i( \qlps ; \kappa) \to H_c^i( \widetilde{\mathcal{Q}}_{\mathbf{L},\mathbf{P}, \boldsymbol{\boldsymbol{\tau}}}  ; \kappa),
\]
such that for any $\boldsymbol{\omega}\in \mathcal{B}$
\[
f_{\boldsymbol{\sigma},\boldsymbol{\boldsymbol{\tau}}}=f_{\boldsymbol{\omega},\boldsymbol{\boldsymbol{\tau}}} \circ f_{\boldsymbol{\sigma},\boldsymbol{\omega}}.
\]
The $W$-equivariance of the local system $\mathcal{H}^i \eta_! \kappa$ implies the following theorem. It can also be proved directly, without referring to equivariance of the local system (see Maffei \cite[Section 5]{Maffei}).
\begin{theorem}\label{th_equivariantce_ls}
Let $\boldsymbol{\sigma},\boldsymbol{\boldsymbol{\tau}}\in \mathcal{B}$, then the following diagram commutes
\begin{equation*}
\begin{tikzcd}
H_c^i( \qlps ; \kappa) \arrow[r," w^* "]  \arrow[d, swap, "f_{\boldsymbol{\sigma},\boldsymbol{\boldsymbol{\tau}}}" ]&  H_c^i (\widetilde{\mathcal{Q}}_{\mathbf{L},\mathbf{P},w\inv .\boldsymbol{\sigma} } ; \kappa)  \arrow[d,  " f_{w\inv . \boldsymbol{\sigma}  , w\inv . \boldsymbol{\boldsymbol{\tau}} }" ] \\
 H_c^i (\widetilde{\mathcal{Q}}_{\mathbf{L},\mathbf{P}, \boldsymbol{\boldsymbol{\tau}} } ; \kappa) \arrow[r," w^* " ]  & H_c^i (\widetilde{\mathcal{Q}}_{\mathbf{L},\mathbf{P},w\inv.\boldsymbol{\boldsymbol{\tau}} } ; \kappa).
\end{tikzcd}
\end{equation*} 
\end{theorem}
\begin{remark}
Note that if $\boldsymbol{\sigma}\in\mathcal{B}$ is not regular, then the map 
\[
w^* : H_c^i\left(\qlps,\kappa\right)\to H_c^i\left(\widetilde{\mathcal{Q}}_{\boldsymbol{L},\boldsymbol{P},w\inv.\boldsymbol{\sigma}},\kappa\right)
\]
is only the map coming from the $W$-equivariant structure of the constant sheaf $\mathcal{H}^i\eta_!\kappa$. It does not necessarily come by functoriality from a morphism of variety.
\end{remark}
This theorem allows to define a $W$-action on the compactly supported cohomology space  $H_c^i( \qlps ; \kappa)$.
\begin{theorem}\label{prop_w_action_quiver}
For $\boldsymbol{\sigma}\in \mathcal{B}$ and for $w\in W$ let
\[
\rho_i(w) = f_{w.\boldsymbol{\sigma} , \boldsymbol{\sigma}} \circ (w\inv)^*.
\] 
This defines an action of $W$ on $H_c^i( \qlps ; \kappa)$, it is called the monodromic Weyl group action.
\end{theorem}
\begin{proof}
For $w_1$ and $w_2$ in $W$, the following diagram commutes by Theorem \ref{th_equivariantce_ls}.
\begin{equation*}
\begin{tikzcd}
H_c^i( \qlps ; \kappa) \arrow[r,"{(w_2\inv)}^*"] & H_c^i( \widetilde{\mathcal{Q}}_{\mathbf{L},\mathbf{P},w_2 .\boldsymbol{\sigma}  } ; \kappa) \arrow[r,"{(w_1 \inv )}^*"]  \arrow[d,"f_{w_2 . \boldsymbol{\sigma} ,\boldsymbol{\sigma}}"]& H_c^i( \widetilde{\mathcal{Q}}_{\mathbf{L},\mathbf{P},w_1 w_2 .\boldsymbol{\sigma}  } ; \kappa) \arrow[d,"f_{w_1 w_2 .\boldsymbol{\sigma} , w_1 . \boldsymbol{\sigma} }"] \\
 & H_c^i( \qlps ; \kappa) \arrow[r,"{(w_1 \inv)}^*"]  & H_c^i( \widetilde{\mathcal{Q}}_{\mathbf{L},\mathbf{P},w_1 . \boldsymbol{\sigma} } ; \kappa) \arrow[d,"f_{w_1 .\boldsymbol{\sigma} , \boldsymbol{\sigma} }"] \\
 & & H_c^i( \qlps ; \kappa)
 \end{tikzcd}
\end{equation*} 
Going from the top left corner to bottom right corner by the top right corner is $\rho(w_1 w_2)$. Going by the middle gives $\rho(w_1)\rho(w_2)$. Therefore $\rho(w_1 w_2) = \rho(w_1) \rho(w_2)$.

\end{proof}

\subsection{Frobenius morphism and monodromic action}
The techniques in this section come from Hausel, Letellier and Rodriguez-Villegas \cite{kac_conj}, however we do no consider regular semisimple values of the moment map. Instead each component of the moment map is central and each leg of the comet-shaped quiver corresponds to a particular adjoint orbit. Comet-shaped quiver varieties were also studied in this context by Letellier \cite{Letellier_uni_12}. A slightly more general situation is considered here, as a leg can represents any adjoint orbit and not only a semisimple regular one.

The representation defined in Theorem \ref{prop_w_action_quiver} when $\KK=\C$ is isomorphic to the representation obtained for $\KK=\Fqbar$ and large enough characteristic. Indeed this can be proved by base change exactly like in \cite[Theorem 2.5]{kac_conj}. 
%does By base change, when the characteristic of $\KK=\Fqbar$ is large enough, as in \cite{kac_conj}, we obtain an isomorphic representation $\rho_{\KK }^i$ on the compactly supported cohomology of the variety $\qlps$ over $\KK$. 
Therefore from now on we assume:
\begin{assumption}
$\mathbb{K}=\Fqbar$ and the characteristic is large enough.
\end{assumption}
This assumption is very convenient as it allows to introduce Frobenius endomorphism and use Grothendiek's trace formula to compute the traces of the action obtained.

We denote by $F$ the Frobenius endomorphism on $\gl_n $ raising the coefficients to the power $q$. The set of $F$-fixed points in $\gl_n$ is $\gl_n(\Fq)$ and similarly for the group $\GL_n$. Assume that the $L^j$ are subgroups of bock diagonal matrices, and that the $P^j$ are subgroups of block upper triangular matrices, then they are $F$-stable. The morphism $F$ induces a Frobenius endomorphism on $\qlp^{\reg}$ and on $\mathcal{B}^{reg}$ also denoted by $F$,
\[
F\left(\boldsymbol{\sigma}, \left(A_i,B_i\right)_{1\le i \le g} , (X_j, g_j L_j)_{1\le j\le k} \right) = \left(F(\boldsymbol{\sigma}), \left(F(A_i),F(B_i)\right)_{1\le i \le g} , (F(X_j), F(g_j) L_j)_{1\le j\le k} \right). 
\]
This Frobenius endomorphism can be twisted by an element $w=(w_1,\dots,w_k)$ in the Weyl group $W$. For $\boldsymbol{\sigma}\in B^{\reg}$, define 
\[
w F (\boldsymbol{\sigma}_1,\dots,\boldsymbol{\sigma}_k) := \left(w_1. F(\boldsymbol{\sigma}_1) , \dots, w_k . F(\boldsymbol{\sigma}_k) \right).
\]
 The set of points fixed by $w F$ is $\left(\mathcal{B}^{\reg}\right)^{w F}$.
Similarly, the $w$-twisted Frobenius on $\qlp^{\reg}$ is
\[
w F := w\circ F.
\]
They are compatible, $p^{\reg}\circ w F = w F \circ p^{\reg}$ so that for $\boldsymbol{\sigma},\boldsymbol{\boldsymbol{\tau}}\in \mathcal{B}^{\reg}$ the following diagram commutes
\[
\begin{tikzcd}
H_c^i( \widetilde{\mathcal{Q}}_{\mathbf{L},\mathbf{P}, \boldsymbol{\sigma}} ; \kappa) \arrow[r,"F^*"] \arrow[d,"f_{\boldsymbol{\sigma},\boldsymbol{\boldsymbol{\tau}}}"] & H_c^i( \widetilde{\mathcal{Q}}_{\mathbf{L},\mathbf{P}, F\inv(\boldsymbol{\sigma}) } ; \kappa)  \arrow[d,"f_{F\inv(\boldsymbol{\sigma}),F\inv(\boldsymbol{\boldsymbol{\tau}})}"]   \\
H_c^i( \widetilde{\mathcal{Q}}_{\mathbf{L},\mathbf{P}, \boldsymbol{\tau}} ; \kappa)  \arrow[r,"F^*"] &  H_c^i( \widetilde{\mathcal{Q}}_{\mathbf{L},\mathbf{P}, F\inv (\boldsymbol{\tau}) } ; \kappa) .
\end{tikzcd}
\] 
\begin{theorem}\label{th_card_twist}
For $\boldsymbol{\tau} \in \left( B^{\reg} \right)^{ F}$ and for $\boldsymbol{\sigma} \in \left( B^{\reg} \right)^{w F} $, the cardinal of the set of $w F$ fixed points in $\qlps$ is
\[
\sharp \qlps^{w F} =  \sum_i \tr \left(  \rho^{2 i}(w ) , H_c^{2 i} (\widetilde{\mathcal{Q}}_{\mathbf{L},\mathbf{P}, \boldsymbol{\tau}} ; \kappa )  \right) q^i.
\]
\end{theorem}
\begin{proof}
%Let $\boldsymbol{\tau} \in \left( B^{\reg} \right)^{ F}$ and $\boldsymbol{\sigma} \in \left( B^{\reg} \right)^{w F} $ and 
Consider the commutative diagram
\[
\begin{tikzcd}
H_c^i( \widetilde{\mathcal{Q}}_{\mathbf{L},\mathbf{P}, \boldsymbol{\tau}} ; \kappa) \arrow[r,"w^*"] \arrow[rd, swap, dashed, "\rho(w\inv)"]&  H_c^i( \widetilde{\mathcal{Q}}_{\mathbf{L},\mathbf{P}, w\inv. \boldsymbol{\tau}  } ; \kappa) \arrow[d,"f_{w\inv .\boldsymbol{\tau} , \boldsymbol{\tau}}"] & \\
 & H_c^i( \widetilde{\mathcal{Q}}_{\mathbf{L},\mathbf{P}, \boldsymbol{\tau}} ; \kappa) \arrow[r,"F^*"] &  H_c^i( \widetilde{\mathcal{Q}}_{\mathbf{L},\mathbf{P}, \boldsymbol{\tau}} ; \kappa)  \\
 H_c^i( \widetilde{\mathcal{Q}}_{\mathbf{L},\mathbf{P}, \boldsymbol{\sigma}} ; \kappa) \arrow[uu,"f_{\boldsymbol{\sigma},\boldsymbol{\tau}}"] \arrow[r,"w^*"] &  H_c^i( \widetilde{\mathcal{Q}}_{\mathbf{L},\mathbf{P}, F(\boldsymbol{\sigma}) } ; \kappa)  \arrow[u,"f_{F(\boldsymbol{\sigma}),\boldsymbol{\tau}}"] \arrow[r,"F^*"] & H_c^i( \widetilde{\mathcal{Q}}_{\mathbf{L},\mathbf{P}, \boldsymbol{\sigma}} ; \kappa) \arrow[u,"f_{\boldsymbol{\sigma},\boldsymbol{\tau}}"].
\end{tikzcd}
\]
Apply Grothendieck trace formula to $w F$,
\begin{eqnarray*}
\sharp \qlps^{w F} &=& \sum_i (-1)^i \tr \left( (w F)^*, H_c^i (\qlps ; \kappa )\right) \\
&=& \sum_i (-1)^i \tr \left( F^* \circ \rho^i(w\inv) , H_c^i (\widetilde{\mathcal{Q}}_{\mathbf{L},\mathbf{P}, \boldsymbol{\tau}} ; \kappa ) \right) .
\end{eqnarray*}
The varieties  $\mathcal{Q}_{\mathbf{L},\mathbf{P}, \boldsymbol{\tau}} $ are pure and polynomial count (see Remark \ref{remark_qss_resol} and \cite[Theorem 1.3.1]{aha}) and $\rho(w\inv)$ commutes with $F$ so that
\begin{eqnarray*}
\sharp \qlps^{w F}  &=& \sum_i \tr \left( F^* \circ \rho^{2 i}(w\inv) , H_c^{2 i} (\widetilde{\mathcal{Q}}_{\mathbf{L},\mathbf{P}, \boldsymbol{\tau}} ; \kappa )  \right) \\
&=&  \sum_i \tr \left(  \rho^{2 i}(w\inv) , H_c^{2 i} (\widetilde{\mathcal{Q}}_{\mathbf{L},\mathbf{P}, \boldsymbol{\tau}} ; \kappa )  \right) q^i.
\end{eqnarray*}
Now as $W$ is isomorphic to a product of symmetric groups, $w$ is conjugated to its inverse $w\inv$ and
\[
\sharp \qlps^{w F} =  \sum_i \tr \left(  \rho^{2 i}(w ) , H_c^{2 i} (\widetilde{\mathcal{Q}}_{\mathbf{L},\mathbf{P}, \boldsymbol{\tau}} ; \kappa )  \right) q^i.
\]
\end{proof}
\begin{notations}\label{remark_type_weyl}
The $j$-th part of $\boldsymbol{L}$ is
\[
L^j \cong \underbrace{ \GL_{c_{1}^j}\times \dots \times \GL_{c_1^j}}_{m_{1}^j} \times \dots \times \underbrace{ \GL_{c_{l_j}^j}\times \dots \times \GL_{c_{l_j}^j}}_{m_{l_j}^j} ,
\]
with $c^j_r \ne c^j_s$ for $r\ne s$. Then the $j$-th part of the relative Weyl group $W$ is
\[
W_{\GL_n} (L^j) \cong \Sym_{m_{1}^j}\times \dots \times \Sym_{m_{l_j}^j}
\]
The symmetric group $\Sym_{m_r^j}$ acts by permuting the blocks of size $c_r^j$. Take $w=(w_1,\dots,w_k)$ in $W$  and choose $\boldsymbol{\sigma}_w=\left(\sigma_1,\dots,\sigma_k\right)$ in $\left(\mathcal{B}^{\reg}\right)^{w F}$. The conjugacy class of the element $w_j$ is determined by a $l_j$-tuple $( \eta^{j,1},\dots,\eta^{j , l_j})$ with $\eta^{j,r}\in \Part_{m_{r}^j}$.  Let $\mathcal{O}_{\sigma_j}$ be the adjoint orbit of $\sigma_j$. This orbit is semisimple, $F$-stable and of the following type (as defined in \ref{def_type_fq_cc_ao}),
\[
\left( \eta^{j,1}_{1},1^{c_{1}^j}  \right)\dots \left( \eta^{j,1}_{l(\eta^{j,1})} , 1^{c_{1}^j} \right) \dots \left( \eta^{j, l_j}_{1}  , 1^{c_{l_j}^j} \right) \dots  \left( \eta^{j, l_j}_{  l(\eta^{j, l_j })  }, 1^{c_{l_j}^j} \right).
\]
Define $\mathcal{O}_{w}:=\left(\mathcal{O}_{\sigma_1},\dots,\mathcal{O}_{\sigma_k}\right)$.
\end{notations}
\begin{lemma}
With the previous notations we have the following identity between cardinals,
\begin{equation}\label{eq_equality_card}
\sharp \qlps^{w F}  = \sharp  \mathcal{Q}_{  \mathbf{{\mathcal{O}} }_{w } }^F.
\end{equation}
\end{lemma}
\begin{proof}
As the orbits $\mathcal{O}_{w}=\left(\mathcal{O}_{\sigma_1},\dots,\mathcal{O}_{\sigma_k}\right)$ are semisimple (hence they are closed), the map $\qlps \to \mathcal{Q}_{ {{\mathcal{O}}}_{w}}$ is an isomorphism compatible with the Frobenius $w F$ on the source and the Frobenius $F$ on the target.
\end{proof}
Letellier \cite{letellier_gl} computed the number of points of comet shaped quiver varieties, in particular of $\mathcal{Q}_{  {{\mathcal{O}} }_{w } }^F$.
\begin{theorem}
With Notations \ref{remark_type_weyl}, the cardinal of $\mathcal{Q}_{  {{\mathcal{O}} }_{w } }^F$ is given by
\[
\sharp \mathcal{Q}_{  \mathcal{O}_{w } }^F = (-1)^{r(\boldsymbol{\eta})}q^{\frac{d_{\boldsymbol{\mu}}}{2}} \left\langle \widetilde{h}_{\boldsymbol{\eta}}, \mathbb{H}_n^{H L V}(0, q^{\frac{1}{2}}) \right\rangle,
\]
where $\widetilde{h}_{\boldsymbol{\eta}}$ is a particular case of the generalized Schur function from \cite{letellier_gl}. This symmetric function can be expressed in terms of complete symmetric functions $h_n$,
\[
\widetilde{h}_{\boldsymbol{\eta}}:=\prod_{j=1}^k\prod_{r=1}^{l_j}\prod_{s=1}^{l(\eta^{j,i})} h_{c^j_r}\left[ X_j^{\eta^{j,r}_s}\right], 
\]
and
\[
r(\boldsymbol{\eta}) := \sum_{j=1}^k\sum_{r=1}^{l_j} c^j_r \sum_{s=1}^{l(\eta^{j,i})}(\eta^{j,r}_s-1).
\]
\end{theorem}
\begin{proof}
As the orbits $\mathcal{O}_{\sigma_j}$ are semisimple, the variety $\mathcal{Q}_{  {\mathcal{O} }_{w } }
$ is smooth so that the characteristic function of the intersection complex is constant with value $1$. The result follows from Letellier \cite[Theorem 6.9.1, Theorem 7.4.1 and Corollary 7.4.3]{letellier_gl}.
\end{proof}

\begin{corollary}\label{corollary_twisted_pc}
For $\boldsymbol{\sigma}\in\mathcal{B}$ and $\boldsymbol{\eta}$ representing a conjugacy class in the Weyl group as described in Notations \ref{remark_type_weyl}, the $\boldsymbol{\eta}$-twisted Poincaré polynomial of $\widetilde{\mathcal{Q}}_{\boldsymbol{L},\boldsymbol{P},\boldsymbol{\sigma}}$ is
\[
\sum_{i}\tr\left(\boldsymbol{\eta},H_c^i(\qlps,\kappa) \right) v^{i} = (-1)^{r(\boldsymbol{\eta})}v^{d_{\boldsymbol{\mu}}} \left\langle \widetilde{h}_{\boldsymbol{\eta}}, \mathbb{H}_n^{H L V}(0, v) \right\rangle.
\]
\end{corollary}
\begin{proof}
The action comes from the $W$-equivariant structure of the constant sheaves $\mathcal{H}^i\eta_! \kappa$. Therefore, up to isomorphism, the representation does not depend on the choice of $\boldsymbol{\sigma}\in\mathcal{B}$ so that the twisted Poincaré polynomial can be computed for $\boldsymbol{\tau}\in(\mathcal{B}^{\reg})^F$. Then from Theorem \ref{th_card_twist} and from \eqref{eq_equality_card},
\[
\sum_i \tr \left(\rho^{2 i}(\boldsymbol{\eta}), H_c^{2 i}\left(\widetilde{\mathcal{Q}}_{\boldsymbol{L},\boldsymbol{P},\boldsymbol{\tau}},\kappa\right)\right)q^i = (-1)^{r(\boldsymbol{\eta})}q^{\frac{d_{\boldsymbol{\mu}}}{2}} \left\langle \widetilde{h}_{\boldsymbol{\eta}}, \mathbb{H}_n^{H L V}(0, q^{\frac{1}{2}}) \right\rangle.
\]
This equality remains true after substituting $q^n$ for $q$ for $n>0$. Thus it is an equality between two polynomials and the corollary is proved.
\end{proof}

\begin{remark}[Comparison between monodromic and Springer action]\label{remark_comparison_action}
Let $\boldsymbol{\sigma}=\left(\sigma_1,\dots,\sigma_k\right)\in\mathcal{B}$, as before $M^j$ is the stabilizer of $\sigma_j$ in $\GL_n$. The relative Weyl group is
\[
W_{\boldsymbol{M}}(\boldsymbol{L}) = \prod_{j=1}^k W_{M^j}(L^j)
\]
with $W_{M^j}(L^j) = N_{M^j}(L^j)/L^j$. Then $W_{\boldsymbol{M}}\left(\boldsymbol{L}\right)$ is a subgroup of the Weyl group $W$ studied in this section. The group $W_{\boldsymbol{M}}\left(\boldsymbol{L}\right)$ is exactly the subgroup of elements $w\in W$  such that $w.\boldsymbol{\sigma}=\boldsymbol{\sigma}$. The monodromic Weyl group action from Theorem \ref{prop_w_action_quiver} induces an action of $W_{\boldsymbol{M}}\left(\boldsymbol{L}\right)$ on $H_c^i\left(\qlps, \kappa \right)$. Interestingly, this action comes only from the $W$-equivariant structure, it does not rely on the constant property of the sheaf: it is given explicitly by $\rho_i(w)=\left(w\inv\right)^*$.

There is another action of $W_{\boldsymbol{M}}\left(\boldsymbol{L}\right)$ on $H_c^i\left(\qlps, \kappa \right)$, the Springer action constructed by Letellier and mentioned in \ref{remark_springer_comet}. Letellier computed the twisted Poincaré polynomial for this Springer action \cite[Corollary 7.4.3]{letellier_gl}, it coincides with the Poincaré polynomial obtained from the monodromic action, therefore both action are isomorphic. It would be interesting to have a direct proof of this fact. We proved it in the character variety setting for just one orbit, regular, with a unique eigenvalue \cite[Chapter 5]{ballandras_thesis}.

\end{remark}

It is also interesting to consider the monodromic action over the regular locus $\mathcal{B}^{\reg}$ as an action on the cohomology of a quiver variety with semisimple adjoint orbits. For $\boldsymbol{\sigma}\in \mathcal{B}^{\reg}$ consider the associated generic $k$-tuple of semisimple adjoint orbits $\mathcal{O}_{\boldsymbol{\sigma}}=\left(\mathcal{O}_{\sigma_1},\dots,\mathcal{O}_{\sigma_k}\right)$. The Weyl group $W_{\GL_n}(L^j)$ is the group of permutation of the distinct eigenvalues of $\mathcal{O}_{\sigma_j}$ with the same multiplicities. This provides another formulation of Corollary \ref{corollary_twisted_pc}.

\begin{corollary}\label{corollary_twisted_pc_ss}
For $\boldsymbol{\eta}$ representing a conjugacy class in the Weyl group as described in Notations \ref{remark_type_weyl}, the $\boldsymbol{\eta}$-twisted Poincaré polynomial of $\mathcal{Q}_{\mathcal{O}_{\boldsymbol{\sigma}}}$ is
\[
\sum_{i}\tr\left(\boldsymbol{\eta},H_c^i(\mathcal{Q}_{\mathcal{O}_{\boldsymbol{\sigma}}},\kappa) \right) v^{i} = (-1)^{r(\boldsymbol{\eta})}v^{d_{\boldsymbol{\eta}}} \left\langle \widetilde{h}_{\boldsymbol{\eta'}}, \mathbb{H}_n^{H L V}(0, v) \right\rangle.
\]
\end{corollary}

The interpretation over the regular locus in terms of semisimple quiver varieties together with Remark \ref{remark_comparison_action} show the advantages of extending the $W$-equivariant structure from $\mathcal{B}^{\reg}$ to $\mathcal{B}$ (Theorem \ref{th_equivariantce_ls}). This provides a uniform description of the Springer action on the cohomology of some resolution $\qlps$ and the monodromic action on the cohomology of semisimple quiver varieties $\mathcal{Q}_{\mathcal{O}_{\boldsymbol{\sigma}}}$.

\begin{remark}\label{remark_isotypical_puncture}
It is also interesting to study the action of a Weyl group relative to a particuar leg $1\le j\le k$, for instance relative to the first one. This will be used in \ref{subsect_coefficient_coho} to describe some structure coefficients of the algebra spanned by Kostka polynomial. A particularly interesting case is when  $L^1$ is a maximal torus and $M^1=\GL_n$. Then the component of the Weyl group relative to the first leg is $W_{M^1}(L^1)\cong \Sym_n$ and
\[
W_{\boldsymbol{M}}(\boldsymbol{L})\cong \Sym_n\times\prod_{j=2}^k W_{M^j}(L^j).
\]
According to this decomposition, consider an element $(w,1,\dots,1)\in W_{\boldsymbol{M}}(\boldsymbol{L})$ with $w\in\Sym_n$ an element of cycle type $\lambda\in\Part_n$. Then
\[
\widetilde{h}_{\boldsymbol{\eta}}= p_{\lambda}[X_1] h_{\underline{\mu'}^2}[X_2]\dots h_{\underline{\mu'}^k}[X_k],
\]
and $(-1)^{r(\boldsymbol{\eta})}=\epsilon(\lambda)$ is the sign of the permutation $w$ with cycle type $\lambda$. Corollary \ref{corollary_twisted_pc} reads
\[
 P_c^{\boldsymbol{\eta}}\left(\qlps,v\right) = v^{d_{\boldsymbol{\mu}}}\epsilon(\lambda)\left\langle p_{\lambda}[X_1] h_{\underline{\mu'}^2}[X_2]\dots h_{\underline{\mu'}^k}[X_k],\mathbb{H}_n^{H L V}\left(0,v \right)  \right\rangle.
\]
This can be understood in terms of Frobenius characteristic, see Definition \ref{def_frobenius_char}. Consider the representation of $\Sym_n$ on the cohomology of $\qlps$ twisted by the sign, $H^{\bullet}(\qlps,\kappa)\otimes \epsilon$. Its $v$-graded Frobenius characteristic is given by the following symmetric function in $X_1$
\[
v^{d_{\boldsymbol{\mu}}}\left\langle h_{\underline{\mu'}^2}[X_2]\dots h_{\underline{\mu'}^k}[X_k],\mathbb{H}_n^{H L V}\left(0,v \right)  \right\rangle_{X_2,\dots,X_k}.
\]
Notice that $V_{\rho}\otimes \epsilon \cong V_{\rho'}$, then by Remark \ref{remark_trace_mult}, the multiplicity of the irreducible representation $V_{\rho}$ in $H^{\bullet}(\qlps,\kappa)$ is given by
\[
v^{d_{\boldsymbol{\mu}}}\left\langle s_{\rho'}[X_1] h_{\underline{\mu'}^2}[X_2]\dots h_{\underline{\mu'}^k}[X_k],\mathbb{H}_n^{H L V}\left(0,v \right)  \right\rangle.
\]

\end{remark}

% \[
% \sharp  \mathcal{Q}_{  \mathbf{\mathcal{O} }_{\omega } }^F = (-1)^{r(\boldsymbol{\omega})} q^{d_{\boldsymbol{\omega}}} \left\langle s_{\boldsymbol{\omega'}}, \mathbb{H}_n^{H L V} \right\rangle 
% \]

\section{Geometric interpretations in the algebra spanned by Kostka polynomials}\label{sect_combi_kostka_quiver}
\subsection{Description of the algebra}
In this section an algebra spanned by Kostka polynomials is studied and some structure coefficients are related to traces of Weyl group action on the cohomology of quiver varieties.
Define a linear map $\Delta^{\#} : \SymF[X]\to \SymF[X,Y]$ such that on the basis of modified Macdonald polynomials,
\[
\Delta^{\#}\left(\Ht_{\lambda}[X]\right):=\Ht_{\lambda}[X] \tilde{H}_{\lambda}[Y] \; \text{ for  all }\lambda\in\Part.
\]
As in \ref{notation_qt_drop}, the variables $(q,t)$ are implicit. Now as the Hall pairing is non-degenerate, there is a uniquely determined bilinear map $\dots \# \dots$ such that for all $F,G$ and $H$ in $\SymF[X]$,
\[
\left\langle F[X] \# G[X], H[X]\right\rangle = \left\langle F[X]G[Y],\Delta^{\#}\left( H[X]\right) \right\rangle .
\]
The product $\#$ defines an associative and commutative algebra structure on $\SymF[X]$.
\begin{definition}\label{c_coeff}
For a $k$-tuple of partitions $\boldsymbol{\mu}=\left(\mu^{1},\dots,\mu^k\right)\in\Part_n^k$ and for $\lambda\in\Part_n$ we denote by $c_{\boldsymbol{\mu}}^{\lambda}$ the structure coefficients of the product $\#$ in the basis of Schur functions 
\begin{equation}
    s_{\mu^1}\#s_{\mu^2}\dots\#s_{\mu^k}=\sum_{|\mu|=n}c_{\boldsymbol{\mu}}^\lambda s_\lambda.
\end{equation}
\end{definition}

\begin{lemma}
For $\boldsymbol{\mu}=(\mu,\nu)$, the coefficients $c_{\mu,\nu}^{\lambda}$ coincide with those defined in the introduction, \textit{i.e.}, the following relation is satisfied
\begin{equation}\label{eq_kostka_c}
    \widetilde{K}_{\mu,\rho} \widetilde{K}_{\nu,\rho} = \sum_{\lambda} c_{\mu,\nu}^{\lambda}  \widetilde{K}_{\lambda,\rho}.
\end{equation}
\end{lemma}
\begin{proof}
First let $\left(\widetilde{L}_{\eta,\lambda}\right)_{\lambda,\eta\in\Part_n}$ be the inverse of the matrix of Kostka polynomials $\left(\widetilde{K}_{\eta,\lambda}\right)_{\lambda,\eta\in\Part_n}$ (see Definition \ref{def_kostka_poly}), then
\[
s_{\lambda} = \sum_{\eta\in\Part_n} \widetilde{L}_{\eta,\lambda} \tilde{H}_{\eta}[X].
\]
Now the coefficient $c_{\mu,\nu}^{\lambda}$ is defined by
\begin{eqnarray*}
c_{\mu,\nu}^{\lambda} &=& \left\langle s_{\mu} \# s_{\nu},s_{\lambda}\right\rangle \\
& =& \left\langle s_{\mu} \# s_{\nu},\sum_{\eta\in\Part_n} \widetilde{L}_{\eta,\lambda} \tilde{H}_{\eta}[X]\right\rangle.
\end{eqnarray*}
By definition of the product $\#$ and of the coproduct $\Delta^{\#}$,
\begin{eqnarray*}
c_{\mu,\nu}^{\lambda} &=& \sum_{\eta\in\Part_n} \widetilde{L}_{\eta,\lambda} \left\langle s_{\mu}[X] s_{\nu}[Y] , \tilde{H}_{\eta}[X]\tilde{H}_{\eta}[Y] \right\rangle,  \\
c_{\mu,\nu}^{\lambda} &=& \sum_{\eta\in\Part_n} \widetilde{L}_{\eta,\lambda} \widetilde{K}_{\mu,\eta} \widetilde{K}_{\nu,\eta}.
\end{eqnarray*}
Multiply the last equation by $\widetilde{K}_{\lambda,\rho}$ and sum over $\lambda\in\Part_n$,
\[
\widetilde{K}_{\mu,\rho} \widetilde{K}_{\nu,\rho} = \sum_{\lambda} c_{\mu,\nu}^{\lambda}  \widetilde{K}_{\lambda,\rho}.
\]
This last relation is exactly the one used in the introduction to define the coefficients $c_{\mu,\nu}^{\lambda}$.
\end{proof}
\begin{example}
We computed some coefficients with the software SageMath
% \begin{eqnarray*}
% c_{(2^2),(2,1^2)}^{(2,1^2)} &=& -q^3 t - q^2 t^2 - q t^3 -q^2 t -t^2 q + q^2 + q t + t^2 \\
% c_{(2^2),(2,1^2)}^{(1^4)} &=& q^3+ q^2 t + q t^2 + t^3 + q^2 + 2 q t + t^2 + q +t
% \end{eqnarray*}
\begin{eqnarray*}
c_{(2,2),(2,1,1)}^{(2,1,1)} &=& -q^3 t - q^2 t^2 - q t^3 -q^2 t -t^2 q + q^2 + q t + t^2, \\
c_{(2,2),(2,1,1)}^{(1, 1, 1, 1)} &=& q^3+ q^2 t + q t^2 + t^3 + q^2 + 2 q t + t^2 + q +t.
\end{eqnarray*}
\end{example}
The next conjecture comes from unpublished notes by Rodriguez-Villegas.
\begin{conjecture}\label{conjecture_1}
The structure coefficients $c_{\boldsymbol{\mu}}^{\lambda}$ lie in $\mathbb{Z}[q,t]$.
\end{conjecture}
Some evidences supporting this conjecture will be provided. The following definition and remark were suggested by François Bergeron.
\begin{definition}
Let $F$ be a symmetric function, consider the operator
\begin{eqnarray*}
\begin{array}{ccccc}
     F\#\dots &:&\SymX &\to & \SymX  \\
     & & G & \mapsto & F\#G.
\end{array}
\end{eqnarray*}
We denote $\psi_F$ its adjoint with respect to the Hall pairing so that for any $G,H\in\SymX$
\begin{equation}\label{adjonction}
    \left\langle F\#G,H\right\rangle=\left\langle G,\psi_F(H)\right\rangle 
\end{equation}
Those operators are diagonal in the basis of modified Macdonald polynomials
\begin{equation}\label{psi_f_diag}
    \psi_F(\Ht_\lambda[X;q,t])=\left\langle F,\Ht_\lambda[X;q,t]\right\rangle \Ht_\lambda[X;q,t].
\end{equation}
\end{definition}
\begin{remark}\label{remark_psi_en}
Applying the relation \eqref{psi_f_diag} with $e_n$, 
\begin{equation*}
    \psi_{e_n}\left(\Ht_\lambda[X;q,t]\right)=q^{n(\lambda')}t^{n(\lambda)}\Ht_\lambda[X;q,t],
\end{equation*}
and we recognize the usual expression of the operator $\nabla$ introduced by Bergeron-Garsia \cite{sf}. The higher $(q,t)$-Catalan sequence from Garsia--Haiman \cite{garsia} (see also Haiman \cite[p.95]{haiman2002}) is defined by
\[
C_{n}^{(m)}(q,t) = \left\langle e_n, \nabla^m e_n\right\rangle ,
\]
but $\nabla = \psi_{e_n}$ is the adjoint of $e_n  \# \dots$ and $s_{1^n}=e_n$ so that
\[
C_{n}^{(m)}(q,t) = c^{1^n}_{\underbrace{1^n,\dots,1^{n}}_{m+1} }.
\]
The higher $(q,t)$-Catalan sequences are particular cases of the coefficients $c_{\boldsymbol{\mu}}^{1^n}$.
\end{remark}
We recall an important theorem which was first conjectured by Garsia--Haiman \cite{garsia}.
\begin{theorem}[\cite{haiman2002} theorem 4.2.5]\label{th_nabla_en}
The symmetric function $\nabla(e_n)$ is obtained as the Frobenius characteristic (see definition  \ref{def_frobenius_char}) of a bigraded representation of $\Sym_n$ called the diagonal harmonics. In particular,
\begin{equation*}
    \left\langle\nabla(e_n),s_\lambda\right\rangle \in \mathbb{Z}_{\ge 0}[q,t].
\end{equation*}
\end{theorem}
\begin{corollary}
For any $\mu\in\Part_n$, the structure coefficients $c_{1^n,\mu}^{1^n}$ gives the multiplicity of the irreducible representation of type $\mu$ in the bigraded representation of $\Sym_n$ on diagonal harmonics. In particular  $c_{1^n,\mu}^{1^n}(q,t)\in\mathbb{Z}_{\ge 0}[q,t]$ so that the conjecture \ref{conjecture_1} is true for those particular coefficients.
\end{corollary}
\begin{proof}
According to Remark \ref{remark_psi_en} and to the adjonction relation \eqref{adjonction},
\begin{equation}\label{eq_proof_nabla}
    \left\langle s_\mu,\nabla(e_n)\right\rangle=\left\langle e_n\#s_\mu,e_n\right\rangle.
\end{equation}
By definition of the structure coefficients $c_{\mu,\nu}^{\lambda}$ and as $e_n=s_{1^n}$, we have
\begin{equation*}
    e_n\# s_\mu=\sum_{\lambda\in\Part_n}c_{1^n\mu}^{\lambda} s_{\lambda},
\end{equation*}
substituting in \eqref{eq_proof_nabla} we obtain
\begin{equation*}
    c_{1^n,\mu}^{1^n}(q,t)=\left\langle s_\mu,\nabla(e_n)\right\rangle.
\end{equation*}
We conclude with the interpretation of $\nabla(e_n)$ as a Frobenius characteristic from Theorem \ref{th_nabla_en}.
\end{proof}

The next theorem and the following corollary come from unpublished notes by Rodriguez-Villegas. They relate particular structure coefficients $c_{\boldsymbol{\mu}}^{1^n}$ to the kernel $\mathbb{H}_n^{H L V}$.

Consider the generating function from Definition \ref{def_omega_hlv} for genus $g=0$, $k+2$ punctures and with variable $z=q^{\frac{1}{2}}$, $w=t^{\frac{1}{2}}$. It is given by

\begin{equation*}
    \Omega^0_{k+2}:=\sum_{\lambda\in\Part} \frac{\prod_{i=1}^{k+2}\Ht_\lambda \left[X_i;q,t\right]}{a_\lambda(q,t)}s^{|\lambda|},
\end{equation*}
with $a_{\lambda}(q,t) = \left( \Htl,\Htl\right)^{q,t}$ as in \ref{prop_charact_mcd}.
\begin{theorem}\label{fernando_th}
The following relation holds,
\begin{equation*}
    \left\langle p_{(n)}[X_{k+1}]h_{(n-1,1)}[X_{k+2}],\Log\left[\Omega_{k+2}^g\right]\right\rangle_{X_{k+1},X_{k+2}}=\sum_{|\lambda|=n}\frac{\phi_\lambda \Pi_{\lambda}'}{a_{\lambda}}\prod_{i=1}^k\Ht_\lambda[X_i]s^{|\lambda|},
\end{equation*}
with
\begin{eqnarray*}
\phi_\lambda&=&\sum_{i,j\in\lambda}q^{j-1}t^{i-1},\\
\Pi'_{\lambda}&=&\prod_{i,j\in\lambda\setminus\left( 1,1\right)}(1-q^{j-1}t^{i-1}).
\end{eqnarray*}
\end{theorem}
\begin{proof}
According to Lemma \ref{h_pairing}, to take the Hall pairing with $h_{(n-1,1)}[X_{k+2}]$ is equivalent to do plethystic substitution $X_{k+2}=1+u$ and to take the degree $n$ coefficient in front of $u$. As the plethystic substitution and the plethystic logarithm commute, we can perform this substitution inside the plethystic logarithm.  We consider terms of order $1$ in $u$ using Lemma \ref{lemma_uplus},
\begin{eqnarray*}
    \Log\left[\Omega_{k+2}^0\right]=\Log\left[\Omega_{k+1}^0+u\sum_{\lambda\in\Part^*}\frac{\phi_\lambda}{a_{\lambda}}\prod_{i=1}^{k+1}\Ht_\lambda[X_i]s^{|\lambda|}+\mathcal{O}(u^2)\right]\\
    =\Log\left[\Omega_{k+1}^0\left(1+u\frac{1}{\Omega_{k+1}^0}\sum_{\lambda\in\Part^*}\frac{\phi_\lambda}{a_{\lambda}}\prod_{i=1}^{k+1}\Ht_\lambda[X_i]s^{|\lambda|}+\mathcal{O}(u^2)\right)\right] \\
    = \Log\left[\Omega_{k+1}^0\right]+\Log\left[1+u\frac{1}{\Omega_{k+1}^0}\sum_{\lambda\in\Part^*}\frac{\phi_\lambda}{a_{\lambda}}\prod_{i=1}^{k+1}\Ht_\lambda[X_i]s^{|\lambda|}+\mathcal{O}(u^2)\right].
\end{eqnarray*}
We used that plethystic logarithm turns products into sums. From the definition of the plethystic logarithm, as $p_n[u]=u^n$, we easily see the coefficient in front of $u$ in the previous expression
\begin{eqnarray*}
      \left.\Log\left[\Omega_{k+2}^0\right]\right|_{u}
      =\frac{1}{\Omega_{k+1}^0}\sum_{\lambda\in\Part^*}\frac{\phi_\lambda}{a_{\lambda}}\prod_{i=1}^{k+1}\Ht_\lambda[X_i]s^{|\lambda|}.
\end{eqnarray*}
Keeping the terms of degree $n$ we obtain
\begin{equation*}
    \left\langle h_{(n-1,1)}[X_{k+2}],\Log\left[\Omega_{k+2}^0\right]\right\rangle_{X_{k+2}}=\left.\frac{1}{\Omega_{k+1}^0}\sum_{\lambda\in\Part^*}\frac{\phi_\lambda}{a_{\lambda}}\prod_{i=1}^{k+1}\Ht_\lambda[X_i]s^{|\lambda|}\right|_{s^n}.
\end{equation*}
Inverting $\Omega_{k+1}^0$ is licit, it is defined by
%in terms of power series since $\Omega_{k+1}^0-1$ lies in $s\SymF[X_1,\dots,X_{k+1}][[s]]$
\begin{equation*}
    \frac{1}{\Omega_{k+1}^0}=\frac{1}{1+\left(\Omega^0_{k+1}-1\right)}=\sum_k \left(1-\Omega_{k+1}^0\right)^k.
\end{equation*}
Now we just have to take the Hall pairing with the power sum $p_{(n)}\left[X_{k+1}\right]$. This is equivalent to picking the coefficient in front of $n^{-1}p_{(n)}\left[X_{k+1}\right]$. But $p_{(n)}$ cannot be written as the product of two symmetric functions of degree strictly smaller than $n$ so that the contribution of $\Omega_{k+1}^0$ in the denominator is irrelevant for the coefficient in front of $n^{-1}p_{(n)}\left[X_{k+1}\right]$ and
\begin{equation*}
    \left\langle p_{(n)}[X_{k+1}]h_{(n-1,1)}[X_{k+2}],\Log\left[\Omega_{k+2}^0\right]\right\rangle_{X_{k+1},X_{k+2}}=\left\langle p_{(n)}[X_{k+1}], \sum_{\lambda\in\Part^*}\frac{\phi_\lambda}{a_{\lambda}}\prod_{i=1}^{k+1}\Ht_\lambda[X_i]s^{|\lambda|}\right\rangle_{X_{k+1}} .
\end{equation*}
We conclude with Lemma \ref{pn_scal} and \eqref{pleth_u}.
\end{proof}
The following corollary allows to obtain a geometric interpretation of the coefficients. Indeed, it relates the coefficient $c_{\boldsymbol{\mu}}^{(1^n)}$ to the generating serie $\Omega_{k+2}^0$ known to encode cohomological information about comet-shaped quiver varieties and character varieties.
\begin{corollary}\label{corollary_c_mu}
With the notations of the previous theorem and Definition \ref{c_coeff},
\begin{equation}\label{c_Log}
(-1)^{n-1}c_{\boldsymbol{\mu}}^{(1^n)}=(q-1)(1-t)\left\langle\prod_{j=1}^k s_{\mu^j}[X_j]p_{(n)}[X_{k+1}]h_{(n-1,1)}[X_{k+2}],\Log\left[\Omega_{k+2}^0\right]\right\rangle_{X_{1},\dots,X_{k+2}}.
%    (-1)^{n-1}c_{\boldsymbol{\mu}}^{(1^n)}=(q-1)(1-t)\left(s_{\mu^{1}}[X_1]\dots s_{\mu^{k}}[X_k]p_{(n)}[X_{k+1}]h_{(n-1,1)}[X_{k+2}],\Log\left[\Omega_{k+2}^0\right]\right)_{X_{1},\dots,X_{k+2}}
\end{equation}
\end{corollary}
\begin{proof}
We apply Theorem \ref{fernando_th} to express the right hand side of \eqref{c_Log} as
\begin{equation*}
    (q-1)(1-t)\left\langle s_{\mu^{1}}[X_1]\dots  s_{\mu^{k}}[X_k],\sum_{|\lambda|=n}\frac{\phi_\lambda \Pi'_{\lambda}}{a_{\lambda}}\prod_{i=1}^k\Ht_\lambda[X_i]\right\rangle_{X_{1},\dots,X_{k}}.
\end{equation*}
By definition of the product $\#$,
\begin{equation*}
    (q-1)(1-t)\left\langle s_{\mu^{1}}\#\dots\#s_{\mu^{k}}[X],\sum_{|\lambda|=n}\frac{\phi_\lambda\Pi'}{a_{\lambda}}\Ht_\lambda[X]\right\rangle_{X}.
\end{equation*}
Here we recognize the expression of Theorem \ref{garsia_m_th} 
\begin{equation*}
    \left\langle s_{\mu^{1}}\#\dots\#s_{\mu^{k}}[X],(-1)^{n-1}s_{(1^n)}\right\rangle_{X}
\end{equation*}
so that if we write
\begin{equation*}
    s_{\mu^{1}}\#\dots\#s_{\mu^{k}}[X]=\sum_{\lambda}c_{\boldsymbol{\mu}}^\lambda s_\lambda[X]
\end{equation*}
the result follows from orthonormality of Schur functions.
\end{proof}

\subsection{Interpretation of certain coefficients as traces of Weyl group actions on the intersection cohomology of quiver varieties}\label{subsect_coefficient_coho}

In this section a cohomological interpretation is given for the coefficients $c_{\boldsymbol{\mu}}^{1^n}$. In order to lighten the notations, the description is only given for the coefficient $c_{\mu,\nu}^{1^n}$. The generalization to any $\boldsymbol{\mu}$ is straightforward.

First let us detail the data to describe the relevant variety $\qlps$. The Levi subgroups are torus of diagonal matrices $L^j=T$ for $1\le j \le 3$ and $L^4=\GL_1\times \GL_{n-1}$. The semisimple part $\boldsymbol{\sigma}=\left(\sigma^1,\dots,\sigma^4\right)$ is such that:
\begin{itemize}
    \item $\sigma^1=\zeta_1 \Id$ is central,
    \item $\sigma^2=\zeta_2 \Id$ is central,
    \item $\sigma^{3}=\begin{pmatrix}
    \alpha_1 & &  &\\
    & \alpha_2 & & \\
    & & \ddots & \\
     & & & \alpha_n
  \end{pmatrix}$ with $\alpha_r\ne \alpha_s$ for $r\ne s$,
    \item $\sigma^{4}=\begin{pmatrix}
    \beta & &  &\\
    & \gamma & & \\
    & & \ddots & \\
     & & & \gamma
  \end{pmatrix}$ has two eigenvalues $\beta\ne \gamma$. The multiplicity of $\beta$ is one and the multiplicity of $\gamma$ is $n-1$.
\end{itemize}
Notice that such a choice can be made in the generic locus, \textit{i.e.}, with $\boldsymbol{\sigma}\in \mathcal{B}$.

First we consider Letellier's construction of the action \textit{à la} Springer in order to compute isotypical components. Let $\boldsymbol{M}=M^1\times\dots\times M^{4}$ with $M^j$ the stabilizer in $\GL_n$ of $\sigma^j$. Then $W_{\boldsymbol{M}}(\boldsymbol{L})\cong \Sym_n^2$. Letellier's construction (Remark \ref{remark_springer_comet}) provides an action of $W_{\boldsymbol{M}}(\boldsymbol{L})$ on the cohomology of $\qlps$. Moreover for $V_{\mu'}$, respectively $V_{\nu'}$, the irreducible representation of $\Sym_n$ associated to the transpose of some partition $\mu$, respectively $\nu$,
\begin{equation}\label{eq_mult_munu}
\Hom_{W_{\boldsymbol{M}}(\boldsymbol{L})}\left(V_{\mu'}\otimes V_{\nu'}, H_c^{i+d_{\qlps}}\left(\qlps,\kappa\right) \right) = H_c^{i+d_{\mathcal{Q}_{\overline{\boldsymbol{\mathcal{O}}}}}}\left(\mathcal{Q}_{\overline{\boldsymbol{\mathcal{O}}}},\kappa\right).
\end{equation}
With $\boldsymbol{\mathcal{O}}=(\mathcal{O}_1,\dots,\mathcal{O}_4)$ the $4$-tuple of generic adjoint orbits defined by,
\begin{itemize}
    \item $\mathcal{O}_1$ has Jordan type $\mu'$ and eigenvalue $\zeta_1$,
    \item $\mathcal{O}_2$ has Jordan type $\nu'$ and eigenvalue $\zeta_2$,
    \item $\mathcal{O}_3$ is the orbit of $\sigma^3$,
    \item $\mathcal{O}_4$ is the orbit of $\sigma^4$.
\end{itemize}

Now with the construction from Theorem \ref{prop_w_action_quiver}, there is an action of the whole group $W\cong \Sym_n^3$ on the cohomology of $\qlps$. The restriction of this $W$-action to $W_{\boldsymbol{M}}(\boldsymbol{L})\cong \Sym_n^2$ is isomorphic to the Springer action (see Remark \ref{remark_comparison_action}). First take the $V_{\mu'}\otimes V_{\nu'}$ isotypical component with respect to the $ \Sym_n^2$-action. There remains an action of the Weyl group $\Sym_n$ relative to the third leg on the intersection cohomology $I H_c^{i}\left(\mathcal{Q}_{\overline{\boldsymbol{\mathcal{O}}}},\kappa\right)$. 
\begin{theorem}\label{th_quiver_c}
Let $w$ be a $n$-cycle in the Weyl group relative to the third leg (this terminology comes from the comet-shaped quiver). The coefficient $c_{\mu,\nu}^{1^n}$, after specialization $q=0$, is given by the $w$-twisted Poincaré polynomial of $\mathcal{Q}_{\overline{\boldsymbol{\mathcal{O}}}}$, namely
\[
c_{\mu,\nu}^{1^n}(0,t) = t^{-\frac{d_{\overline{\boldsymbol{\mathcal{O}}}}}{2}} \sum_{i} \tr\left(w,I H_c^{2 i} \left(\mathcal{Q}_{\overline{\boldsymbol{\mathcal{O}}}},\kappa\right)\right) t^{i}.
\]
\end{theorem}
\begin{proof}
Combining \eqref{eq_mult_munu}, Corollary \ref{corollary_twisted_pc} and Remark \ref{remark_isotypical_puncture},
\[
\sum_{i} \tr\left(w,I H_c^i \left(\mathcal{Q}_{\overline{\boldsymbol{\mathcal{O}}}},\kappa\right)\right) v^i= (-1)^{n-1} v^{d_{\overline{\mathcal{O}}}} \left\langle s_{\mu}[X_1]s_{\nu}[X_2] p_{(n)}[X_3] h_{(n-1,1)}[X_4], \mathbb{H}_n^{H L V}(0,v) \right\rangle.
\]
The theorem now follows from Corollary \ref{corollary_c_mu}.
\end{proof}

\subsection{Cohomological interpretation in the multiplicative case}\label{subsect_conj_mixed_c}
There are similar interpretations in the multiplicative case. A conjectural one involving $c_{\mu,\nu}^{1^n}(q,t)$ which is a theorem after specializing to $c_{\mu,\nu}^{1^n}(1,t)$. Unfortunately in the multiplicative case the monodromic action is not defined in the general case so that we have to rely only on the Springer action. Therefore the statements involve partial resolutions of character varieties instead of actual character varieties.

First introduce the relevant parameters. The Levi subgroups are torus of diagonal matrices $L^j=T$ for $1\le j \le 3$ and $L^4=\GL_1\times \GL_{n-1}$. The semisimple part $\boldsymbol{\sigma}=(\sigma^1,\dots,\sigma^4)$ is such that:
\begin{itemize}
    \item $\sigma^1=\zeta_1 \Id$ is central,
    \item $\sigma^2=\zeta_2 \Id$ is central,
    \item $\sigma^3=\zeta_3 \Id$ is central,
    \item $\sigma^{4}=\begin{pmatrix}
    \beta & &  &\\
    & \gamma & & \\
    & & \ddots & \\
     & & & \gamma
  \end{pmatrix}$ has two eigenvalues $\beta\ne \gamma$. %The multiplicity of $\beta$ is one and the multiplicity of $\alpha$ is $n-1$.
\end{itemize}
This $4$-tuple is chosen to be generic (in the multiplicative sense of Definition \ref{def_generic_cc}). This is the case for instance if $\zeta_1 \zeta_2 \zeta_3 = 1$ and $\gamma^{n-1}=\beta\inv \ne 1$.
%Moreover this $4$-tuple can be chosen in the regular locus. Note that the notations for the parameters are the same as in the previous section, however objects are different as we now consider resolutions of character varieties. For instance the eigenvalues are now necessarily non zero and the genericity condition is the multiplicative one.
The relative Weyl group is $W_{\boldsymbol{M}}(\boldsymbol{L})\cong \Sym_n^3$. Now consider the following conjugacy classes
\begin{itemize}
    \item $\mathcal{C}_1$ has Jordan type $\mu '$ and eigenvalue $\zeta_1$,
    \item $\mathcal{C}_2$ has Jordan type $\nu '$ and eigenvalue $\zeta_2$,
    \item $\mathcal{C}_3$ has one Jordan block of size $n$ with eigenvalue $\zeta_3$,
    \item $\mathcal{C}_4$ is the conjugacy class of $\sigma_4$.
\end{itemize}
Then $\mlps$ is the resolution of $\mathcal{M}_{\overline{\boldsymbol{\mathcal{C}}}}$ with $\boldsymbol{\mathcal{C}}=\left(\mathcal{C}_1,\dots,\mathcal{C}_4\right)$ (see Definition \ref{def_char_var}). An intermediate in between $\mlps$ and $\mathcal{M}_{\overline{\boldsymbol{\mathcal{C}}}}$ is given by the variety
\begin{multline*}
  \mathcal{M}_{\mu,\nu}=\left\lbrace (X_1,\dots, X_4)\in \overline{\mathcal{C}}_1\times\dots \times \overline{\mathcal{C}}_4, g B \in \GL_n /B \left| g\inv X_3 g \in \zeta_3 U \right.\right.\\ 
  \left. X_1 \dots X_4= \Id \right\rbrace   \git \PGL_n,
\end{multline*}
with $B$ the Borel subgroup of upper triangular matrices in $\GL_n$ and $U$ its unipotent radical. Then the resolution $\mlps\to \mathcal{M}_{\overline{\boldsymbol{\mathcal{C}}}}$ factors through $\mathcal{M}_{\mu,\nu}$. This is a particular case of the \textit{partial} resolutions of character varieties studied by Letellier \cite{zclosure}. The result we recalled about Springer theory for resolutions of character varieties \eqref{eq_springer_char} admit a more general version for partial resolutions.  In particular considering the action of $\Sym_n^2$ with respect to the first two punctures and taking the $V_{\mu'}\otimes V_{\nu '}$ isotypical component of the cohomology $H_c^{\bullet}\left(\mlps,\kappa\right)$ we obtain,
\[
\Hom_{W_{\boldsymbol{M}}(\boldsymbol{L})}\left(V_{\mu'}\otimes V_{\nu'}, H_c^{i+d_{\mlps}}\left(\mlps,\kappa\right) \right) = H_c^{i+d_{\mathcal{M}_{\mu,\nu}}}\left(\mathcal{M}_{\mu,\nu},\kappa\right).
\]
There remains an action of the Weyl group $\Sym_n$ relative to the third puncture. For $w$ in this Weyl group define the $w$-twisted mixed Hodge polynomial by
\[
IH_c^{w}\left( \mathcal{M}_{\mu,\nu}; u, v \right):=   \sum_{i,r} u^r v^i \tr\left(w, IH^{r, r,i}_c\left(\mathcal{M}_{\mu,\nu},\kappa \right)\right).
\]
Conjecture \ref{conjecture_ic_char} admits a generalization describing the Weyl group action on the intersection cohomology of partial resolutions of character varieties (Letellier \cite[Conjecture 5.5]{zclosure}). In particular this conjecture predicts the following formula for the $w$-twisted mixed Hodge polynomial for $w$ a $n$-cycle,
\begin{multline*}
 IH_c^{w}\left( \mathcal{M}_{\mu,\nu}; u, v \right) = \\ (-1)^{n-1}\left(v\sqrt{u}\right)^{\dim \mathcal{M}_{\mu,\nu}}\left\langle s_{\mu}[X_1] s_{\nu} [X_2] p_{(n)}[X_3] h_{(n-1,1)}[X_4]  , \mathbb{H}_n^{H L V}\left(\frac{-1}{\sqrt{u}}, v \sqrt{u}\right)\right\rangle  .
\end{multline*}

The next conjecture follows from this conjectural formula, just like Theorem \ref{th_quiver_c} is deduced from Corollary \ref{corollary_twisted_pc}. 

\begin{conjecture}\label{conjecture_c}
Let $w$ be a $n$-cycle in the Weyl group relative to the third puncture. The coefficient $c_{\mu,\nu}^{1^n}$ relates to the $w$-twisted mixed Hodge polynomial of $\mathcal{M}_{\mu',\nu'}$ by
\[
c_{\mu,\nu}^{1^n}(q,t) = t^{\frac{{-\dim {\mathcal{M}_{\mu,\nu}}}}{2}} IH_c^{w}\left( \mathcal{M}_{\mu,\nu}, \frac{1}{q}, \sqrt{q t} \right).
\]
\end{conjecture}

We will prove that the right handside of this conjecture is indeed a polynomial in $q,t$ thus supporting Conjecture \ref{conjecture_1}. 

\begin{lemma}\label{lemma_poincare}
The Poincaré polynomial (for compactly supported intersection cohomology) of the variety $\mathcal{M}_{\mu,\nu}$ is
\begin{equation}    \label{eq_poincare_mmunu}
\sum_i v^i\dim IH_c^i\left(\mathcal{M}_{\mu,\nu},\kappa\right) = v^d \left\langle s_{\mu}[X_1]s_{\nu}[X_2]h_{1^n}[X_3]h_{(n-1,1)}[X_4] , \mathbb{H}_n^{H L V}\left(-1, v\right)\right\rangle ,
\end{equation}
and for $w$ a $n$-cycle in the Weyl group relative to the third puncture the $w$-twisted Poincaré polynomial is
\[
\sum_i v^i\dim \tr\left(w,IH_c^i\left(\mathcal{M}_{\mu,\nu},\kappa\right) \right)= v^d \left\langle s_{\mu}[X_1]s_{\nu}[X_2]p_{(n)}[X_3]h_{(n-1,1)}[X_4] , \mathbb{H}_n^{H L V}\left(-1, v\right)\right\rangle .
\]
\end{lemma}
\begin{proof}
Letellier \cite[Theorem 5.4]{zclosure} studied Springer theory for partial resolutions of character varieties. This allows to describe the intersection cohomology of the partial resolutions, together with its Weyl group action, in terms of intersection cohomology of character varieties. The Poincaré polynomials for intersection cohomology of character varieties are computed in \cite{ballandras_ic}. The formula for (twisted) Poincaré polynomials of the resolutions are therefore a consequence of the Poincaré polynomial specialization of \cite[Proposition 5.7]{zclosure}.
\end{proof}

\begin{proposition}
The following expression
\[
t^{\frac{{-\dim {\mathcal{M}_{\mu,\nu}}}}{2}} IH_c^{w}\left( \mathcal{M}_{\mu,\nu}, \frac{1}{q}, \sqrt{q t} \right),
\]
which is the value of $c_{\mu,\nu}^{1^n}(q,t)$ according to Conjecture \ref{conjecture_c}, is a polynomial in $q,t$ with integer coefficients.
\end{proposition}
\begin{proof}
First note that only integer powers of $q$ and $t$ appear because of \eqref{eq_poincare_mmunu} and the fact that in genus $g=0$ the kernel $\mathbb{H}_n^{HLV}(z,w)$ contains only terms in $z^2, w^2$.

Let $\mathcal{S}$ be a regular semisimple conjugacy class such that the $4$-tuple $\boldsymbol{\mathcal{{C}'}}=\left(\mathcal{C}_1,\mathcal{C}_2,\mathcal{S},\mathcal{C}_4\right)$ is generic. By \cite{ballandras_ic} and \eqref{eq_poincare_mmunu}, the Poincaré polynomial of the character variety $\mathcal{M}_{\boldsymbol{\mathcal{\overline{C}'}}}$ is the same as the Poincaré polynomial of the variety $\mathcal{M}_{\mu,\nu}$. But the variety $\mathcal{M}_{\boldsymbol{\mathcal{\overline{C}'}}}$ is affine, hence its compactly supported intersection cohomology vanishes in degree strictly smaller than its dimension and only positive power of $t$ appear in the expression. Only positive power of $q$ appear because the weight on the compactly supported intersection cohomology is smaller than the cohomological degree.
\end{proof}

The second equality in Lemma \ref{lemma_poincare} implies the Poincaré polynomial specialisation of Conjecture \ref{conjecture_c}.
\begin{theorem}
Let $w$ be a $n$-cycle in the Weyl group relative to the third puncture. The coefficient $c_{\mu,\nu}^{1^n}$ relates to the $w$-twisted Poincaré polynomial of $\mathcal{M}_{\mu,\nu}$:
\[
c_{\mu,\nu}^{1^n}(1,t) = t^{\frac{{-\dim {\mathcal{M}_{\mu,\nu}}}}{2}} \sum_i t^{\frac{i}{2}} \tr \left( w , \IH_c^i\left(\mathcal{M}_{\mu,\nu},\kappa\right)\right).
\]
\end{theorem}

\bibliographystyle{alpha}
\bibliography{biblio_these}

\end{document}